\documentclass{article}[10pt]
\usepackage{a4}
\usepackage[T1]{fontenc}
\usepackage[latin1]{inputenc}
\usepackage[francais]{babel}
\usepackage{hyperref}
\usepackage[dvips]{graphics}
\usepackage{amsmath}
\usepackage{amssymb}
\usepackage{amsopn}
\usepackage{amsbsy}
\usepackage{amstext}
\usepackage{latexsym}
\usepackage{amsfonts}
\usepackage{array}
\usepackage{epsfig}
\usepackage[mathscr]{eucal}

 \newtheorem{theorem}{Th\'eor\`eme}[section]
 \newtheorem{corollary}[theorem]{Corollaire}
 \newtheorem{lemma}[theorem]{Lemme}
 \newtheorem{proposition}[theorem]{Proposition}
 \newtheorem{definition}[theorem]{D\'efinition}
 \newtheorem{definitions}[theorem]{D\'efinitions}

 \newtheorem{thm*}{Th\'eor\`eme}
 \numberwithin{equation}{section}
\def\EMdash{\leavevmode\hbox to 7.5mm{\vrule height .63ex depth -.59ex
    width 5.4mm\hfill}}
\long\def\parag#1{\vspace{0.2cm} \\ {\bf #1}}

\def\Mbar{{\, \overline{M}}}
\def\Gbar{{\, \overline{G}}}
\def\Kbar{{\, \overline{K}}}
\def\hbar{{\, \overline{h}}}
\def\Qbar{{\, \overline{\mathbb{Q}}}}
\def\vvec{{\, \vec{\rm{v}}}}
\def\deg{{\rm{deg}}}
\def\rang{{\rm{rang}}}
\def\Ass{{\rm{Ass}}}

\def\plongement{\hookrightarrow}
\def\Xsoul{{\underline{X}}}
\def\Ysoul{{\underline{Y}}}

\def\Asoul{{\underline{A}}}
\def\Tsoul{{\underline{T}}}
\def\Zsoul{{\underline{Z}}}
\def\esoul{{\underline{e}}}
\def\alphasoul{{\underline{\alpha}}}
\def\betasoul{{\underline{\beta}}}
\def\gammasoul{{\underline{\gamma}}}
\def\Xsoull{{\underline{{\underline{X}}}}}
\def\Tsoull{{\underline{{\underline{T}}}}}
\def\Asoull{{\underline{{\underline{A}}}}}
\def\alphasoull{\underline{{\underline{\alpha}}}}
\def\l{{\ell}}
\def\c{{\rm{c}}}

\def\EMts{\mspace{.3mu}}  
\def\nb#1{{\left\vert{\EMts\EMts #1 \EMts\EMts}\right\vert}}        
\def\norm#1{{\left\Vert{\EMts\EMts #1 \EMts\EMts}\right\Vert}}      

\def\:{\!\!:}

\def\x{{\mathbf{x}}}

\def\0{{\mathbf{0}}}
\def\longu{{\rm long}}
\let\epsilon=\varepsilon
\font\septpt = cmr7

\begin{document}

\title{Un lemme de Roth sur les groupes algébriques commutatifs}

\author{Bakir FARHI \\
\footnotesize{D\'epartement de Math\'ematiques, Universit\'e du
Maine,} \\ \footnotesize{Avenue Olivier Messiaen, 72085 Le Mans
Cedex 9, France.} \\ \footnotesize{Bakir.Farhi@univ-lemans.fr}}
\date{}

\maketitle \tableofcontents \newpage
\section{{\bf HAUTEURS DE POLYN\^OMES ET DE VARI\'ET\'ES}}
\subsection{Mesures absolues:}
Soit $P$ un polyn\^ome de $\mathbb{C}[X_1 , \dots , X_n]$ $(n \in
{\mathbb N}^*)$, on définit la mesure de Gauss-Weil de $P$, notée
$\widetilde{M}(P)$, comme le maximum des normes de tous les
coefficients de $P$ et on définit la mesure de Mahler de $P$,
notée $\Mbar(P)$, par la relation:
$$\log{\Mbar(P)} := \int_{0}^{1} \!\!\!\int_{0}^{1} \!\!\!\dots \!\!\!\int_{0}^{1} \!\!\!\log{\nb{P\left(e^{2 \pi i x_1} , \dots , e^{2 \pi i x_n}\right)} d x_1 \dots d x_n} .$$
Si maintenant $P$ est un polyn\^ome homog\`ene de $\mathbb{C}[X_0
, \dots , X_n]$ $(n \in {\mathbb N}^*)$, on définit -comme dans
\cite{Ph4}3- la mesure unitaire homog\`ene de $P$, notée $M(P)$,
par la relation:
$$\log{M(P)} := \int_{S_{n + 1}(1)} \!\!\!\log{\nb{P\left(X_0 , \dots , X_n\right)}} \sigma_{n + 1} + {d°}_{\rm{tot}}P . \sum_{j = 1}^{n} \frac{1}{2 j}$$
o\`u $S_{n + 1}(1)$ est la sph\`ere unité de ${\mathbb C}^{n + 1}$ et $\sigma_{n + 1}$ est la mesure invariante de masse totale $1$ sur $S_{n + 1}(1)$. \\
Finalement, si $P$ est un polyn\^ome multihomog\`ene en des
ensembles de variables ${\underline X}^{(1)} , \dots , {\underline
X}^{(p)}$ (avec ${\underline X}^{(i)} = (X_{0}^{(i)} , \dots ,
X_{m_i}^{(i)})$ pour $i = 1 , \dots , p$), on définit aussi (comme
dans \cite{Ph4}3) la mesure unitaire multihomog\`ene $M(P)$ de $P$
par:
\begin{equation*}
\begin{split}
\log{M(P)} &:= \int_{S_{m_1 + 1}(1) \times \dots \times S_{m_p + 1}(1)} \!\!\!\log{\nb{P\left({\underline X}^{(1)} , \dots , {\underline X}^{(p)}\right)}} \sigma_{m_1 + 1} \times \dots \times \sigma_{m_p + 1} \\
&\quad+ \sum_{i = 1}^{p} \left({d°}_{{\underline X}^{(i)}} P .
\sum_{j = 1}^{m_i} \frac{1}{2 j}\right) .
\end{split}
\end{equation*}
La comparaison entre les deux mesures $\widetilde{M}$ et $\Mbar$
est donnée par le lemme suivant:
\begin{lemma}\label{a.7}
Soient $n$ un entier $\geq 1$ et $P$ un polyn\^ome de
$\mathbb{C}[X_1 , \dots , X_n]$ (ou un polyn\^ome homog\`ene de
$\mathbb{C}[X_0 , \dots , X_n]$) de degré total $d$. Alors on a:
\begin{equation}
{\binom{d + n}{n}}^{\!\!- 1/2} \Mbar(P) \leq \widetilde{M}(P) \leq
2^{n d} \Mbar(P) . \label{1.7}
\end{equation}
\end{lemma}
{\bf Démonstration.---} Démontrons le lemme \ref{a.7} dans le cas
o\`u $P$ est un polyn\^ome quelconque de $\mathbb{C}[X_1 , \dots ,
X_n]$ (la démonstration est la m\^eme pour le cas homog\`ene).
Commençons par démontrer l'inégalité de droite de (\ref{1.7}) en
raisonnant par récurrence sur $n \in {\mathbb N}^*$.

Pour $n = 1$, $P$ est un polyn\^ome d'une variable de degré $d$,
donc, en appelant $\alpha_1 , \dots , \alpha_d$ ses racines dans
$\mathbb C$, $P$ se factorise en:
$$
P(X_1) = a_0 X_{1}^{d} + a_1 X_{1}^{d - 1} + \dots + a_d = a_0
\prod_{i = 1}^{d} \left(X_1 - \alpha_i\right) .
$$
D'une part, la mesure de Gauss-Weil de $P$ vaut par définition:
$$\widetilde{M}(P) := \max\left\{\nb{a_0} , \dots , \nb{a_d}\right\}$$
et d'autre part, d'apr\`es \cite{Mah} (pp. 5-6), la mesure de
Mahler de $P$ vaut:
$$\Mbar(P) := \nb{a_0} . \prod_{i = 1}^{d} \max\left\{1 , \nb{\alpha_i}\right\} .$$
La formule -bien connue- exprimant les $a_j$ ($j = 0 , \dots , d$)
en fonction de $a_0$ et des $\alpha_1 , \dots , \alpha_d$ est
donnée par:
$$a_j = (- 1)^j a_0 . \sum_{1 \leq s_1 < \dots < s_j \leq d} \!\!\!\alpha_{s_1} \dots \alpha_{s_j} ~~~~ (\forall j \in \{0 , \dots , d\}) ;$$
donc, pour tout $j \in \{0 , \dots , d\}$ on a:
\begin{equation*}
\begin{split}
\nb{a_j} &\leq \nb{a_0} \binom{d}{j} \prod_{i = 1}^{d} \max\left\{1 , \nb{\alpha_i}\right\} \\
&\leq 2^d \nb{a_0} \prod_{i = 1}^{d} \max\left\{1 ,
\nb{\alpha_i}\right\} = 2^d \widetilde{M}(P) ,
\end{split}
\end{equation*}
d'o\`u:
$$\widetilde{M}(P) := \max_{0 \leq j \leq d} \nb{a_j} \leq 2^d \Mbar(P)$$
ce qui est l'inégalité de droite de (\ref{1.7}) pour $n = 1$. Soit
maintenant $n$ un entier $\geq 2$. Supposons que l'inégalité de
droite de (\ref{1.7}) est vraie pour l'entier $n - 1$ et montrons
la pour l'entier $n$. Soit, pour cela, $P$ un polyn\^ome de
$\mathbb{C}[X_1 , \dots , X_n]$ de degré total $d$ et écrivons $P$
sous la forme:
$$P\left(X_1 , \dots , X_n\right) = \sum_{j \in J} P_j\left(X_1 , \dots , X_{n - 1}\right) X_{n}^{j}$$
pour un certain sous-ensemble fini $J$ de $\mathbb N$ et certains
polyn\^omes $P_j$ ($j \in J$) de $\mathbb{C}[X_1 , \dots , X_{n -
1}]$. En fait, comme $P$ est de degré total égal à $d$, on a $J
\subset \{0 , \dots , d\}$ et les polyn\^omes $P_j$ ($j \in J$) de
$\mathbb{C}[X_1 , \dots , X_{n - 1}]$ sont tous de degrés totaux
majorés par $d$. On a clairement:
$$\widetilde{M}(P) = \max_{j \in J} \widetilde{M}(P_j)$$
et, d'apr\`es l'hypoth\`ese de récurrence, pour tout $j \in J$:
$$\widetilde{M}(P_j) \leq 2^{(n - 1) d} \Mbar(P_j) ;$$
d'o\`u
\begin{equation}
\widetilde{M}(P) \leq 2^{(n - 1) d} \max_{j \in J} \Mbar(P_j).
\label{1.8}
\end{equation}
On a par ailleurs:
\begin{equation*}
\begin{split}
\log\max_{j \in J} \Mbar(P_j) &= \max_{j \in J} \log\Mbar(P_j) \\
&= \max_{j \in J} \int_{0}^{1} \!\!\!\dots \!\!\!\int_{0}^{1} \!\!\!\log{\nb{P_j\left(e^{2 \pi i x_1} , \dots , e^{2 \pi i x_{n - 1}}\right)}} d x_1 \dots d x_{n - 1} \\
&\leq \int_{0}^{1} \!\!\! \dots \!\!\!\int_{0}^{1} \!\!\!
\log\left(\max_{j \in J} \nb{P_j\left(e^{2 \pi i x_1} , \dots ,
e^{2 \pi i x_{n - 1}}\right)}\right) d x_1 \dots d x_{n - 1}
\end{split}
\end{equation*}
et, d'apr\`es le cas $n = 1$ (déja démontré):
\begin{equation*}
\begin{split}
\max_{j \in J} \nb{P_j \left(e^{2 \pi i x_1} , \dots , e^{2 \pi i x_{n - 1}}\right)} &= \widetilde{M}\!\left(P\left(e^{2 \pi i x_1} , \dots , e^{2 \pi i x_{n - 1}} , X_n\right)\right) \\
&\leq 2^d \Mbar\!\left(P\left(e^{2 \pi i x_1} , \dots , e^{2 \pi i
x_{n - 1}} , X_n\right)\right)
\end{split}
\end{equation*}
o\`u -dans ces deux inégalités- $x_1 , \dots , x_{n - 1}$ sont
considérés comme param\`etres et $X_n$ considéré comme l'unique
variable. D'o\`u:
$$\log\max_{j \in J} \Mbar(P_j) \leq \int_{0}^{1} \!\!\!\dots \!\!\!\int_{0}^{1} \!\!\!\log\left(2^d \Mbar\!\left(P\left(e^{2 \pi i x_1} , \dots , e^{2 \pi i x_{n - 1}} , X_n\right)\right)\right) d x_1 \dots d x_{n - 1} ,$$
c'est-à-dire:
$$\log\max_{j \in J} \Mbar(P_j) \leq d \log{2} + \int_{0}^{1} \!\!\!\dots \!\!\!\int_{0}^{1} \!\!\!\log{\nb{P\left(e^{2 \pi i x_1} , \dots , e^{2 \pi i x_n}\right)}} d x_1 \dots d x_n$$
ou encore:
$$\log\max_{j \in J} \Mbar(P_j) \leq d \log{2} + \log\Mbar(P)$$
(car $\log\Mbar(P(e^{2 \pi i x_1} , \dots , e^{2 \pi i x_{n - 1}} , X_n)) := \int_{0}^{1} \log{\nb{P(e^{2 \pi i x_1} , \dots , e^{2 \pi i x_n})}} d x_n$). \\
D'o\`u, en passant à l'exponentielle:
\begin{equation}
\max_{j \in J} \Mbar(P_j) \leq 2^d \Mbar(P) . \label{1.9}
\end{equation}
De (\ref{1.8}) et (\ref{1.9}) suit l'inégalité:
$$\widetilde{M}(P) \leq 2^{n d} \Mbar(P)$$
qui n'est rien d'autre que l'inégalité de droite de (\ref{1.7})
pour l'entier $n$. L'inégalité de droite de (\ref{1.7}) est
démontrée.

Montrons maintenant l'inégalité de gauche de (\ref{1.7}). On utilise pour cela le théor\`eme de
l'inégalité arithmético-géometrique dont l'énoncé est: \\
\og Pour toute fonction réguli\`ere \footnote{La régularité ici
consiste en l'existence des intégrales constituants notre
inégalité arithmético-géométrique.} $f$ de ${\mathbb C}^n$ on a:
\begin{equation*}
\begin{split}
\exp\!\left(\int_{0}^{1} \!\!\!\dots \!\!\!\int_{0}^{1} \!\!\!\log{\nb{f\left(e^{2 \pi i x_1} , \dots , e^{2 \pi i x_n}\right)}} d x_1 \dots d x_n\right) \\
&\!\!\!\!\!\!\!\!\!\!\!\!\!\!\!\!\!\!\!\!\!\!\!\!\!\!\!\!\!\!\!\!\!\!\!\!\!\!\leq
\int_{0}^{1} \!\!\!\dots \!\!\!\int_{0}^{1} \!\!\!\nb{f\left(e^{2
\pi i x_1} , \dots , e^{2 \pi i x_n}\right)} d x_1 \dots d x_n .
\fg
\end{split}
\end{equation*}
En appliquant cette derni\`ere à $f^2$ (au lieu de $f$) on obtient
l'inégalité:
\begin{equation}
\begin{split}
\exp\!\left(\int_{0}^{1} \!\!\!\dots \!\!\!\int_{0}^{1} \!\!\!\log{\nb{f\left(e^{2 \pi i x_1} , \dots , e^{2 \pi i x_n}\right)}} d x_1 \dots d x_n\right) \\
&\!\!\!\!\!\!\!\!\!\!\!\!\!\!\!\!\!\!\!\!\!\!\!\!\!\!\!\!\!\!\!\!\!\!\!\!\!\!\!\!\!\!\!\!\!\!\!\!\leq
\left(\int_{0}^{1} \!\!\!\dots \!\!\!\int_{0}^{1} \!\!\!
{\nb{f\left(e^{2 \pi i x_1} , \dots , e^{2 \pi i x_n}\right)}}^2 d
x_1 \dots d x_n\right)^{\!\!\frac{1}{2}} . \label{1.10}
\end{split}
\end{equation}
\'Etant donné maintenant un polyn\^ome $P$ de $\mathbb{C}[X_1 ,
\dots , X_n]$ de degré total $d$, appliquant (\ref{1.10}) à $f =
P$ on a:
$$\Mbar(P) \leq \left(\int_{0}^{1} \!\!\! \dots \!\!\!\int_{0}^{1} \!\!\! {\nb{P\left(e^{2 \pi i x_1} , \dots , e^{2 \pi i x_n}\right)}}^2 d x_1 \dots d x_n\right)^{\!\!\frac{1}{2}} = {\norm{P}}_2$$
o\`u ${\norm{P}}_2$ désigne la norme euclidienne du vecteur formé
de tous les coefficients de $P$. Comme $P$ est à $n$ variables de
degré total $d$, $P$ contient au plus $\binom{d + n}{n}$
mon\^omes, d'o\`u:
$${\norm{P}}_2 \leq {\binom{d + n}{n}}^{1/2} \widetilde{M}(P) ,$$
ce qui donne finalement:
$$\Mbar(P) \leq {\norm{P}}_2 \leq {\binom{d + n}{n}}^{1/2} \widetilde{M}(P)$$
qui n'est rien d'autre que l'inégalité de gauche de (\ref{1.7}).
La démonstration est achevée. $~~~~\blacksquare$\vspace{1mm}
\subsection{Mesures $v$-adiques et hauteurs de polyn\^omes:}
\subsubsection{{\bf Mesures $v$-adiques:}}
Soient $K$ un corps de nombres et $v$ une place de $K$. On définit
${\widetilde{M}}_v$ et $\Mbar_v$ les mesures $v$-adiques des
polyn\^omes à coefficients dans $K$, associant à tout $P \in K[X_1
, \dots , X_n]$ $(n \in {\mathbb N}^*)$ les réels positifs notés
respectivement ${\widetilde{M}}_v (P)$ et $\Mbar_v (P)$ et
définis, selon le cas $v$ finie ou $v$ infinie, par:
\begin{description}
\item[1)] Si $v$ est finie, ${\widetilde{M}}_v (P)$ co\"{\i}ncide avec $\Mbar_v(P)$ et est définie comme le maximum des valeurs absolues $v$-adiques des coefficients de $P$; \\
\item[2)] Si $v$ est infinie et associée à un plongement
$\sigma\!\!: K \plongement \mathbb C$, ${\widetilde{M}}_v (P)$
désigne toujours le maximum des valeurs absolues $v$-adique des
coefficients de $P$ et:
$$
\Mbar_v(P) := \Mbar(P^{\sigma}) ,
$$
o\`u $P^{\sigma}$ désigne le polyn\^ome de $\mathbb{C}[X_1 , \dots
, X_n]$ obtenu à partir de $P$ en remplaçant chaque coefficient de
$P$ par son image par $\sigma$.
\end{description}
On définit aussi $M_v$ la mesure $v$-adique des formes
(homog\`enes ou multihomog\`enes) à coefficients dans $K$,
associant à toute telle forme $P$ le réel positif:
$$M_v (P) := \begin{cases}
{\widetilde{M}}_v (P) & \text{si $v$ est finie} \\
M(P^{\sigma}) & \text{si $v$ est infinie, associée à un plongement
$\sigma\: K \plongement \mathbb C$}
\end{cases} .
$$
{\bf Avertissement.---} $~~$On désigne -dans les prochains lemmes
qui suivent- par $H_v$ la mesure $v$-adique ${\widetilde{M}}_v$
définie ci-dessus et par $L_v$ l'application de $K[X_1 , \dots ,
X_n]$ dans ${\mathbb R}^+$ associant à tout polyn\^ome $P$, la
somme des valeurs absolues $v$-adiques de ses coefficients. $H_v$
et $L_v$ sont appelées ``hauteur'' et ``longueur'' $v$-adique et
leurs logarithmes -désignés respectivement par $h_v$ et $\ell_v$-
sont appelés hauteur logarithmique $v$-adique et longueur
logarithmique $v$-adique.
\subsubsection{{\bf Hauteurs:}}
Soit $K$ un corps de nombres. On définit $\widetilde{h}$ et
$\hbar$ les hauteurs dites respectivement de Gauss-Weil (ou
simplement de Weil) et de Mahler, associant à tout polyn\^ome $P$
de $K[X_1 , \dots , X_n]$ $(n \in {\mathbb N}^*)$ les réels:
\begin{align}
\widetilde{h}(P) &:= \sum_{v \in M_K} \frac{[K_v : {\mathbb Q}_v]}{[K : \mathbb Q]} \log{{\widetilde{M}}_v (P)} \notag \\
\hbar(P) &:= \sum_{v \in M_K} \frac{[K_v : {\mathbb Q}_v]}{[K :
\mathbb Q]} \log{\Mbar_v(P)} \notag
\end{align}
o\`u $M_K$ désigne l'ensemble des places de $K$, normalisées de
sorte que ${\nb{2}}_v = 2$ lorsque $v$ est infinie et ${\nb{p}}_v
= p^{-1}$ lorsque $v$ est finie et étend la place $p$ de
$\mathbb{Q}$. \\
On définit aussi $h$ la hauteur dite unitaire, qui
associe à toute forme $P$ (homog\`ene ou multihomog\`ene) à
coefficients dans $K$, le réel $h(P)$ défini par:
$$h(P) := \sum_{v \in M_K} \frac{[K_v : {\mathbb Q}_v]}{[K : \mathbb Q]} \log{M_v(P)} .$$
Notons que ces trois hauteurs $\widetilde{h}$, $\hbar$ et $h$ sont positives, indépendantes du corps de nombres $K$ choisi et invariantes par homothétie d'un élément de $K$. \\
Pour tout entier $n \geq 1$ et tout $P \in \Qbar[X_1 , \dots ,
X_n]$ (ou un polyn\^ome homog\`ene de $\Qbar[X_0 , \dots , X_n]$)
de degré total $d$, le lemme \ref{a.7} nous permet d'obtenir
immédiatement la relation de comparaison entre $\widetilde{h}(P)$
et $\hbar(P)$ suivante:
\begin{equation}
\hbar(P) - \frac{1}{2} \log{\!\binom{d + n}{n}} \leq
\widetilde{h}(P) \leq \hbar(P) + n d \log{2} . \label{1.19}
\end{equation}
Par ailleurs, dans le cas o\`u $P$ est une forme de $\Qbar[X_0 ,
\dots , X_n]$ $(n \geq 1)$, il est montré dans \cite{Le}
(théor\`eme $4$) qu'on a:
\begin{equation}
\hbar(P) \leq h(P) \leq \hbar(P) + d \sum_{j = 1}^{n} \frac{1}{2
j} . \label{1.20}
\end{equation}
Finalement, toujours dans le cas o\`u $P$ est une forme de degré
$d$ de $\Qbar[X_0 , \dots , X_n]$ $(n \geq 1)$, la relation de
comparaison entre $\widetilde{h}(P)$ et $h(P)$ résultant
immédiatement de (\ref{1.19}) et (\ref{1.20}) est:
\begin{equation}
h(P) - \frac{1}{2} \log{\!\binom{d + n}{n}} - d \sum_{j = 1}^{n}
\frac{1}{2 j} \leq \widetilde{h}(P) \leq h(P) + n d \log{2} .
\label{1.21}
\end{equation}
\subsection{Hauteurs de variétés projectives:}
Soit $V$ une variété algébrique définie sur un corps de nombres
$K$ et plongée dans un espace projectif ${\mathbb P}_N$ via un
plongement $\varphi$. La hauteur de Gauss-Weil (resp de Mahler,
unitaire) de $V$ est définie comme étant la hauteur de Gauss-Weil
(resp de Mahler, unitaire) de la forme éliminante $f_V$ de l'idéal
de définition de $\varphi(V)$ dans $K[X_0 , \dots , X_N]$. Cette
hauteur dépend évidemment du plongement $\varphi$ mais elle ne
dépend pas du choix de la forme éliminante gr\^ace à la formule du
produit sur $K$, on la note ${\widetilde{h}}_{\varphi}(V)$ (resp
$\hbar_{\varphi}(V)$, $h_{\varphi}(V)$). On a donc:
$$
{\widetilde{h}}_{\varphi}(V) ~\!:=~\! \widetilde{h}(f_V) ~,~
\hbar_{\varphi}(V) ~\!:=~\! \hbar(f_V) ~~\!\text{et}~~\!
h_{\varphi}(V) := h(f_V) .
$$
\section{{\bf ESTIMATIONS SUR LES COEFFICIENTS DES FONCTIONS ALG\'EBRIQUES}}
\subsection*{Quelques notations:}~\\ $\bullet$ Lorsque $n \in {\mathbb{N}}^*$ et $I = (i_1 ,
\dots i_n) \in {\mathbb{N}}^n$, on désignera par $\nb{I}$ la somme
des coordonnées de $I$, c'est-à-dire:
$$\nb{I} := i_1 + \dots + i_n .$$
On appellera ce dernier nombre: {\it longueur de $I$}.\\
$\bullet$ On désignera aussi par $I!$ le nombre:
$$I! := i_1! \dots i_n! .$$
$\bullet$ Si de plus $\Xsoul = (X_1 , \dots , X_n)$ est un
ensemble de $n$ variables indépendantes, on notera $\Xsoul^I$ le
mon\^ome:
$$\Xsoul^I := X_{1}^{i_1} \dots X_{n}^{i_n} ,$$
comme on posera:
$$\partial \Xsoul^I := \partial X_{1}^{i_1} \dots \partial X_{n}^{i_n}$$
et on notera enfin par $D^I$ l'opérateur de la {\it dérivation
divisée d'ordre $I$}:
$$D^I := \frac{1}{I!} \frac{\partial^{i_1}}{\partial X_{1}^{i_1}} \dots \frac{\partial^{i_n}}{\partial X_{n}^{i_n}} .$$
\begin{lemma}\label{a.1}
Soit $K$ un corps de nombres et $P$ un polyn\^ome de $K[T , Y_1 ,
\dots , Y_n] ~(n \in {\mathbb N}^*)$ de degré total majoré par $d
~(d \in {\mathbb N}^*)$. On pose $Y := (Y_1 , \dots , Y_n)$.
Alors, en tout point $\underline{x} = (t , \underline{y})$ de la
sous-vari\'et\'e d\'efinie par $P$ o\`u la d\'eriv\'ee partielle
de $P$ par rapport à $T$ est non nulle (i.e: $\frac{\partial
P}{\partial T} (\underline x) \neq 0$), on a pour tout $m \in
{\mathbb N}^*$ et tout $I \in {\mathbb N}^n$ tel que $\nb{I} = m$:
$$\frac{1}{I!} \frac{\partial^m T}{{\partial Y}^I} (
\underline{x}) = \frac{P_I ( \underline{x} )}{
{\left(\frac{\partial P}{\partial T}( \underline{x}
)\right)}^{2m-1}} .$$ o\`u $P_I$ est un polyn\^ome de $K[T , Y_1 ,
\dots , Y_n]$, de degré total majoré par $(2 m - 1) (d - 1)$ et
pour toute place infinie $v$ de $K$, $P_I$ est de longueur
$v$-adique $L_v$ majorée par:
$$L_v(P_I) \leq (8 n)^{m - 1} d^{3 m - 2} {L_v(P)}^{\! 2 m - 1} .$$
\end{lemma}
{\bf Démonstration.---} On proc\`ede par r\'ecurrence sur
$m\in\mathbb{N}^*$. Pour $m=1$, on a $\nb{I} = 1$ donc $I$ est de
la forme $I = (0 , \ldots , 0 , 1 , 0 , \ldots , 0)$ et $I!=1$.
Ainsi, on a:
$$\frac{1}{I!}\frac{\partial^1 T}{{\partial Y}^I} = \frac{\partial T
}{\partial Y_k} = \frac{-\frac{\partial P}{\partial
Y_k}}{\frac{\partial P}{\partial T}} .$$ Il suffit de poser $P_I =
- \frac{\partial P}{\partial Y_k}$ qui est de degr\'e total majoré
par $d$ et de longueur $v$-adique:
$$L_v \left(-\frac{\partial P}{\partial Y_k}\right) \leq d L_v(P) .$$
Donc le lemme \ref{a.1} est vrai pour $m=1$. \\

Supposons maintenant que les majorations du lemme \ref{a.1} sont
vraies jusqu'à un certain entier $m \geq 1$ et montrons ces
majorations pour $m+1$. Soit $I'=(i_1 , \ldots , i_n) \in {\mathbb
N}^n$ quelconque tel que $\nb{I'} = m + 1$. Choisissons $k \in \{1
, \dots , n\}$ tel que $i_k = \max(i_1 , \ldots , i_n)$. Puisque
$I' \neq 0$ (car $\nb{I'} = m + 1$), il est clair que $i_k \geq
1$. Posons $I = (i_1 , i_2 , \ldots , i_k-1 , \ldots , i_n) \in
{\mathbb N}^n$, donc $\nb{I} = m$. \\ Remarquons qu'on a $i_k \geq
\frac{i_1 + \ldots + i_n}{n} = \frac{m+1}{n}$, d'o\`u
\begin{equation}
\frac{1}{i_k} \leq \frac{n}{m+1} . \label{1.1}
\end{equation}
\begin{eqnarray*}
\text{On a:}~~~~~~~~~~~~~~~~~~~~~~~~~~~~~~~~\frac{1}{I'!}
\frac{\partial^{m + 1} T}{{\partial Y}^{I'}} & = & \frac{1}{I'!}
\frac{\partial}{\partial Y_k} \left(\frac{\partial^m T}{{\partial
Y}^I}\right) ~~~~~~~~~~~~~~~~~~~~~~~~~~~~~~~~~~~~~~\\
& = & \frac{I!}{I'!} \frac{\partial}{\partial Y_k} \left(\frac{1}{I!} \frac{\partial^m T}{{\partial Y}^I}\right)\\
& = & \frac{1}{i_k} \frac{\partial}{\partial Y_k}
\left(\frac{1}{I!} \frac{\partial^m T}{{\partial Y}^I}\right) .
\end{eqnarray*}
Or, d'apr\`es l'hypoth\`ese de récurrence: $$\frac{1}{I!}
\frac{\partial^m T}{{\partial Y}^I} = \frac{P_I}{(\frac{\partial P
}{\partial T})^{2m-1}} .$$ D'o\`u:
$$\frac{\partial}{\partial Y_k} \!\left(\frac{1}{I!} \frac{\partial^m T}{{\partial Y}^I}\right) = \frac{\left(\frac{\partial P_I}{\partial
Y_k}+\frac{\partial T}{\partial Y_k} \frac{\partial P_I}{\partial
T }\right)\!\!\left(\frac{\partial P}{\partial T}\right)^{2m-1}
\!-(2m-1)P_I\left(\frac{\partial P}{\partial T}\right)^{2m-2}
\!\!\left(\frac{\partial^2 P}{\partial Y_k \partial
T}+\frac{\partial T }{\partial Y_k} \frac{\partial^2 P}{{\partial
T}^2}\right)}{\left(\frac{\partial P}{\partial
T}\right)^{\!2(2m-1)}} .$$ En remplaçant, dans le membre de droite
de cette égalité, $\frac{\partial T}{\partial Y_k}$ par
$\frac{-\frac{\partial P}{\partial Y_k}}{\frac{\partial
P}{\partial T}}$ et en simplifiant la fraction obtenue -en
divisant son numérateur et le dénominateur par $(\frac{\partial
P}{\partial T})^{2m-3}$- on trouve:
$$\frac{\partial}{\partial Y_k} \left(\frac{1}{I!} \frac{\partial^m T}{{\partial Y}^I}\right) = \frac{\frac{\partial P_I}{\partial Y_k}
\left(\frac{\partial P}{\partial T}\right)^2 -\frac{\partial P_I
}{\partial T} \frac{\partial P}{\partial T} \frac{\partial P
}{\partial Y_k} - (2m-1)P_I \left[\frac{\partial P}{\partial T }
\frac{\partial^2 P}{\partial Y_k \partial T} - \frac{\partial
P}{\partial Y_k} \frac{\partial^2 P}{{\partial T}^2}
\right]}{\left(\frac{\partial P}{\partial T}\right)^{2m+1}} .$$ Il
suffit alors de poser:
$$P_{I'} = \frac{1}{i_k}
\!\left[\frac{\partial P_I}{\partial Y_k} \!\left(\frac{\partial P
}{\partial T}\right)^{\!2} \!- \frac{\partial P_I}{\partial T}
\frac{\partial P}{\partial T} \frac{\partial P}{\partial Y_k } -
(2m-1) P_I \!\left(\frac{\partial P}{\partial T} \frac{\partial^2
P }{\partial Y_k \partial T} - \frac{\partial P}{\partial Y_k}
\frac{\partial^2 P}{{\partial T}^2} \right) \!\right] \!.$$ Les
majorations sur les degr\'es se vérifient facilement puisque les
quatres polyn\^omes ${\left(\frac{\partial P}{\partial
T}\right)}^2 , \frac{\partial P}{\partial T} \frac{\partial
P}{\partial Y_k} , \frac{\partial P}{\partial T} \frac{\partial^2
P}{\partial Y_k \partial T}$ et $\frac{\partial P}{\partial Y_k}
\frac{\partial^2 P}{{\partial
T}^2}$ sont tous de degr\'e total majoré par $2 (d - 1)$. \\

Maintenant, étant donné une place infinie $v$ de $K$, pour obtenir
la longueur $v$-adique $L_v$ de $P_I$, on utilise les estimations
bien connues:
\begin{eqnarray*}
L_v(A + B) & \leq & L_v(A) + L_v(B) , \\
L_v(A . B) & \leq & L_v(A) . L_v(B)
\end{eqnarray*}
pour tout $A , B \in K[T , Y_1 , \dots , Y_n]$. En tenant compte,
de plus, de (\ref{1.1}) on a:
$$
L_v(P_{I'}) ~\!\leq~\! \frac{n}{m + 1} \left[4 (2 m - 1) d^3
{L_v(P)}^2 L_v(P_I) \right] ~\!\leq~\! 8 n d^3 {L_v(P)}^2 L_v(P_I)
.
$$
Comme d'apr\`es notre hypoth\`ese de récurrence:
$$L_v(P_I) \leq (8 n)^{m - 1} d^{3 m - 2} {L_v(P)}^{2 m - 1} ,$$
alors:
$$L_v(P_{I'}) \leq (8 n)^m d^{3 m + 1} {L_v(P)}^{2 m + 1} .$$
Ce qui montre l'estimation du lemme pour la longueur $v$-adique de
$P_{I'}$ et ach\`eve cette démonstration.
$~~~~\blacksquare$\vspace{1mm}
\begin{lemma}\label{a.2}
Soit $A$ un anneau commutatif unitaire int\`egre, $B := A[Y_1 ,
\ldots , Y_D][T]$ et $P$ un polyn\^ome de $B$ tel que
$\frac{\partial P}{\partial T} \neq 0$. On pose $\Xsoul := (X_1 ,
\dots , X_D)$ et $\Ysoul := (Y_1 , \dots , Y_D)$. Alors, il existe
une unique s\'erie formelle $U(\underline{X})\in S^{-1}B[[X_1 ,
\ldots , X_D]]$ satisfaisant:
$$P\left(\underline{Y}+\underline{X} , T + U(\underline{X})
\right) = P\left(\underline{Y} , T\right) ,$$ o\`u $S$ est la
partie multiplicative de $B$ engendr\'ee par $\frac{\partial
P}{\partial T}$.
\end{lemma}
{\bf Démonstration.---} La construction de la s\'erie formelle
$U(\underline{X})$ se fait pas \`a pas en construisant au
$n$-i\`eme pas la partie homog\`ene de degr\'e $n$ de
$U(\underline{X})$. En appelant, pour tout $n$ dans $\mathbb{N}$,
$V_n(\underline{X})$ la partie homog\`ene de degr\'e $n$ de
$U(\underline{X})$, on a:
$$U(\underline{X}) = V_0(\underline{X}) + V_1(\underline{X}) + \ldots + V_n(\underline{X}) + \ldots$$
Plus précisément, on construit pour tout $n$ dans $\mathbb{N}$ un
polyn\^ome $U_n(\underline{X})$ dans $S^{-1}B[X_1, \ldots, X_D]$
de degr\'e $\leq n$ tel que $\left(P(\underline{Y} + \underline{X}
, T + U_n(\underline{X})) - P(\underline{Y} , T )\right)$ soit
d'ordre $\geq (n+1)$ en $\underline{X}$, et tel que le
$(n+1)$-i\`eme polyn\^ome $U_{n+1}(\underline{X})$ soit une somme
de $U_n(\underline{X})$ et d'un polyn\^ome homog\`ene de
$S^{-1}B[X_1, \ldots, X_D]$ de degr\'e $n+1$. Ainsi
$U_{n+1}(\underline{X}) - U_n(\underline{X})$ est exactement la
partie homog\`ene de degr\'e $n+1$ de $U(\underline{X})$. Il
suffit alors de poser:
$$U(\underline{X}) = U_0(\underline{X}) + \left(U_1(\underline{X})
- U_0(\underline{X}) \right) + \left(U_2(\underline{X}) -
U_1(\underline{X}) \right) + \ldots$$ On aura
$P\left(\underline{Y} + \underline{X} , T +
U(\underline{X})\right) - P\left(\underline{Y} , T\right)$ de
degr\'e infini en $\underline{X}$, d'o\`u:\\ $P\left(\underline{Y}
+ \underline{X} , T + U(\underline{X})\right) =
P\left(\underline{Y} , T\right)$.

On constate d'abord que $U_0(\underline{X}) = U(\underline{0}) =
0$. Construisons les $U_n(\underline{X})$: Pour tout $G$ dans
$S^{-1}B[[X_1 , \ldots , X_D]]$, on pose
$$\Phi_G\left(\underline{Y} , T ,
\underline{X}\right) = P\left(\underline{Y} + \underline{X} , T +
G(\underline{X}) \right) - P\left(\underline{Y} , T\right) .$$ On
va montrer par r\'ecurrence que, pour tout $n$ dans $\mathbb{N}$,
il existe un polyn\^ome $U_n(\underline{X})$ dans $S^{-1}B[X_1,
\ldots, X_D]$ v\'erifiant:
\begin{description}
\item[$~~~1)$] ${\rm{ord}}_{\underline{X}}\Phi_{U_n} \geq n+1$;
\item[$~~~2)$] $U_0(\underline{X}) = 0 , U_n(\underline{X}) =
U_{n-1}(\underline{X}) + \sum_{\nb{\underline{\gamma}} = n}
 A_{\underline{\gamma}}{\underline{X}}^{\underline{\gamma}}
 $.
\end{description}
En effet, remarquons que pour $G$ dans $S^{-1}B[[\underline{X}]]$
on a:
\begin{equation}
\begin{split}
\Phi_G\left(\underline{Y} , T , \underline{X}\right) &= P\left(\underline{Y} + \underline{X} , T + G(\underline{X}) \right) - P\left(\underline{Y} , T\right) \\
&= \sum_{\Lambda = (\underline{\theta} , m)} D^\Lambda
P\left(\underline{Y} ,
T\right){\underline{X}}^{\underline{\theta}}G^m -
P\left(\underline{Y}
, T \right) \\
&= \sum_{\Lambda = (\underline{\theta}, m)\neq (0, 0)} D^\Lambda
P\left(\underline{Y} ,
T\right){\underline{X}}^{\underline{\theta}}G^m . \label{1.2}
\end{split}
\end{equation}
Pour $n = 0$ on a $U_0(\underline{X}) = 0$, d'o\`u d'apr\`es
(\ref{1.2}):
$$\Phi_{U_0}\left(\underline{Y} , T , \underline{X}\right) =
\sum_{\underline{\theta}\neq \underline{0}} D^{(\underline{\theta}
, 0)}P\left(\underline{Y} ,
T\right)\underline{X}^{\underline{\theta}} .$$ Donc
${\rm{ord}}_{\underline{X}}\Phi_{U_0} \geq 1$, par suite $1)$ et
$2)$ sont satisfaites pour $n = 0$. Supposons pour $n \geq 1$ que
$U_{n-1}(\underline{X})$ est construit satisfaisant $1)$,
c'est-\`a-dire: ${\rm{ord}}_{\underline{X}}\Phi_{U_{n-1}} \geq n$.
Par la contrainte $2)$ on pose:
$$U_n(\underline{X}) = U_{n-1}(\underline{X}) + \sum_{\nb{\underline{\gamma}} = n} A_{\underline{\gamma}}{\underline{X}}^{\underline{\gamma}} .$$
Il suffit alors de chercher les $A_\gamma$ dans $S^{-1}B$ de façon
 \`a avoir ${\rm{ord}}_{\underline{X}}\Phi_{U_n} \geq n+1$. Or, on a
 d'apr\`es (\ref{1.2}):
\begin{equation*}
\begin{split}
\Phi_{U_n} &= \sum_{\Lambda = (\underline{\theta} , m)
  \neq (\underline 0 , 0)} D^\Lambda P\left(\underline{Y} , T\right) \underline{X}^{\underline{\theta}}
  \left(U_{n-1}(\underline{X}) + \Sigma_{\nb{\underline{\gamma}} = n}
  A_{\underline{\gamma}}{\underline{X}}^{\underline{\gamma}}\right)^m \\
&\!\!\!\!\!\!\!\!\!\!\!= \Phi_{U_{n-1}} + \!\!\!\sum_{\Lambda =
(\underline{\theta} , m)
  \neq (\underline 0 , 0)}\!\!\!\!\!\!\!\!\!\!D^\Lambda P\left(\underline{Y} , T\right)
  \underline{X}^{\underline{\theta}}\left[
   \left(U_{n-1}(\underline{X}) + \Sigma_{\nb{\underline{\gamma}} = n} A_{\underline{\gamma}}
   {\underline{X}}^{\underline{\gamma}}\right)^m \!\!\!- \left(U_{n-1}(\underline{X})\right)^m \right] .
\end{split}
\end{equation*}
Comme l'expression
$$\sum_{\Lambda = (\underline{\theta} , m)
  \neq (0 , 0)}D^\Lambda P\left(\underline{Y} , T\right)
  \underline{X}^{\underline{\theta}}\left[
   \left(U_{n-1}(\underline{X}) + \Sigma_{\nb{\underline{\gamma}} = n} A_{\underline{\gamma}}
   {\underline{X}}^{\underline{\gamma}}\right)^m - \left(U_{n-1}(\underline{X})\right)^m \right]$$
vaut
$$\frac{\partial P}{\partial T}\left(\underline{Y} , T\right)\Sigma_{\nb{\underline{\gamma}} = n} A_{\underline{\gamma}}
   {\underline{X}}^{\underline{\gamma}}$$
pour $\Lambda = (\underline 0 , 1)$ et est une série en
$\underline X$ d'ordre $\geq n + 1$ pour tout $\Lambda =
(\underline{\theta} , m)$ \\ $\not\in \{(\underline 0 , 0) ,
(\underline 0 , 1)\}$, alors $\Phi_{U_n}$ s'écrit:
$$\Phi_{U_n} = \Phi_{U_{n-1}} + \frac{\partial P}{\partial T}\left(\underline{Y} , T\right)\Sigma_{\nb{\underline{\gamma}} = n} A_{\underline{\gamma}}
   {\underline{X}}^{\underline{\gamma}} + R_1(\underline{X})$$
avec $R_1(\underline{X})$ une série d'ordre $\geq n + 1$.
Mais, d'apr\`es l'hypoth\`ese de récurrence, \\
   ${\rm{{ord}}}_{\underline{X}}$$\!{\Phi}_{U_{n-1}} \geq n$, d'o\`u
   $\Phi_{U_{n-1}}(\underline{X})$ s'écrit sous la forme: $$\Phi_{U_{n-1}}(\underline{X}) = \sum_{\nb{\underline{\gamma}} = n} B_{\underline{\gamma}}{\underline{X}}^{\underline{\gamma}} + R_0(\underline{X})$$
   o\`u les $B_{\underline{\gamma}}$ sont dans $S^{-1}B$ et $R_0$
   est une série d'ordre $\geq n + 1$ en $\underline{X}$. D'o\`u: $$\Phi_{U_{n}}(\underline{X}) = \sum_{\nb{\underline{\gamma}} = n}
   B_{\underline{\gamma}}{\underline{X}}^{\underline{\gamma}} + P'_T\sum_{\nb{\underline{\gamma}} = n} A_{\underline{\gamma}}{\underline{X}}^{\underline{\gamma}} + R_2(\underline{X})$$ o\`u
   $R_2$ est une série d'ordre $\geq n + 1$ en $\underline{X}$ $(R_2 = R_1 + R_0)$. Il faut et il suffit alors de prendre, pour tout $\underline{\gamma}$: $$A_{\underline{\gamma}} = - \frac{B_{\underline{\gamma}}}{P'_T}$$
En r\'esum\'e, la construction des $U_n(\underline{X})$ se
   fait tout simplement comme suit: $U_0(\underline{X}) = 0$
   et une fois $U_{n - 1}(\underline{X})$ construit, on
   consid\`ere la s\'erie $\Phi_{U_{n - 1}}(\underline{X})$ qui
   est d'ordre $\geq n$ en $\underline{X}$ et soit $\sum_{\nb{\underline{\gamma}} = n}B_{\underline{\gamma}}{\underline{X}}^{\underline{\gamma}}$
   sa partie homog\`ene de degr\'e $n$, on prend alors:
\begin{equation}
    U_n(\underline{X}) = U_{n - 1}(\underline{X}) + \sum_{\nb{\underline{\gamma}} = n} - \frac{B_\gamma}{P'_T}{\underline{X}}^{\underline{\gamma}} ~. \label{1.3}
\end{equation}
L'unicit\'e des $U_n(\underline{X})$ est \'evidente car on est
oblig\'e de prendre $U_0(\underline{X}) = 0$ et une fois que $U_{n
- 1}(\underline{X})$ est construit on a une et une seule façon de
choisir $U_n(\underline{X})$. La d\'ecomposition de
$U(\underline{X})$ en somme de ses parties homog\`enes est alors:
\begin{equation*}
\begin{split}
U(\underline{X}) &= U_0(\underline{X}) + \left(U_1(\underline{X})
- U_0(\underline{X}) \right) + \left(U_2(\underline{X}) -
U_1(\underline{X}) \right) \\
&\quad+ \ldots + \left(U_n(\underline{X}) - U_{n -
1}(\underline{X}) \right) + \ldots
\end{split}
\end{equation*}
ce qui termine la preuve du lemme. $~~~~\blacksquare$\vspace{1mm}

On note $V_n(\underline{X})$ $\left(\in S^{-1}B\left[X_1 , \ldots
, X_D \right]\right)$ la partie homog\`ene de $U(\underline{X})$
de degr\'e $n$, c'est-\`a-dire: $V_n(\underline{X}) =
U_n(\underline{X}) - U_{n-1}(\underline{X})$. Pour tout $k\geq 1$
entier, on note $\frac{1}{{P'_T}^k} B$ les \'el\'ements de
$S^{-1}B$ qui sont le produit de $\frac{1}{{P'_T}^k}$ par un
\'el\'ement de $B$. On remarquera que $\frac{1}{{P'_T}^k} B$ est
stable par l'addition, mais ce n'est pas un anneau.
\begin{corollary}\label{a.3}
Pour tout entier $n\geq 1$, les coefficients de
$U_n(\underline{X})$ sont dans $\frac{1}{{P'_T}^{2n - 1}} B$ et
pour tout entier $d \geq n + 1$, les coefficients de la partie
homog\`ene de degr\'e $~d~$ de $~\Phi_{U_n}(\underline{X})~$ sont
dans $~\frac{1}{{P'_T}^{2d - 2 }} B$.
\end{corollary}
{\bf Démonstration.---} On proc\`ede par récurrence sur $n$. Pour
$n = 1$, comme on a $U_0(\underline{X}) = 0$ alors
$$\Phi_{U_0}(\underline{X}) = \sum_{\underline{\theta} \neq
\underline{0}} D^{(\underline{\theta} , 0)} P\left(\underline{Y} ,
T\right) {\underline{X}}^{\underline{\theta}} .$$ Et, comme pour
tout $\underline{\theta} \neq \underline{0}$:
$D^{(\underline{\theta} , 0)} P\left(\underline{Y} , T\right)$ est
dans $B$, alors $\Phi_{U_0}(\underline{X})$ est \`a coefficients
dans $B$. Donc, d'apr\`es (\ref{1.3}) les coefficients de
$U_1(\underline{X})$ sont dans $\frac{1}{P'_T} B$ et
$U_1(\underline{X})$ s'écrit:
$$U_1(\underline{X}) = U_0(\underline{X}) + \sum_{\nb{\underline{\gamma}}
= 1 } A_{\underline{\gamma}} {\underline{X}}^{\underline{\gamma}}
,$$ avec $A_{\underline{\gamma}} \in \frac{1}{P'_T} B$ pour tout
$\gammasoul$. Donc:
$$U_1(\underline{X}) = \sum_{\nb{\underline{\gamma}} = 1 }
A_{\underline{\gamma}} {\underline{X}}^{\underline{\gamma}} .$$
Par suite:
$$\Phi_{U_1}(\underline{X}) = \sum_{\Lambda =
(\underline{\theta} , m) \neq (0 , 0)} D^\Lambda
P\left(\underline{Y} , T\right)
{\underline{X}}^{\underline{\theta}}
\left(\Sigma_{\nb{\underline{\gamma}} = 1} A_{\underline{\gamma}}
{\underline{X}}^{\underline{\gamma}} \right)^m .$$ Soit maintenant
$d \geq 2$ un entier. Les coefficients de la partie homog\`ene de
degr\'e $d$ de $\Phi_{U_1}(\underline{X})$ sont des
coefficients des polyn\^omes homog\`enes de la forme: \\
$D^\Lambda P\left(\underline{Y} , T\right)
{\underline{X}}^{\underline{\theta}}
\left(\Sigma_{\nb{\underline{\gamma}} = 1} A_{\underline{\gamma}}
{\underline{X}}^{\underline{\gamma}} \right)^m$ avec
$\nb{\underline{\theta}} + m = d$. Donc, comme les
$A_{\underline{\gamma}}$ sont dans $\frac{1}{P'_T} B$, on déduit
que ces coefficients sont dans un certain $\frac{1}{{P'_T}^m} B$
$(m \leq d)$ donc aussi dans $\frac{1}{{P'_T}^d} B$. Et, comme $d
\leq 2d - 2$ car $(d \geq 2)$, on déduit que les coefficients de
la partie homog\`ene de degr\'e $d$ de $\Phi_{U_1}(\underline{X})$
sont -a fortiori- dans $\frac{1}{{P'_T}^{2d - 2}} B$, le
corollaire est vrai pour $n = 1$.

Soit maintenant $n\geq 2$ un entier. Supposons que le corollaire
est vrai pour tout $k \leq n -1$ et montrons le pour l'entier $n$.
L'hypoth\`ese de récurrence pour $n - 1$ implique, en prenant $d =
n \geq (n - 1) + 1$, que les coefficients de la partie homog\`ene
de degr\'e $n$ de $\Phi_{U_{n - 1}}(\underline{X})$ sont dans
$\frac{1}{{P'_T}^{2n - 2}} B$ d'o\`u, d'apr\`es (\ref{1.3}), les
coefficients de $U_n(\underline{X})$ sont dans
$\frac{1}{{P'_T}^{2n - 1}} B$ ce qui montre la premi\`ere
assertion du corollaire pour l'entier $n$.

D'autre part, soit $d \geq n + 1$ un entier. Ecrivons la
d\'ecomposition homog\`ene de $U_n(\underline{X})$:
$$U_n(\underline{X}) = V_1(\underline{X}) + \dots
V_n(\underline{X}) ,$$ o\`u pour tout $\alpha$ $(1 \leq \alpha
\leq n )$, $V_\alpha (\underline{X})$ désigne la partie homog\`ene
de degr\'e $\alpha$ de $U_n(\underline{X})$. On sait que:
$V_1(\underline{X}) = U_1(\underline{X}) , V_2(\underline{X}) =
U_2(\underline{X}) - U_1(\underline{X}) , \dots ,
V_n(\underline{X}) = U_n(\underline{X}) - U_{n -
1}(\underline{X})$. D'apr\`es l'hypoth\`ese de r\'ecurrence, pour
tout $\alpha$ $(1 \leq \alpha \leq n - 1)$, les coefficients de
$V_\alpha (\underline{X})$ sont dans $\frac{1}{{P'_T}^{2\alpha - 1
}} B$ et d'apr\`es la premi\`ere assertion du corollaire pour
l'entier $n$ d\'eja prouv\'ee ci-dessus, les coefficients de
$V_n(\underline{X})$ sont dans
$\frac{1}{{P'_T}^{2n - 1}} B$. \\
On écrit:
\begin{equation*}
\begin{split}
\Phi_{U_n}(\underline{X}) &= \sum_{\Lambda = (\underline{\theta} , m) \neq (0 , 0)} D^\Lambda P\left(\underline{Y} , T\right) {\underline{X}}^{\underline{\theta}} \left(U_n(\underline{X})\right)^m \\
&= \sum_{\Lambda = (\underline{\theta} , m) \neq (0 , 0)}
D^\Lambda P\left( \underline{Y} , T\right)
{\underline{X}}^{\underline{\theta}} \left(V_1(\underline{X}) +
\dots + V_n(\underline{X})\right)^m .
\end{split}
\end{equation*}
La partie homog\`ene de degr\'e $d \geq n + 1$ de
$\Phi_{U_n}(\underline{X})$ est donc une somme de polyn\^omes
homog\`enes de degr\'e $d$ de la forme:
$$D^\Lambda P\left( \underline{Y} , T\right)
{\underline{X}}^{\underline{\theta}} {V_1(\underline{X})}^{k_1}
\dots {V_n(\underline{X})}^{k_n}$$ avec: $$k_1 + \dots + k_n = m$$
et
$$\nb{\underline{\theta}} + k_1 + 2 k_2 + \dots + n k_n = d .$$
Les coefficients de tels polyn\^omes sont dans
$$\frac{1}{{P'_T}^{\left(2.1 - 1\right)k_1 + \left(2.2 -
1\right) k_2 + \dots + \left(2\alpha - 1\right) k_\alpha + \dots +
\left(2n - 1\right)k_n }} B ;$$ or:
\begin{equation*}
\begin{split}
\left(2.1 - 1\right) k_1 + \left(2.2 - 1\right) k_2 + \dots +
\left(2\alpha -
1\right) k_\alpha + \dots + \left(2n - 1\right) k_n \\
&\!\!\!\!\!\!\!\!\!\!\!\!\!\!\!\!\!\!\!\!\!\!\!\!\!\!\!\!\!\!\!\!\!\!\!\!\!\!\!\!\!\!\!\!\!\!\!\!\!\!\!\!\!\!\!\!\!\!\!\!\!\!\!\!\!\!\!\!\!\!\!\!\!\!\!\!\!\!\!\!\!\!\!\!\!=
2\left(k_1 + 2 k_2 + \dots + n k_n\right) - \left(k_1 +
\dots + k_n \right) \\
&\!\!\!\!\!\!\!\!\!\!\!\!\!\!\!\!\!\!\!\!\!\!\!\!\!\!\!\!\!\!\!\!\!\!\!\!\!\!\!\!\!\!\!\!\!\!\!\!\!\!\!\!\!\!\!\!\!\!\!\!\!\!\!\!\!\!\!\!\!\!\!\!\!\!\!\!\!\!\!\!\!\!\!\!\!=
2\left(d - \nb{\underline{\theta}}\right)
- m \\
&\!\!\!\!\!\!\!\!\!\!\!\!\!\!\!\!\!\!\!\!\!\!\!\!\!\!\!\!\!\!\!\!\!\!\!\!\!\!\!\!\!\!\!\!\!\!\!\!\!\!\!\!\!\!\!\!\!\!\!\!\!\!\!\!\!\!\!\!\!\!\!\!\!\!\!\!\!\!\!\!\!\!\!\!\!=
2 d - \left(2 \nb{\underline{\theta}} + m\right) .
\end{split}
\end{equation*}
On remarque que $2 \nb{\underline{\theta}} + m \geq 2$ car soit
$\nb{\underline{\theta}} \geq 1$ ou bien $\underline{\theta} =
\underline{0}$ et $m \geq 2$ (car $d \geq n + 1$). D'o\`u:
$$\left(2.1 - 1\right) k_1 + \left(2.2 - 1\right) k_2 + \dots +
\left(2 \alpha - 1\right) k_\alpha + \dots + \left(2n - 1\right)
k_n \leq 2d - 2 .$$ Donc -a fortiori- les coefficients de tels
polyn\^omes sont dans $\frac{1}{{P'_T}^{2d - 2 }} B$, et enfin les
coefficients de la partie homog\`ene de degr\'e $d$ de
$\Phi_{U_n}(\underline{X})$ sont dans $\frac{1}{{P'_T}^{2d - 2}}
B$. Ceci termine la preuve du corollaire.
$~~~~\blacksquare$\vspace{1mm}
\begin{corollary}\label{a.4}
Si $$U(\underline{X}) = \sum_\Lambda \alpha_\Lambda
{\underline{X}}^\Lambda$$ o\`u les $\alpha_\Lambda$ sont dans
$S^{-1} B$, alors pour tout $\Lambda$ on a: $$\alpha_\Lambda
\left(\frac{\partial P}{\partial T}\right)^{2 \nb{\Lambda} - 1}
\in B .$$
\end{corollary}
{\bf Démonstration.---} C'est une conséquence immédiate du
corollaire \ref{a.3}. $~~~~\blacksquare$\vspace{1mm}
\begin{lemma}\label{a.5}
Soient $K$ un corps de nombres et $P$ un polyn\^ome de $K[Y_1 ,
\dots , Y_n][T]$ de degré total majoré par $d$. Soit aussi
$\underline{x} = (\underline{y} , t)$ un point à coordonnées dans
$K$ de la sous-vari\'et\'e affine d\'efinie par $P$ tel que:
$$\frac{\partial p}{\partial T}(\underline{x}) \neq 0 ,$$
et $T = T(\underline{Y})$ la branche de la sous-vari\'et\'e
d'\'equation $P=0$ contenant le point $\underline{x}$. En posant
pour tout $I \in \mathbb{N}^n$:
$$a_I := \frac{1}{I!} \frac{\partial^{\nb{I}} T}{{\partial \underline{Y}}^I}(\underline{y}) \in K ,$$
on a pour tout $m \in \mathbb N$ et toute place $v$ de $K$:
\begin{equation}
\begin{split}
\max\left(1 , {\nb{a_I}}_v ; \nb{I} \leq m\right) \\
&\!\!\!\!\!\!\!\!\!\!\!\!\!\!\!\!\!\!\!\!\!\!\!\!\!\!\!\!\!\!\!\!\!\!\!\!\!\!\!\!\!\!\!\!\!\!\!\!\!\leq
\left(8 n d^3 (d + 1)^{2 (n + 1)} {H_v(P)}^2 {H_v(1 ,
\underline{x})}^{2 (d - 1)} . {\max\!\left\{1 ,
{\nb{\frac{1}{\frac{\partial P}{\partial
T}(\underline{x})}}}_v\right\}}^{\!2} \right)^{\!\!m} . H_v(1 ,
\underline{x}) \label{1.4}
\end{split}
\end{equation}
si $v$ est infinie et:
\begin{equation}
\max\!\left(1 , {\nb{a_I}}_v ; \nb{I} \leq m\right) \leq
\left(\!{H_v(P)}^2 {H_v(1 , \underline{x})}^{2 (d - 1)}
{\max\!\left\{1 , {\nb{\frac{1}{\frac{\partial P}{\partial
T}(\underline{x})}}}_v\right\}}^{\!\!2}\right)^{\!\!m} \!\!. H_v(1
, \underline{x}) \label{1.5}
\end{equation}
si $v$ est finie. \\
De plus pour tout $m \in \mathbb N$ on a:
$$h(1 , a_I ; \nb{I} \leq m) \leq \!\left[4 \widetilde{h}(P) + 4 (d - 1) h(1 : \underline{x}) + (4 n + 9) (\log d + 1)\right] \!m + h(1 : \underline{x}) .$$
\end{lemma}
{\bf Démonstration.---} Supposons sans restreindre la généralité
que l'un des coefficients de $P$ vaut $1$ et écrivons:
$$P \left(\underline{Y} , T\right) := \sum_{\underline{j} = (j_0 , j_1 , \dots , j_n) \in \Lambda} b_{\underline{j}} T^{j_0} Y_{1}^{j_1} \dots Y_{n}^{j_n}$$
pour un certain sous-ensemble fini $\Lambda$ de ${\mathbb N}^{n + 1}$ formé des uplets de longueurs $\leq d$ et certains nombres $b_{\underline j}$ $(\underline j \in \Lambda)$ de $K$. \\
Pour une place donnée $v$ de $K$ et un entier $m \in \mathbb{N}$
on estime dans ce qui suit, selon le cas $v$ finie ou $v$ infinie,
le maximum entre $1$ et les valeurs absolues $v$-adiques des $a_I$
$(I \in \mathbb{N} , \nb{I} \leq m)$.
\parag{\bf{Premier cas (si $v$ est infinie)}}
D'apr\`es le lemme \ref{a.1}, pour tout $I \in {\mathbb N}^n$, $I
\neq 0$ et $\ell := \nb{I}$ on a:
$$a_I = \frac{P_I(\underline x)}{\left(\frac{\partial P}{\partial T}(\underline x)\right)^{2 \ell - 1}}$$
d'o\`u:
$${\nb{a_I}}_v = \frac{{\nb{P_I(\underline x)}}_v}{{\nb{\frac{\partial P}{\partial T}(\underline x)}}_{v}^{2 \ell - 1}} .$$
En écrivant:
$$P_I (\underline Y , T) = \sum_{\underline j \in J} c_{\underline j} T^{j_0} Y_{1}^{j_1} \dots Y_{n}^{j_n}$$
pour un certain sous-ensemble fini $J$ de ${\mathbb N}^{n + 1}$ et
certains nombres $c_{\underline j}$ $(\underline j \in J)$ dans
$K$, on a:
$$P_I (\underline x) = \sum_{\underline j \in \Lambda} c_{\underline j} t^{j_0} y_{1}^{j_1} \dots y_{n}^{j_n}$$
d'o\`u:
\begin{equation*}
\begin{split}
{\nb{P_I (\underline x)}}_v &\leq \sum_{\underline j \in J} {\nb{c_{\underline j}}}_v {\nb{t}}_{v}^{j_0} {\nb{y_1}}_{v}^{j_1} \dots {\nb{y_n}}_{v}^{j_n} \\
&\leq L_v(P_I) . {\max\left\{1 , {\nb{t}}_v , {\nb{y_1}}_v , \dots , {\nb{y_n}}_v\right\}}^{{d°}_{\!\!\!\rm{tot}} P_I} \\
&\leq (8 n)^{\ell - 1} d^{3 \ell - 2} {L_v(P)}^{2 \ell - 1} {H_v(1 , \underline x)}^{(2 \ell - 1) (d - 1)} ~~~~\text{(d'apr\`es le lemme \ref{a.1})} \\
&\leq \left(8 n d^3 {L_v(P)}^2 {H_v(1 , \underline{x})}^{2 (d -
1)}\right)^{\ell} .
\end{split}
\end{equation*}
Par suite:
\begin{equation*}
\begin{split}
{\nb{a_I}}_v &= {\nb{P_I(\underline{x})}}_v . {\nb{\frac{1}{\frac{\partial P}{\partial T}(\underline{x})}}}_{v}^{2 \ell - 1} \\
&\leq \left(8 n d^3 {L_v(P)}^2 {H_v(1 , \underline{x})}^{2 (d - 1)} . {\max\!\left\{1 , {\nb{\frac{1}{\frac{\partial P}{\partial T}(\underline{x})}}}_v\right\}}^{\!2} \right)^{\!\!\ell} \\
&\leq \left(8 n d^3 (d + 1)^{2 (n + 1)} {H_v(P)}^2 {H_v(1 ,
\underline{x})}^{2 (d - 1)} . {\max\!\left\{1 ,
{\nb{\frac{1}{\frac{\partial P}{\partial
T}(\underline{x})}}}_v\right\}}^{\!2} \right)^{\!\!\ell}
\end{split}
\end{equation*}
car, comme ${d°}_{\rm{tot}} P \leq d$, le polyn\^ome $P$ contient au plus $\binom{d + n + 1}{n + 1} \leq (d + 1)^{n + 1}$ mon\^omes et $L_v(P) \leq (d + 1)^{n + 1} H_v(P)$. \\
Par ailleurs, si $I = 0$ on a: $a_I = a_0 = t$ donc ${\nb{a_I}}_v
= {\nb{t}}_v \leq H_v(1 , \underline x)$. D'o\`u, pour tout $m \in
\mathbb{N}$:
\begin{equation*}
\begin{split}
\max\left(1 , {\nb{a_I}}_v ; \nb{I} \leq m\right) \\
&\!\!\!\!\!\!\!\!\!\!\!\!\!\!\!\!\!\!\!\!\!\!\!\!\!\!\!\!\!\!\!\!\!\!\!\!\!\!\!\!\!\!\!\!\!\!\!\!\!\leq
\left(8 n d^3 (d + 1)^{2 (n + 1)} {H_v(P)}^2 {H_v(1 ,
\underline{x})}^{2 (d - 1)} . {\max\!\left\{1 ,
{\nb{\frac{1}{\frac{\partial P}{\partial
T}(\underline{x})}}}_v\right\}}^{\!2} \right)^{\!\!m} . H_v(1 ,
\underline{x}) .
\end{split}
\end{equation*}
Ce qui est la relation (\ref{1.4}) du lemme.
\parag{\bf{Deuxi\`eme cas (si $v$ est finie)}}
Soit, dans ce cas, $A_v$ l'anneau des entiers $v$-adique de $K$:
$$A_v = \{\alpha \in K /~\! {\nb{\alpha}}_v \leq 1\} .$$
Choisissons un nombre algébrique $\ell_v(P) \in K$ tel que:
$$\ell_v(P) \in \left\{\frac{1}{b_{\underline j}} ; ~\underline j \in \Lambda , b_{\underline j} \neq 0\right\}$$
et
$${\nb{\ell_v(P)}}_v = \min\left\{{\nb{\frac{1}{b_{\underline j}}}}_v ; ~\underline j \in \Lambda , b_{\underline j} \neq 0\right\} .$$
Ainsi choisi, $\ell_v(P) \in A_v$ et pour tout $\underline j \in
\Lambda$: $\ell_v(P) b_{\underline j} \in A_v$, ce qui
entra{\sf\^\i}ne que le polyn\^ome $q(\underline Y , T) :=
\ell_v(P) P(\underline Y , T)$ appartient à $A_v[\underline Y ,
T]$ et de plus:
$$\frac{\partial q}{\partial T}(\underline x) = \ell_v(P) \frac{\partial P}{\partial T} (\underline x) \neq 0$$
d'apr\`es les hypoth\`eses du lemme \ref{a.5}. \\
En appliquant le lemme \ref{a.2} avec $A$ remplac\'e par $A_v$, on
d\'eduit l'existence et l'unicit\'e d'une s\'erie
$$U(\underline{X}) \in S^{-1} A_v \left[\underline{Y} , T\right] \left[[\underline{X}]\right]$$
(avec $\underline X := (X_1 , \dots , X_D)$ et $S$ est la partie
multiplicative de $A_v \left[\underline{Y} , T\right]$ engendr\'ee
par $\left(\frac{\partial q}{\partial T}\right)$) telle que:
$$q(\underline{Y} + \underline{X} , T + U(\underline{X})) =
q(\underline{Y} , T) ,$$ d'o\`u, en simplifiant les deux membres
de cette \'egalit\'e par $\ell_v(P)$, on a aussi:
$$P(\underline{Y} + \underline{X} , T + U(\underline{X})) =
P(\underline{Y} , T) .$$ En sp\'ecialisant dans cette derni\`ere
\'egalit\'e $T$ en $T(\underline{Y})$ (d\'efini dans le lemme
\ref{a.5}), on obtient:
$$P(\underline{Y} + \underline{X} , T(\underline{Y}) + U(\underline{X})) = 0 .$$
Et, comme la sp\'ecialisation: $(\underline{X} , \underline{Y}) =
(0,\underline{y})$ donne $\underline{Y} + \underline{X} =
\underline{y}$ et $T(\underline{Y}) + U(\underline{X})) = t$, on
d\'eduit que: $T(\underline{Y}) + U(\underline{X}) =
T(\underline{Y} + \underline{Y})$, d'o\`u:
\begin{equation*}
\begin{split}
U(\underline{X}) &= T(\underline{Y} + \underline{X}) - T(\underline{Y}) \\
&= \sum_{I \in {\mathbb N}^n} D^I T(\underline{Y}) . {\underline
X}^I - T(\underline Y) ,
\end{split}
\end{equation*}
o\`u $D^I T(\underline{Y})$ d\'esigne la d\'erivation divis\'ee
d'ordre $I$ de $T(\underline{Y})$, i.e $$ D^I T(\underline{Y}) :=
\frac{1}{I!}\frac{\partial^{\nb{I}}
T(\underline{Y})}{\partial\underline{Y}^I} .$$ On a:
$$U(\underline{X}) = \sum_{I\in {\mathbb N}^n, I \neq 0 } D^I T(\underline{Y}) . {\underline X}^I ,$$
mais, d'apr\`es le corollaire \ref{a.4}, étant donné $I \in
{\mathbb N}^n$, $I \neq 0$, si l'on pose $m := \nb{I}$, on a:
$$D^I T(\underline{Y}) . {\left(\frac{\partial q}{\partial T
}\right)}^{2m-1} \in A_v\left[\underline{Y} , T\right] ,$$
c'est-\`a-dire:
$$\left(\ell_v(P) \frac{\partial P}{\partial T}\right)^{2m-1} . D^I
T(\underline{Y}) \in A_v\left[\underline{Y} , T\right] .$$ Or,
d'apr\`es le lemme \ref{a.1}, ${\left(\ell_v(P) \frac{\partial
P}{\partial T}\right)}^{2m-1} . D^I T(\underline{Y})$ est un
polyn\^ome de degré total major\'es par $(2m-1)(d - 1)$. D'o\`u,
en sp\'ecialisant $(\underline{Y} , T)$ en $(\underline{y} , t) =
\underline{x} \in K^{n+1}$:
$${\left(\ell_v(P) \frac{\partial P}{\partial T}(\underline{x})\right)}^{2m-1} . a_I \in A_v\left[\underline{y} , t\right]$$ et plus pr\'ecisement $a_I$ s'\'ecrit sous la forme:
$$a_I = \frac{f(\underline{y} , t)}{{\left(\ell_v(P) \frac{\partial P}{\partial T}(\underline{x})\right)}^{2m-1}}$$ o\`u $f$ est un polyn\^ome de $A_v[\underline{Y} , T]$ de degré total major\'e par $(2m-1)(d - 1)$. Donc
$$
{\nb{a_I}}_v ~\!=~\!
\frac{{\nb{f(\underline{y},t)}}_v}{{\nb{\ell_v(P) \frac{\partial
P}{\partial T}(\underline{x})}}_{v}^{2m-1}} ~\!\leq~\!
\frac{{\max\left\{1 , {\nb{t}}_v , {\nb{y_1}}_v , \dots ,
{\nb{y_n}}_v\right\}}^{(2 m - 1)(d - 1)}}{{\nb{\ell_v(P)}}_{v}^{2
m - 1} . {\nb{\frac{\partial P}{\partial T}(\underline x)}}_{v}^{2
m - 1}} .
$$
Comme, par d\'efinition de $\ell_v(P)$,
$$\frac{1}{{\nb{\ell_v(P)}}_v} = \max\left\{1 , {\nb{b_{\underline{j}}}}_v , ~\underline{j} \in \Lambda\right\} =
H_v(P) ,$$ alors pour tout $I \in {\mathbb N}^n \setminus \{0\}$,
on a:
\begin{equation*}
\begin{split}
{\nb{a_I}}_v &\leq {H_v(P)}^{2 \nb{I} - 1} . {H_v(1 , \underline x)}^{(2 \nb{I} - 1)(d - 1)} . {\nb{\frac{1}{\frac{\partial P}{\partial T}(\underline x)}}}_{v}^{2 \nb{I} - 1} \\
&\leq \left({H_v(P)}^2 {H_v(1 , \underline{x})}^{2 (d - 1)}
{\max\!\left\{1 , {\nb{\frac{1}{\frac{\partial P}{\partial
T}(\underline{x})}}}_v\right\}}^{\!\!2}\right)^{\!\!\nb{I}} .
\end{split}
\end{equation*}
Maintenant, si $I = 0$, on a ${\nb{a_I}}_v = {\nb{t}}_v \leq H_v(1
, \underline x)$ d'o\`u, pour tout $m \in \mathbb{N}$:
\begin{equation*}
\max\!\left(1 , {\nb{a_I}}_v ; \nb{I} \leq m\right) \leq
\left(\!{H_v(P)}^2 {H_v(1 , \underline{x})}^{2 (d - 1)}
{\max\!\left\{1 , {\nb{\frac{1}{\frac{\partial P}{\partial
T}(\underline{x})}}}_v\right\}}^{\!\!2}\right)^{\!\!m} \!\!. H_v(1
, \underline{x}) .
\end{equation*}
Ce qui est la relation (\ref{1.5}) du lemme.

La partie restante du lemme s'obtient en reportant, pour toute
place $v$ de $K$, les estimations (\ref{1.4}) et (\ref{1.5}) (déjà
démontrées ci-dessus) dans la définition de la hauteur
logarithmique de la famille $\{1 , a_I ; \nb{I} \leq m\}$ $(m \in
\mathbb{N})$. En effet, on obtient ainsi pour tout $m \in
\mathbb{N}$:
\begin{equation*}
\begin{split}
h\!\left(1 , a_I ; \nb{I} \leq m\right) &\leq \left[2 \widetilde{h}(P) + 2 (d - 1) h(1 : \underline{x}) + 2 h\!\!\left(\!\!1 : \frac{1}{\frac{\partial P}{\partial T}(\underline x)}\right) \right. \\
&\quad ~~~+ \left. \log(8 n) + 3 \log{d} + 2 (n + 1) \log{(d + 1)}
\phantom{h\!\!\left(\!\!1 : \frac{1}{\frac{\partial P}{\partial
T}(\underline x)}\right)}
\!\!\!\!\!\!\!\!\!\!\!\!\!\!\!\!\!\!\!\!\!\!\!\!\!\!\!\!\!\!\!\!\right]
\! m + h(1 : \underline{x}) .
\end{split}
\end{equation*}
Comme on montre facilement, en raisonnant place par place, qu'on
a:
\begin{equation*}
\begin{split}
h\!\!\left(\!\!1 : \frac{1}{\frac{\partial P}{\partial T}(\underline x)}\right) &= h\!\!\left(1 : \frac{\partial P}{\partial T}(\underline x)\right) \\
&\leq \widetilde{h}(P) + (d - 1) h(1 : \underline x) + \log{\!\left(\!d \binom{d + n}{n + 1}\!\right)} \\
&\leq \widetilde{h}(P) + (d - 1) h(1 : \underline x) + (n + 2)
\log{d} ,
\end{split}
\end{equation*}
alors pour tout $m$ dans $\mathbb{N}$:
\begin{equation*}
\begin{split}
h\!\left(1 , a_I ; \nb{I} \leq m\right) \leq \left[4 \widetilde{h}(P) + 4 (d - 1) h(1 : \underline{x}) + \log{(8 n)} + (2 n + 7) \log{d} \right. \\
\left.+ (2 n + 2) \log{(d +
1)}\phantom{\widetilde{h}(P)}\!\!\!\!\!\!\!\!\!\!\!\right] m + h(1
: \underline{x}) .
\end{split}
\end{equation*}
Il ne reste qu'à remarquer que l'expression:\\ $\log{(8 n)} + (2 n
+ 7) \log{d} + (2 n + 2) \log{(d + 1)}$ est majorée
grossi\`erement par:\\ $(4 n + 9)(\log{d} + 1)$ pour avoir
finalement pour tout $m \in \mathbb{N}$:
\begin{equation*}
h\!\left(1 , a_I ; \nb{I} \leq m\right) \!\leq \!\left[4
\widetilde{h}(P) + 4 (d - 1) h(1 : \underline x) + (4 n +
9)(\log{d} + 1)\right] \!m + h(1 : \underline{x}).
\end{equation*}
Ce qui ach\`eve la démonstration. $~~~~\blacksquare$\vspace{1mm}
\section{{\bf APPLICATION AUX PARAM\'ETRISATIONS DE VARI\'ET\'ES}}
Soit $G$ une variété algébrique définie sur un corps de nombres
$K$, irréductible, plongée dans ${\mathbb P}_N$ $(N \in {\mathbb
N}^*)$, de dimension $g$ et $\mathbf e$ un point de $G$ représenté
dans ${\mathbb P}_N$ par un syst\`eme de coordonnées projectives
$\underline e = (e_0 : e_1 : \dots : e_N)$. Soient aussi
$K[\underline X]$ ($\underline X := (X_0 , \dots , X_N)$) l'anneau
des coordonnées de ${\mathbb P}_N$, $\mathfrak{G}$ l'idéal de
définition de l'adhérence de Zariski $\Gbar$ de $G$ dans
$K[\underline X]$ et $A := K[\underline X] / \mathfrak{G}$
l'anneau des coordonnées de $\Gbar$.

En faisant un changement de coordonnées convenable, nous pouvons
nous ramener aux hypoth\`eses suivantes:
\parag{\bf{Hypoth\`eses de normalisation}}
\begin{description}
\item[1)] $e_0 = 1$ et l'espace tangent $T_{\mathbf e} G$ de $G$
en $\mathbf e$ a pour équations:
$$X_{g + 1} = \dots = X_N = 0 .$$
\item[2)] Un point générique $(y_0 : y_1 : \dots : y_N)$ de $G
\hookrightarrow {\mathbb P}_N$ vérifie:
$$\mbox{deg}{\mbox{tr}}_K K(y_0 , \dots , y_g) = \mbox{deg}{\mbox{tr}}_K K(y_0 , \dots , y_N) = g + 1 ,$$
c'est-à-dire que $y_0 , \dots , y_g$ sont $K$-algébriquement
indépendants.
\end{description}
L'hypoth\`ese $2)$ entra{\sf\^\i}ne (suivant la terminologie de
\cite{Ph2}) que $G$ est incompl\`etement défini dans
${\mathbb{P}}_N$ par des équations de la forme:
$${\widetilde{P}}_i (X_0 , \dots , X_g , X_i) = 0 ~~~~ (g + 1 \leq i \leq N)$$
o\`u, pour tout $i \in \{g + 1 , \dots , N\}$, ${\widetilde{P}}_i$ est un polyn\^ome homog\`ene (non identiquement nul) de $K[X_0 , \dots , X_g , X_i]$ et vérifie ${\widetilde{P}}_i (e_0 , \dots , e_g , e_i) = 0$ (puisque $\mathbf e \in G$). Nous supposons sans perte de généralité que pour tout $i = g + 1 , \dots , N$, l'un au moins des coefficients de ${\widetilde{P}}_i$ vaut $1$. \\
Par suite, la forme de l'espace tangent de $G$ en $\mathbf e$
supposée en $1)$ entra\^{\i}ne qu'on a:
$$\frac{\partial {\widetilde P}_i}{\partial X_i} (e_0 , \dots , e_g , e_i) \neq 0 ~~~~~~~~ (\text{pour}~
g + 1 \leq i \leq N) .$$ Posons pour tout $g + 1 \leq i \leq N$:
$$P_i (X_1 , \dots , X_g , X_i) := {\widetilde{P}}_i (1 , X_1 , \dots , X_g , X_i) ,$$
alors on a aussi:
$$\frac{\partial P_i}{\partial X_i} (e_1 , \dots , e_g , e_i) = \frac{\partial {\widetilde{P}}_i}{\partial X_i}
(e_0 , \dots , e_g , e_i) \neq 0 ~~~~~~~~ (\text{pour}~ g + 1 \leq
i \leq N) .$$ Soit maintenant pour tout $g + 1 \leq i \leq N$,
$f_i$ la fonction algébrique en $X_1 , \dots , X_g$ définie
implicitement au voisinage de $(e_1 , \dots , e_g)$ par:
$$\left\{\!\!\!
\begin{array}{l}
f_i (e_1 , \dots , e_g) = e_i \\
P_i \!\left(X_1 , \dots , X_g , f_i (X_1 , \dots , X_g)\right) = 0
\end{array}
\right..$$ Le théor\`eme des fonctions implicites implique
l'existence et l'unicité des fonctions $f_i$ $(g + 1 \leq i \leq
N)$ gr\^ace au fait qu'on a pour tout $i \in \{g + 1 , \dots ,
N\}$: $P_i (e_1 , \dots , e_g , e_i) = 0$ et $\frac{\partial
P_i}{\partial X_i} (e_1 , \dots , e_g , e_i) \neq 0$. Ce
théor\`eme montre de plus que ces fonctions $f_i$ $(g + 1 \leq i
\leq N)$ sont développables en séries de Taylor au voisinage de
$(e_1 , \dots , e_g)$. Notons ces développements:
$$f_i (X_1 , \dots , X_g) = \sum_{I \in {\mathbb N}^g} a_{I}^{(i)} {\underline T}^I ~~~~ (g + 1 \leq i \leq N)$$
avec $\underline T = (T_1 , \dots , T_g) := (X_1 - e_1 , \dots , X_g - e_g)$. \\
En désignant par $({\vvec}_1 , \dots , {\vvec}_g)$ la base
canonique de ${\mathbb C}^g$, considérons de m\^eme les variables
$X_1 , \dots , X_g$ comme des polyn\^omes (donc des séries) en
$\underline{T}$, en écrivant pour tout $i = 1 , \dots , g$:
$$X_i = \sum_{I \in {\mathbb N}^g} a_{I}^{(i)} {\underline{T}}^I$$
avec
$$a_{I}^{(i)} = \begin{cases}
e_1 & \text{si $I = 0$} \\
1 & \text{si $I = {\vvec}_i$} \\
0 & \text{sinon}
\end{cases} .$$
Posons finalement:
$$a_{I}^{(0)} := \begin{cases}
1 & \text{si $I = 0$} \\
0 & \text{sinon}
\end{cases} .$$
On a ainsi une paramétrisation de $G$ au voisinage de $\mathbf e$,
utilisant comme param\`etres $T_1 , \dots , T_g$ et donnée par le
monomorphisme de paramétrisation $\varphi$ de $A = K[\underline X]
/ \mathfrak{G}$ dans $K[[\underline T]]$ suivant:
$$\varphi\: A \longrightarrow K[[\underline T]]$$
tel que:
$$\varphi(X_0) = 1 ,$$
pour $i = 1 , \dots , g$:
$$\varphi(X_i) = X_i = e_i + T_i = \sum_{I \in {\mathbb N}^g} a_{I}^{(i)} {\underline T}^I$$
et pour $i = g + 1 , \dots , N$:
$$\varphi(X_i) = f_i(\underline T) = \sum_{I \in {\mathbb N}^g} a_{I}^{(i)} {\underline T}^I .$$
Donc, pour tout $i = 0 , \dots , N$:
$$\varphi(X_i) = \sum_{I \in {\mathbb N}^g} a_{I}^{(i)} {\underline T}^I .$$
Pour une place donnée $v$ de $K$, l'estimation des valeurs
absolues $v$-adiques des nombres $a_{I}^{(i)}$ $(I \in {\mathbb
N}^g , i = 0 , \dots , N)$ est donnée par le lemme suivant qui
n'est qu'une conséquence immédiate du lemme \ref{a.5}.
\begin{lemma}\label{a.6}
Pour tout $I \in {\mathbb N}^g$, tout $i \in \{0 , \dots , N\}$ et
toute place $v$ de $K$ on a:
\begin{equation}
\begin{split}
\max\left(1 , {\nb{a_{I}^{(i)}}}_v\right) &\leq \left(\!\!\!\!\!\!\!\!\!\!\!\!\!\!\!\!\!\!\!\!\!\!\!\!\!\!\!\!\!\!\!\!\!\!\!\!\!\!\!\!\!\!\!\!\!\!\!\!\!\!\!\!\!\!\!\!\!\!\!\!\!\!\!\!\!\!\!\!\!\!\!\!\!\!\!\!\!\!\!\!\!\!\!\phantom{{\max\!\left\{1 , {\nb{\frac{1}{\frac{\partial P_{i_r}}{\partial X_{i_r}}(e_1 , \dots , e_g , e_{i_r})}}}_v\right\}}^{\!\!2}}8 g {d(G)}^3 (d(G) + 1)^{2 (g + 1)} . \!\!\!\prod_{j = g + 1}^{N} \!\!\!{H_v({\widetilde P}_j)}^2 . {H_v(\underline e)}^{2 (d(G) - 1)} \right.\\
&\quad \left.\times \!\!\!\prod_{j = g + 1}^{N}
\!\!\!{\max\!\left\{1 , {\nb{\frac{1}{\frac{\partial P_j}{\partial
X_j}(e_1 , \dots , e_g ,
e_j)}}}_v\right\}}^{\!\!2}\right)^{\!\!\nb{I}} \!\!\!.
H_v(\underline e)
\end{split} \label{1.23}
\end{equation}
si $v$ est infinie et
\begin{equation}
\begin{split}
\max\left(1 , {\nb{a_{I}^{(i)}}}_v\right) &\leq \left(\!\!\!\!\!\!\!\!\!\!\!\!\!\!\!\!\!\!\!\!\!\!\!\!\!\!\!\!\!\!\!\!\!\!\!\!\!\!\!\!\!\!\!\!\!\!\!\!\!\!\!\!\!\!\!\!\!\!\!\!\!\!\!\!\!\!\!\!\!\!\!\!\!\!\!\!\!\!\!\!\!\!\!\phantom{{\max\!\left\{1 , {\nb{\frac{1}{\frac{\partial P_{i_r}}{\partial X_{i_r}}(e_1 , \dots , e_g , e_{i_r})}}}_v\right\}}^{\!\!2}} \prod_{j = g + 1}^{N} \!\!\!{H_v({\widetilde P}_j)}^2 . {H_v(\underline e)}^{2 (d(G) - 1)} \right.\\
&\quad \left.\times \!\!\!\prod_{j = g + 1}^{N}
\!\!\!{\max\!\left\{1 , {\nb{\frac{1}{\frac{\partial P_j}{\partial
X_j}(e_1 , \dots , e_g ,
e_j)}}}_v\right\}}^{\!\!2}\right)^{\!\!\nb{I}} \!\!\!.
H_v(\underline e)
\end{split} \label{1.24}
\end{equation}
si $v$ est finie.
\end{lemma}
{\bf Démonstration.---} \'Etant donné une place $v$ de $K$ et $I
\in {\mathbb N}^g$, pour $i = 0 , \dots , g$, les estimations du
lemme \ref{a.6} sont triviales et pour $i = g + 1 , \dots , N$,
ces derni\`eres resultent de l'application du lemme \ref{a.5} pour
$P = P_i$, $n = g$, $(Y_1 , \dots , Y_n) = (X_1 , \dots , X_g)$,
$T = X_i$ et $\underline x = (e_1 , \dots , e_g , e_i)$. Il suffit
juste de remarquer que l'entier $d$ du lemme \ref{a.5} vaut dans
cette application:
$$d := {d°}_{\rm{tot}}P_i \leq {d°}_{\rm{tot}}{\widetilde{P}}_i \leq d(G) ,$$
que $H_v(P_i) = H_v({\widetilde{P}}_i) \leq \prod_{j = g + 1}^{N}
H_v({\widetilde{P}}_j)$ et que $H_v(1 , e_1 , \dots , e_g , e_i)
\leq H_v(\underline e)$. Le lemme \ref{a.6} s'ensuit.
$~~~~\blacksquare$\vspace{1mm}
\section{{\bf UN LEMME DE Z\'EROS DANS LES GROUPES ALG\'EBRIQUES COMMUTATIFS SOUS UNE PARAM\'ETRISATION
NON EXPLICITE}} Soit $G$ un groupe algébrique commutatif de
dimension $g$ $(g \geq 1)$, défini sur un corps de nombres $K$ et
plongé dans un espace projectif ${\mathbb P}_N$ $(N \geq 1)$.
Soient aussi $\mathbf e$ l'élément neutre de $G$ représenté dans
${\mathbb P}_N$ par un syst\`eme de coordonnées projectives
$\underline e = (e_0 : e_1 : \dots : e_N)$ et $\mathbf x$ un point
fixé de $G$. Désignons par $K[\underline X]$ ($\underline X = (X_0
, \dots , X_N)$) l'anneau des coordonnées de ${\mathbb P}_N$, par
$\mathfrak{G}$ l'idéal de définition de l'adhérence de Zariski
$\Gbar$ de $G$ dans $K[\underline X]$ et par $A$ l'anneau des
coordonnées de $\Gbar$:
$$A := K[\underline X] / \mathfrak{G} .$$
\subsection{Hypoth\`eses sur la paramétrisation du groupe $G$ et estimations de base:}
Nous supposons qu'il existe un entier strictement positif $p$, des
entiers strictement positifs $n_1 , \dots , n_p , g_1 , \dots ,
g_p , \delta_1 , \dots , \delta_p$, un réel strictement positif
$s$, des syst\`emes de coordonnées projectives ${\underline e}_1 ,
\dots , {\underline e}_p$ ne dépendants que de $\underline e$ et
appartenant aux espaces projectifs ${\mathbb P}_{n_1}(K) , \dots ,
{\mathbb P}_{n_p}(K)$ respectivement, des fonctions additives $t_1
, \dots , t_p$ de ${\mathbb N}^g$ dans ${\mathbb R}^+$
satisfaisant pour tout $I \in {\mathbb N}^g$:
$$t_1(I) + \dots + t_p(I) = \nb{I}$$
et un syst\`eme $\underline T = (T_1 , \dots , T_g)$ de variables
paramétrisant $G$ au voisinage de l'origine $\mathbf e$, dont la
paramétrisation induit un monomorphisme:
$$\varphi\: A \longrightarrow K[[\underline T]]$$
tel que pour tout $i = 0 , \dots , N$, l'image $\varphi(X_i)$ de
la coordonnée $X_i$ de $A$ s'exprime comme une série en
$\underline T$ à coefficients dans $K$:
$$\varphi(X_i) = \sum_{I \in {\mathbb N}^g} a_{I}^{(i)} {\underline T}^I$$
avec les $a_{I}^{(i)}$ $(I \in {\mathbb N}^g , i = 0 , \dots , N)$
satisfaisant pour tout $I \in {\mathbb N}^g$, tout $i = 0 , \dots
, N$ et toute place $v$ de $K$:
\begin{equation}
{\nb{a_{I}^{(i)}}}_v \leq s_{v}^{g_1 \delta_1 + \dots + g_p
\delta_p} E_{1 , v}^{t_1(I)} \dots E_{p , v}^{t_p(I)}
{H_v({\underline e}_1)}^{\delta_1} \dots {H_v({\underline
e}_p)}^{\delta_p} ,\label{1.25}
\end{equation}
o\`u
$$s_v := \begin{cases}
1 & \text{si $v$ est finie} \\
s & \text{si $v$ est infinie}
\end{cases}$$
et les $E_{\ell , v} \geq 1$ $(\ell = 1 , \dots , p)$ sont des
expressions suppérieures ou égales à $1$ ne dépendants que de
$n_{\ell} , G$ et ${\underline e}_{\ell}$.

Posons pour tout $\ell = 1 , \dots , p$:
$$E_{\ell} := \prod_{v \in M_K} E_{\ell , v}^{\frac{[K_v : {\mathbb Q}_v]}{[K : \mathbb Q]}} .$$

\'Etant donné $k \in {\mathbb N}^*$ et $(m_1 , \dots , m_p) \in
{{\mathbb N}^*}^p$, en utilisant (\ref{1.25}), on estime dans la
proposition qui suit les hauteurs locales ainsi que la hauteur de
Gauss-Weil de la famille de nombres de $K$ constituée des produits
$a_{I_1}^{(i_1)} \dots a_{I_k}^{(i_k)}$ pour $0 \leq i_1 \leq N ,
\dots , 0 \leq i_k \leq N$ et $t_1(I_1 + \dots + I_k) \leq m_1 ,
\dots , t_p(I_1 + \dots + I_k) \leq m_p$.
\begin{proposition}\label{a.8}
Sous toutes les hypoth\`eses précédentes, on a pour tout $k \in
{\mathbb N}^*$, tout $(m_1 , \dots , m_p) \in {{\mathbb N}^*}^p$
et toute place $v$ de $K$:
\begin{equation*}
\begin{split}
H_v\!\left(a_{I_1}^{(i_1)} \dots a_{I_k}^{(i_k)} ;~ 0 \leq i_1 \leq N , \dots , 0 \leq i_k \leq N ,\right.& \\
&\!\!\!\!\!\!\!\!\!\!\!\!\!\!\!\!\!\!\!\!\!\!\!\!\!\!\!\!\!\!\!\!\!\!\!\!\!\!\!\!\!\!\!\!\!\!\!\!\!\!\!\!\!\!\!\!\!\!\!\!\!\!\!\!\!\!\!\!\!\!\!\!\!\!\!\!\!\!\!\!\!\!\!\!\!\!\!\!\!\!\!\!\!\!\!\!\!\!\!\!\!\!\!\!\!\!\!\!\!\!\!\!\!\!\!\!\!\!\!\!\!\!\!\!\!\!\!\!\!\!\!\!\!\!\!\!\!\!\!\!\!\!\!\!\!\!\!\left.\phantom{a_{I_1}^{(i_1)} \dots a_{I_k}^{(i_k)} ;~ 0 \leq i_1 \leq N , \dots , 0 \leq i_k \leq N ,}t_1(I_1 + \dots + I_k) \leq m_1 , \dots , t_p(I_1 + \dots + I_k) \leq m_p\right) \\
&\!\!\!\!\!\!\!\!\!\!\!\!\!\!\!\!\!\!\!\!\!\!\!\!\!\!\!\!\!\!\!\!\!\!\!\!\!\!\!\!\!\!\!\!\!\!\!\!\!\leq
s_{v}^{k (g_1 \delta_1 + \dots + g_p \delta_p)} E_{1 , v}^{m_1}
\dots E_{p , v}^{m_p} {H_v({\underline e}_1)}^{k \delta_1} \dots
{H_v({\underline e}_p)}^{k \delta_p} .
\end{split}
\end{equation*}
Par conséquent, pour tout $k \in {\mathbb N}^*$ et tout $(m_1 ,
\dots , m_p) \in {{\mathbb N}^*}^p$, on a:
\begin{equation*}
\begin{split}
h\!\left(a_{I_1}^{(i_1)} \dots a_{I_k}^{(i_k)} ;~ 0 \leq i_1 \leq N , \dots , 0 \leq i_k \leq N ,\right.& \\
&\!\!\!\!\!\!\!\!\!\!\!\!\!\!\!\!\!\!\!\!\!\!\!\!\!\!\!\!\!\!\!\!\!\!\!\!\!\!\!\!\!\!\!\!\!\!\!\!\!\!\!\!\!\!\!\!\!\!\!\!\!\!\!\!\!\!\!\!\!\!\!\!\!\!\!\!\!\!\!\!\!\!\!\!\!\!\!\!\!\!\!\!\!\!\!\!\!\!\!\!\!\!\!\!\!\!\!\!\!\!\!\!\!\!\!\!\!\!\!\!\!\!\!\!\!\!\!\!\!\!\!\!\!\!\!\!\!\!\!\!\!\!\!\!\!\!\!\left.\phantom{a_{I_1}^{(i_1)} \dots a_{I_k}^{(i_k)} ;~ 0 \leq i_1 \leq N , \dots , 0 \leq i_k \leq N ,}t_1(I_1 + \dots + I_k) \leq m_1 , \dots , t_p(I_1 + \dots + I_k) \leq m_p\right) \\
&\!\!\!\!\!\!\!\!\!\!\!\!\!\!\!\!\!\!\!\!\!\!\!\!\!\!\!\!\!\!\!\!\leq
\sum_{\ell = 1}^{p} \log{(E_{\ell})} . m_{\ell} + k \sum_{\ell =
1}^{p} \delta_{\ell} \!\left(h({\underline e}_{\ell}) + g_{\ell}
\log{s}\right) .
\end{split}
\end{equation*}
\end{proposition}
{\bf Démonstration.---} Soient $k \in {\mathbb N}^*$ et $(m_1 ,
\dots , m_p) \in {{\mathbb N}^*}^p$. Pour tous $0 \leq i_1 \leq N
, \dots , 0 \leq i_k \leq N$, tous $I_1 , \dots , I_k$ dans
${\mathbb N}^g$ tels que $t_1(I_1 + \dots + I_k) \leq m_1 , \dots
, t_p(I_1 + \dots + I_k) \leq m_p$ et toute place $v$ de $K$, on a
d'apr\`es les relations (\ref{1.25}) et en tenant compte de
l'additivité des fonctions $t_{\ell}$ $(1 \leq \ell \leq p)$:
\begin{equation*}
\begin{split}
{\nb{a_{I_1}^{(i_1)} \dots a_{I_k}^{(i_k)}}}_v &\leq s_{v}^{k (g_1 \delta_1 + \dots + g_p \delta_p)} E_{1 , v}^{t_1(I_1 + \dots + I_k)} \!\dots E_{p , v}^{t_p(I_1 + \dots + I_k)} {H_v({\underline e}_1)}^{k \delta_1} \!\dots {H_v({\underline e}_p)}^{k \delta_p} \\
&\leq s_{v}^{k (g_1 \delta_1 + \dots + g_p \delta_p)} E_{1 ,
v}^{m_1} \dots E_{p , v}^{m_p} {H_v({\underline e}_1)}^{k
\delta_1} \dots {H_v({\underline e}_p)}^{k \delta_p} .
\end{split}
\end{equation*}
Ce qui entra{\sf\^\i}ne l'estimation de la proposition \ref{a.8}
pour toute hauteur locale $v$-adique de la famille
$\{a_{I_1}^{(i_1)} \dots a_{I_k}^{(i_k)} ; 0 \leq i_1 \leq N ,
\dots , 0 \leq i_k \leq N , t_1(I_1 + \dots + I_k) \leq m_1 ,
\dots , t_p(I_1 + \dots + I_k) \leq m_p\}$. On n'a qu'à reporter
les estimations locales dans la définition de la hauteur de
Gauss-Weil pour conclure. $~~~~\blacksquare$\vspace{1mm}
\begin{corollary}\label{a.15}
Pour tout $\underline{i} \in {\mathbb N}^{N + 1}$, on a:
$$\varphi\left({\underline X}^{\underline i}\right) = \sum_{I \in {\mathbb N}^g} \mathcal{C}(\underline i , I) {\underline T}^I ,$$
o\`u les $\mathcal{C}(\underline i , I)$ $(I \in {\mathbb N}^g)$
sont des nombres de $K$ satisfaisant pour tout $k \in {\mathbb
N}^*$, tout $(m_1 , \dots , m_p) \in {{\mathbb N}^*}^p$ et toute
place $v$ de $K$:
\begin{equation*}
\begin{split}
H_v\!\left(\mathcal{C}(\underline i , I) ; \nb{\underline i} = k , t_1(I) \leq m_1 , \dots , t_p(I) \leq m_p\right) \\
&\!\!\!\!\!\!\!\!\!\!\!\!\!\!\!\!\!\!\!\!\!\!\!\!\!\!\!\!\!\!\!\!\!\!\!\!\!\!\!\!\!\!\!\!\!\!\!\!\!\!\!\!\!\!\!\!\!\!\!\!\!\!\!\!\!\!\!\!\!\!\!\!\!\!\!\!\!\!\!\!\!\!\!\!\!\!\leq
2^{m_1 + \dots + m_p + g (k - 1)} s^{k (g_1 \delta_1 + \dots + g_p
\delta_p)} E_{1 , v}^{m_1} \dots E_{p , v}^{m_p} {H_v({\underline
e}_1)}^{k \delta_1} \dots {H_v({\underline e}_p)}^{k \delta_p}
\end{split}
\end{equation*}
si $v$ est infinie et:
\begin{equation*}
\begin{split}
H_v\!\left(\mathcal{C}(\underline i , I) ; \nb{\underline i} = k , t_1(I) \leq m_1 , \dots , t_p(I) \leq m_p\right) \\
&\!\!\!\!\!\!\!\!\!\leq E_{1 , v}^{m_1} \dots E_{p , v}^{m_p}
{H_v({\underline e}_1)}^{k \delta_1} \dots {H_v({\underline
e}_p)}^{k \delta_p}
\end{split}
\end{equation*}
si $v$ est finie. \\
Par conséquent, pour tout $k \in {\mathbb N}^*$ et tout $(m_1 ,
\dots , m_p) \in {{\mathbb N}^*}^p$, on a:
\begin{equation*}
\begin{split}
h\!\left(\mathcal{C}(\underline i , I) ; \nb{\underline i} = k , t_1(I) \leq m_1 , \dots , t_p(I) \leq m_p\right) \\
&\!\!\!\!\!\!\!\!\!\!\!\!\!\!\!\!\!\!\!\!\!\!\!\!\!\!\!\!\!\!\!\!\!\!\!\!\!\!\!\!\!\!\!\!\!\!\!\!\!\!\!\!\!\!\!\!\!\!\!\!\!\!\!\!\!\!\!\!\!\!\!\!\!\!\!\!\!\!\!\!\!\leq
\sum_{\ell = 1}^{p} \left(\log{(E_{\ell})} + \log{2}\right)
m_{\ell} + k \!\left(\!g \log{2} + \sum_{\ell = 1}^{p}
\delta_{\ell} \!\left(h({\underline e}_{\ell}) + g_{\ell}
\log{s}\right)\!\right) .
\end{split}
\end{equation*}
\end{corollary}
{\bf Démonstration.---} Pour tout $\underline{i} = (i_0 , \dots ,
i_N) \in {\mathbb N}^{N + 1}$, on a:
\begin{equation*}
\begin{split}
\varphi({\underline X}^{\underline i}) &= {\varphi(X_0)}^{i_0} {\varphi(X_1)}^{i_1} \dots {\varphi(X_N)}^{i_N} \\
&= \left(\sum_{I \in {\mathbb N}^g} a_{I}^{(0)} {\underline
T}^I\right)^{\!\!i_0} \dots \left(\sum_{I \in {\mathbb N}^g}
a_{I}^{(N)} {\underline T}^I\right)^{\!\!i_N} ,
\end{split}
\end{equation*}
c'est-à-dire:
$$\varphi({\underline X}^{\underline i}) = \sum_{I \in {\mathbb N}^g} \left(\sum_{I_{0 , 1} + \dots + I_{0 , i_0} + \dots + I_{N , 1} + \dots + I_{N , i_N} = I} a_{I_{0 , 1}}^{(0)} \dots a_{I_{0 , i_0}}^{(0)} \dots a_{I_{N , 1}}^{(N)} \dots a_{I_{N , i_N}}^{(N)}\right) {\underline T}^I .$$
Les $\mathcal{C}(\underline i , I)$ $(\underline i \in {\mathbb
N}^{N + 1} , I \in {\mathbb N}^g)$ sont alors les nombres de $K$
définis par les relations:
\begin{equation}
\mathcal{C}(\underline i , I) = \!\!\!\!\!\!\sum_{I_{0 , 1} +
\dots + I_{0 , i_0} + \dots + I_{N , 1} + \dots + I_{N , i_N} = I}
\!\!\!\!\!\!\!\!\!\!\!\!a_{I_{0 , 1}}^{(0)} \dots a_{I_{0 ,
i_0}}^{(0)} \dots a_{I_{N , 1}}^{(N)} \dots a_{I_{N , i_N}}^{(N)}
~~(\underline i \in {\mathbb N}^{N + 1} , I \in {\mathbb N}^g) .
\label{1.26}
\end{equation}
Pour tout $k \in {\mathbb N}^*$, tout $(m_1 , \dots , m_p) \in
{{\mathbb N}^*}^p$ et toute place $v$ de $K$, de (\ref{1.26})
découle immédiatement les estimations:
\begin{equation*}
\begin{split}
H_v\!\left(\mathcal{C}(\underline i , I) ; \nb{\underline i} = k , t_1(I) \leq m_1 , \dots , t_p(I) \leq m_p\right) \\
&\!\!\!\!\!\!\!\!\!\!\!\!\!\!\!\!\!\!\!\!\!\!\!\!\!\!\!\!\!\!\!\!\!\!\!\!\!\!\!\!\!\!\!\!\!\!\!\!\!\!\!\!\!\!\!\!\!\!\!\!\!\!\!\!\!\!\!\!\!\!\!\!\!\!\!\!\!\!\!\!\!\!\!\!\!\!\!\!\!\!\!\!\!\!\!\!\!\!\!\!\!\leq H_v\!\left(a_{I_1}^{(i_1)} \dots a_{I_k}^{(i_k)} ;~ 0 \leq i_1 \leq N , \dots , 0 \leq i_k \leq N ,\right. \\
&\!\!\!\!\!\!\!\!\!\!\!\!\!\!\!\!\!\!\!\!\!\!\!\!\!\!\!\!\!\!\!\!\!\!\!\!\!\!\!\!\!\!\!\!\!\!\!\!\!\!\!\!\!\!\!\!\!\!\!\!\!\!\!\!\!\!\!\!\!\!\!\!\!\!\!\!\!\!\!\!\!\!\!\!\!\!\!\!\!\!\!\!\!\!\!\!\!\!\!\!\!\!\!\!\!\!\!\!\!\!\!\!\!\!\!\!\!\!\!\!\!\!\!\!\!\!\!\!\!\!\!\!\!\!\!\!\!\!\!\!\!\!\!\!\!\left.\phantom{a_{I_1}^{(i_1)}
\dots a_{I_k}^{(i_k)} ;~ 0 \leq i_1 \leq N , \dots , 0 \leq i_k
\leq N ,}t_1(I_1 + \dots + I_k) \leq m_1 , \dots , t_p(I_1 + \dots
+ I_k) \leq m_p\right)
\end{split}
\end{equation*}
si $v$ est finie,
\begin{equation*}
\begin{split}
H_v\!\left(\mathcal{C}(\underline i , I) ; \nb{\underline i} = k , t_1(I) \leq m_1 , \dots , t_p(I) \leq m_p\right) \\
&\!\!\!\!\!\!\!\!\!\!\!\!\!\!\!\!\!\!\!\!\!\!\!\!\!\!\!\!\!\!\!\!\!\!\!\!\!\!\!\!\!\!\!\!\!\!\!\!\!\!\!\!\!\!\!\!\!\!\!\!\!\!\!\!\!\!\!\!\!\!\!\!\!\!\!\!\!\!\!\!\!\!\!\!\!\!\!\!\!\!\!\!\!\!\!\!\!\!\!\!\!\leq R . H_v\!\left(a_{I_1}^{(i_1)} \dots a_{I_k}^{(i_k)} ;~ 0 \leq i_1 \leq N , \dots , 0 \leq i_k \leq N ,\right. \\
&\!\!\!\!\!\!\!\!\!\!\!\!\!\!\!\!\!\!\!\!\!\!\!\!\!\!\!\!\!\!\!\!\!\!\!\!\!\!\!\!\!\!\!\!\!\!\!\!\!\!\!\!\!\!\!\!\!\!\!\!\!\!\!\!\!\!\!\!\!\!\!\!\!\!\!\!\!\!\!\!\!\!\!\!\!\!\!\!\!\!\!\!\!\!\!\!\!\!\!\!\!\!\!\!\!\!\!\!\!\!\!\!\!\!\!\!\!\!\!\!\!\!\!\!\!\!\!\!\!\!\!\!\!\!\!\!\!\!\!\!\!\!\!\!\!\left.\phantom{a_{I_1}^{(i_1)}
\dots a_{I_k}^{(i_k)} ;~ 0 \leq i_1 \leq N , \dots , 0 \leq i_k
\leq N ,}t_1(I_1 + \dots + I_k) \leq m_1 , \dots , t_p(I_1 + \dots
+ I_k) \leq m_p\right)
\end{split}
\end{equation*}
si $v$ est infinie, avec:
$$R := \max_{I \in {\mathbb N}^g , t_1(I) \leq m_1 , \dots , t_p(I) \leq m_p} {\rm{card}} \left\{(I_1 , \dots , I_k) \in ({\mathbb N}^g)^k / I_1 + \dots + I_k = I\right\} .$$
De plus,
\begin{equation*}
\begin{split}
h\!\left(\mathcal{C}(\underline i , I) ; \nb{\underline i} = k , t_1(I) \leq m_1 , \dots , t_p(I) \leq m_p\right) \\
&\!\!\!\!\!\!\!\!\!\!\!\!\!\!\!\!\!\!\!\!\!\!\!\!\!\!\!\!\!\!\!\!\!\!\!\!\!\!\!\!\!\!\!\!\!\!\!\!\!\!\!\!\!\!\!\!\!\!\!\!\!\!\!\!\!\!\!\!\!\!\!\!\!\!\!\!\!\!\!\!\!\!\!\!\!\!\!\!\!\!\!\!\!\!\!\!\!\leq h\!\left(a_{I_1}^{(i_1)} \dots a_{I_k}^{(i_k)} ;~ 0 \leq i_1 \leq N , \dots , 0 \leq i_k \leq N ,\right. \\
&\!\!\!\!\!\!\!\!\!\!\!\!\!\!\!\!\!\!\!\!\!\!\!\!\!\!\!\!\!\!\!\!\!\!\!\!\!\!\!\!\!\!\!\!\!\!\!\!\!\!\!\!\!\!\!\!\!\!\!\!\!\!\!\!\!\!\!\!\!\!\!\!\!\!\!\!\!\!\!\!\!\!\!\!\!\!\!\!\!\!\!\!\!\!\!\!\!\!\!\!\!\!\!\!\!\!\!\!\!\!\!\!\!\!\!\!\!\!\!\!\!\!\!\!\!\!\!\!\!\!\!\!\!\!\!\!\!\!\!\!\!\!\!\!\!\!\!\!\!\!\!\!\!\!\!\!\!\!\!\!\left.\phantom{a_{I_1}^{(i_1)}
\dots a_{I_k}^{(i_k)} ;~ 0 \leq i_1 \leq N , \dots , 0 \leq i_k
\leq N ,}t_1(I_1 + \dots + I_k) \leq m_1 , \dots , t_p(I_1 + \dots
+ I_k) \leq m_p\right) + \log{R}
\end{split}
\end{equation*}
(pour le m\^eme $R$). \\
le corollaire \ref{a.15} suit alors de l'application de la
proposition \ref{a.8} et de la majoration de $R$ par:
\begin{equation*}
\begin{split}
R &= \max_{\begin{array}{c}
\scriptstyle{I = ({\mathcal I}_1 , \dots , {\mathcal I}_g) \in {\mathbb N}^g} \\
\scriptstyle{t_1(I) \leq m_1 , \dots , t_p(I) \leq m_p}
\end{array}
}
\!\!\!\binom{{\mathcal I}_1 + k - 1}{k - 1} \dots \binom{{\mathcal I}_g + k - 1}{k - 1} \\
&\leq \max_{\begin{array}{c}
\scriptstyle{I = ({\mathcal I}_1 , \dots , {\mathcal I}_g) \in {\mathbb N}^g} \\
\scriptstyle{t_1(I) \leq m_1 , \dots , t_p(I) \leq m_p}
\end{array}
} \!\!\!2^{{\mathcal I}_1 + k - 1} \dots 2^{{\mathcal I}_g + k - 1} \\
&\leq \max_{I \in {\mathbb N}^g , t_1(I) \leq m_1 , \dots , t_p(I) \leq m_p} 2^{\nb{I} + g (k - 1)} \\
&\leq 2^{m_1 + \dots + m_p + g (k - 1)} ~~~~~~~~~~\text{(car
$\nb{I} = t_1(I) + \dots + t_p(I)$).}
\end{split}
\end{equation*}
La démonstration est achevée. $~~~~\blacksquare$\vspace{1mm}
\subsection{Opérateurs:}
Fixons pour tout point $\mathbf y$ de $G$ une famille ${\underline
A}_{\mathbf y} = (A_{\mathbf y , 0} , \dots , A_{\mathbf y , N})$
de formes de $k[\underline X , \underline Y]$ (avec $\underline X
:= (X_0 , \dots , X_N)$ et $\underline Y := (Y_0 , \dots , Y_N)$)
représentant l'addition dans $G \hookrightarrow {\mathbb P}_N$ au
voisinage de $\{\mathbf e\} \times \{\mathbf y\}$.

Pour toute forme $P$ de $A$, tout point $\mathbf y$ de $G$ et tout
$g$-uplet $I$ dans ${\mathbb N}^g$, on pose $\Delta_{\mathbf
y}^{I} P$ le coefficient dans $A$ du mon\^ome ${\underline T}^I$
de la série $(\varphi \otimes {\rm{id}}) (P({\underline
A}_{\mathbf y}))$ de $K[\underline Y][[\underline Z]]$ qu'on
obtient en prenant la forme $P({\underline A}_{\mathbf y})$ de
$K[\underline X , \underline Y]$, définie naturellement par:
$$P\!\left({\underline A}_{\mathbf y}\right) \!\left(\underline X , \underline Y\right) := P\!\left(A_{\mathbf y , 0}\left(\underline X , \underline Y\right) , \dots , A_{\mathbf y , N}\left(\underline X , \underline Y\right)\right) ,$$
et en lui appliquant $\varphi \otimes \rm{id}$, qui revient
simplement à substituer dans $P({\underline A}_{\mathbf y})$ aux
variables $X_0 , \dots , X_N$ leurs images par $\varphi$. La série
$(\varphi \otimes \rm{id}) (P({\underline A}_{\mathbf y}))$
s'écrit:
$$\left(\varphi \otimes \rm{id}\right) \!\left(P\!\left({\underline A}_{\mathbf y}\right)\right) \!\left(\underline Y\right) = \sum_{I \in {\mathbb N}^g} \left(\Delta_{\mathbf y}^{I} P\right) \!\left(\underline Y\right) . {\underline T}^I .$$
On utilise la notion de {\it{dessous d'escalier}} définie dans
\cite{Ph3}, page 1071 et pour un dessous d'escalier fini $W$ de
${\mathbb N}^g$, on désigne par $t_i(W)$ $(i = 1 , \dots , p)$ les
quantités positives:
$$t_i(W) := \max_{I \in W} t_i(I) ~~~~ (i = 1 , \dots , p) .$$
Pour toute sous-variété $V$ de $G$ et tout dessous d'escalier fini $W$ de ${\mathbb N}^g$,
on note aussi par $m_W (V)$ la quantité de \cite{Ph3} définie à la page 1072. \vspace{0.2cm}\\
{\bf N.B.}---\label{A} Bien que les formes $\Delta_{\mathbf
y}^{I}P$ $(I \in {\mathbb N}^g)$ dépendent du choix de
${\underline A}_{\mathbf y}$, un sous-ensemble algébrique de $G$
défini par une famille de formes $\Delta_{\mathbf y}^{I}P$, pour
$I$ décrivant un certain dessous d'escalier $W$ de ${\mathbb
N}^g$, est indépendant de ce choix. Ceci justifie les définitions
suivantes:
\begin{definitions}\label{a.14}~\vspace{0.1cm}
\begin{description}
\item[1)] On dit qu'une forme $P$ de $A$ s'annule en un point
$\mathbf y$ de $G$ avec une multiplicité définie par un dessous
d'escalier $W$ de ${\mathbb N}^g$, si pour tout $I \in W$ on a:
$$\left(\Delta_{\mathbf y}^{I} P\right)(\underline y) = 0 ,$$
avec $\underline y$ un syst\`eme de coordonnées projectives
représentant $\mathbf y \in {\mathbb P}_N$. \item[2)] On dit
qu'une forme $P$ de $A$ s'annule sur un sous-ensemble algébrique
$E$ de $G$ avec une multiplicité définie par un dessous d'escalier
$W$ de ${\mathbb N}^g$, si $P$ s'annule en tout point $\mathbf y$
de $E$ avec la multiplicité définie par $W ~\!(\text{en se
référant à} ~\!1))$.
\end{description}
\end{definitions}
On a le lemme suivant:
\begin{lemma}\label{a.11}
Soient $P$ une forme de $A$ et $V$ une sous-variété de $G$ de
dimension $\geq 1$. On a:
\begin{description}
\item[1)] Si $P$ s'annule sur presque\footnote{Le mot presque veut
dire ici: ``à l'exception d'un nombre fini de diviseurs''.} toute
$V$, alors $P$ s'annule sur $V$. \item[2)] Si $P$ s'annule sur
presque toute $V$ avec une multiplicité définie par un dessous
d'escalier $W$ de ${\mathbb N}^g$, alors $P$ s'annule sur $V$ avec
la multiplicité définie par $W$.
\end{description}
\end{lemma}
{\bf Démonstration.---} Démontrons $1)$: supposons que $P$
s'annule sur presque toute $V$ et soit $X$ le sous-ensemble de $V$
constitué de tous les points d'exception. En appelant $H$
l'hypersurface de $G$ définie par $P$, on a donc:
$$V \subset H \cup X .$$
Comme $H \cup X$ est lui aussi un sous-ensemble algébrique de $G$ et $V$ est irréductible (car c'est une sous-variété) alors $V$ doit \^etre contenue dans l'une au moins des composantes irréductibles de $H \cup X$; or chacune de ces derni\`eres est soit une composante irréductible de $H$, soit une composante de dimension $< X$. Mais, $V$ ne peut pas \^etre contenue dans une composante irréductible de dimension $< \dim(X)$. Donc, $V$ est forcément contenue dans une composante irréductible de $H$; par conséquent $V \subset H$, c'est-à-dire que $P$ s'annule sur $V$. Ceci démontre $1)$. \\
Démontrons $2)$: Supposons que $P$ s'annule sur presque toute $V$ avec une multiplicité définie par un dessous d'escalier fini $W$ de ${\mathbb N}^g$ et soit $Y$ le sous-ensemble de $V$ constitué de tous les points d'exceptions. \\
Soient $\mathbf y$ un point quelconque de $V$ représenté dans
${\mathbb P}_N$ par $\underline y$ et $I$ un $g$-uplet quelconque
de $W$. On peut trouver $\Omega_{\mathbf e}$ un voisinage ouvert
de Zariski dans $G$ du point $\mathbf e$ et $\Omega_{\mathbf y}$
un voisinage ouvert de Zariski dans $G$ du point $\mathbf y$ de
façon à ce que la représentation de l'addition dans $G$ au
voisinage de $\{\mathbf e\} \times \{\mathbf y\}$ par la famille
de formes ${\underline A}_{\mathbf y}$ soit valable sur tout le
produit $\Omega_{\mathbf e} \times \Omega_{\mathbf y}$. Ainsi,
d'apr\`es la note de la page \pageref{A}, l'hypoth\`ese de $2)$
entra{\sf\^\i}ne que la forme $\Delta_{\mathbf y}^{I} P$ s'annule
en tout point de $(V \cap \Omega_{\mathbf y}) \setminus Y$, donc
s'annule sur presque toute $V$ (puisque la différence ensembliste
entre $V$ et $(V \cap \Omega_{\mathbf y}) \setminus Y$ est
contenue dans un diviseur); d'o\`u, d'apr\`es $1)$ (déja
démontré), $\Delta_{\mathbf y}^{I} P$ doit s'annuler sur $V$ toute
enti\`ere et en particulier elle doit s'annuler en $\mathbf y$,
c'est-à-dire qu'on a: $(\Delta_{\mathbf y}^{I} P) (\underline y) =
0$. Comme ceci est vrai pour tout $\mathbf y \in V$ et tout $I \in
W$, cela revient à dire (par définition m\^eme) que $P$ s'annule
sur $V$ avec la multiplicité définie par $W$. Ce qui établit $2)$
et ach\`eve cette démonstration. $~~~~\blacksquare$\vspace{1mm}

Pour tous $P , Q$ dans $A$, tout point $\mathbf y$ de $G$ et tout
$I$ dans ${\mathbb N}^g$, on vérifie immédiatement les identités
suivantes:
\begin{equation}
\begin{split}
\Delta_{\mathbf y}^{I}(P + Q) &= \Delta_{\mathbf y}^{I} P + \Delta_{\mathbf y}^{I} Q \\
\Delta_{\mathbf y}^{I}(P.Q) &= \sum_{I_1 + I_2 = I}
\Delta_{\mathbf y}^{I_1} P . \Delta_{\mathbf y}^{I_2} Q .
\end{split} \label{1.22}
\end{equation}
Une autre propriété intéressante de ces opérateurs
$\Delta_{\mathbf y}^{I}$ $(\mathbf y \in G , I \in {\mathbb
N}^g)$, qui suit principalement de l'associativité de la loi
d'addition dans $G$, est donnée par le lemme suivant:
\begin{lemma}\label{a.12}
Pour toute forme $P$ de $A$, tout point $\mathbf y$ de $G$ et tous
$I , J$ dans ${\mathbb N}^g$; la forme $\Delta_{\mathbf
y}^{I}(\Delta_{\mathbf y}^{J} P)$ appartient à l'idéal engendré
par la famille de formes $\Delta_{\mathbf y}^{L}P$; $L \leq I + J$
(o\`u $L \leq I + J$ veut dire ici que chaque composante de $L$
est inférieure ou égale à la composante correspendante de $I +
J$).
\end{lemma}
\begin{definition}
Plus généralement, pour tout dessous d'escalier $W$ de ${\mathbb
N}^g$, tout point $\mathbf y$ de $G$ et tout idéal homog\`ene
$\mathfrak I$ de $A$, on définit $\Delta_{\mathbf y}^{W} \mathfrak
I$ comme l'idéal homog\`ene de $A$ engendré par la famille de
formes:
$$\Delta_{\mathbf y}^{I} P ; ~~ I \in W , P \in \mathfrak I .$$
\end{definition}
On vérifie aisément, en utilisant les relations (\ref{1.22}), que
si $\mathfrak I$ est engendré par des formes $P_1 , \dots , P_s$
de $A$, alors $\Delta_{\mathbf y}^{W}\mathfrak I$ est engendré par
les familles de formes: $\Delta_{\mathbf y}^{I}P_1$ $(I \in W) ,
\dots , \Delta_{\mathbf y}^{I}P_s$ $(I \in W)$.

Afin d'alléger l'écriture, on note $\underline A = (A_0 , \dots ,
A_N)$ pour ${\underline A}_{\mathbf x}$, pour $I \in {\mathbb
N}^g$ et $P \in A$, on note $\Delta^I P$ pour $\Delta_{\mathbf
x}^{I} P$ et pour un dessous d'escalier $W$ de ${\mathbb N}^g$ et
un idéal $\mathfrak I$ de $A$, on note $\Delta^W \mathfrak I$ pour
$\Delta_{\mathbf x}^{W} \mathfrak I$. On fixe aussi ${\mathcal
O}_{\mathbf e}$ et ${\mathcal O}_{\mathbf x}$ deux voisinages
ouverts de Zariski dans $G$ des points $\mathbf e$ et $\mathbf x$
respectivement, tels que la représentation de l'addition dans $G$
au voisinage de $\{\mathbf e\} \times \{\mathbf x\}$ soit valable
sur tout le produit ${\mathcal O}_{\mathbf e} \times {\mathcal
O}_{\mathbf x}$. Notons finalement $(c , c')$ $(c , c'$ des
entiers $\geq 1)$ le bidegré des formes $A_0 , \dots , A_N$ de
$K[\underline X , \underline Y]$ constituant la famille
$\underline A$.

Etant donné $P$ une forme de $K[\underline X]$ de degré $\delta$
($\delta \in \mathbb N$), les $\Delta^I P$ $(I \in {\mathbb N}^g)$
sont -par définition- les coefficients de la série $(\varphi
\otimes {\rm{id}})(P(\underline A))(\underline Y)$ de
$K[\underline Y][[\underline Z]]$ et sont donc des formes de
$K[\underline Y]$ de degré $c' \delta$, puisque $P(\underline A)$
est une forme de $K[\underline X , \underline Y]$ de bidegré $(c
\delta , c' \delta)$. Pour tout $(m_1 , \dots , m_p) \in {{\mathbb
N}^*}^p$, on estime dans la proposition qui suit la hauteur de
Gauss-Weil de la famille de formes $\Delta^I P$, $t_1(I) \leq m_1
, \dots , t_p(I) \leq m_p$ en fonction de $\delta = {d°}P$,
$\widetilde{h}(P)$, $c$, $c'$, $\widetilde{h}(\underline A)$, $N$,
$g = \dim G$, $d(G)$, $h(G)$, $\mathbf e$, $g_1 , \dots , g_p$,
$\delta_1 , \dots , \delta_p$, $E_1 , \dots , E_p$, $m_1 , \dots ,
m_p$ et $s$.
\begin{proposition}\label{a.16}
Pour toute forme $P \in K[\underline X]$ de degré $\delta$
$(\delta \in \mathbb N)$ et tout $(m_1 , \dots , m_p)$ \\ $\in
{{\mathbb N}^*}^p$, on a:
\begin{equation*}
\begin{split}
\widetilde{h}\!\left(\Delta^I P ;~ t_1(I) \leq m_1 , \dots , t_p(I) \leq m_p\right) \leq \\
&\!\!\!\!\!\!\!\!\!\!\!\!\!\!\!\!\!\!\!\!\!\!\!\!\!\!\!\!\!\!\!\!\!\!\!\!\!\!\!\!\!\!\!\!\!\!\!\!\!\!\!\!\!\!\!\!\!\!\!\!\!\!\!\!\!\!\!\!\!\!\!\!\!\!\!\!\!\!\!\!\!\!\!\!\!\!\!\!\!\!\!\!\!\!\!\sum_{\ell = 1}^{p}\left(\log{(E_{\ell})} + \log{2}\right) m_{\ell} + c \delta \!\!\left(\!\!g \log{2} + \sum_{\ell = 1}^{p} \delta_{\ell}\!\left(h({\underline e}_{\ell}) + g_{\ell} \log{s}\!\right)\!\!\right) + \widetilde{h}(P) + \delta \widetilde{h}(\underline A) \\
&~~~~~~~~~~~~\!+ \delta (2 c + c' + 1) \log{(N + 1)} .
\end{split}
\end{equation*}
\end{proposition}
{\bf Démonstration.---} $P(\underline A)$ est une forme de
$K[\underline X , \underline Y]$ de bidegré $(c \delta , c'
\delta)$ et de hauteur de Gauss-Weil facilement estimée en
raisonnant place par place (et en utilisant les longueurs pour les
places infinies et les hauteurs pour les places finies) par:
\begin{equation}
\widetilde{h}\!\left(P(\underline A)\right) \leq \widetilde{h}(P)
+ \delta \widetilde{h}(\underline A) + \log\!\left({\!\binom{N +
\delta}{\delta} \!\!{\binom{N + c}{c}}^{\!\!\delta}} \!\!{\binom{N
+ c'}{c'}}^{\!\!\delta}\right) . \label{1.29}
\end{equation}
En écrivant:
$$P(\underline A)(\underline X , \underline Y) = \sum_{\nb{\underline i} = c \delta , \nb{\underline j} = c' \delta} \!\!\!\!\!\!\rho(\underline i , \underline j) {\underline X}^{\underline i} {\underline Y}^{\underline j}$$
et:
$$\varphi({\underline X}^{\underline i}) = \sum_{I \in {\mathbb N}^g} \mathcal{C}(\underline i , I) {\underline T}^I ~~~~~~ \text{(pour $\underline i \in {\mathbb N}^{N + 1}$)}$$
(comme dans le corollaire \ref{a.15}), on a:
\begin{equation*}
\begin{split}
(\varphi \otimes {\rm{id}})\left(P(\underline A)\right) &:= P(\underline A)\left(\varphi(\underline X) , \underline Y\right) \\
&~= \sum_{\nb{\underline i} = c \delta , \nb{\underline j} = c' \delta} \!\!\!\!\!\!\rho(\underline i , \underline j) \left(\sum_{I \in {\mathbb N}^g} \mathcal{C}(\underline i , I) {\underline T}^I\right) {\underline Y}^{\underline j} \\
&~= \sum_{I \in {\mathbb N}^g} \!\!\left(\sum_{\nb{\underline i} =
c \delta , \nb{\underline j} = c' \delta}
\!\!\!\!\!\!\rho(\underline i , \underline j)
\mathcal{C}(\underline i , I) {\underline Y}^{\underline j}\right)
{\underline T}^I
\end{split}
\end{equation*}
(en permutant les signes de sommations). D'o\`u, par définition
m\^eme des $\Delta^I P$ $(I \in {\mathbb N}^g)$, pour tout $I \in
{\mathbb N}^g$:
\begin{equation*}
\begin{split}
\Delta^I P(\underline Y) &= \sum_{\nb{\underline i} = c \delta , \nb{\underline j} = c' \delta} \!\!\!\!\!\!\rho(\underline i , \underline j) \mathcal{C}(\underline i , I) {\underline Y}^{\underline j} \\
&= \sum_{\nb{\underline j} = c' \delta}
\!\!\left(\sum_{\nb{\underline i} = c \delta} \rho(\underline i ,
\underline j) \mathcal{C}(\underline i , I)\right) {\underline
Y}^{\underline j} ,
\end{split}
\end{equation*}
et par suite, pour tout $(m_1 , \dots , m_p) \in {{\mathbb
N}^*}^p$:
\begin{equation}
\begin{split}
\widetilde{h}\left(\Delta^I P ;~ t_1(I) \leq m_1 , \dots , t_p(I) \leq m_p\right) \\
&\!\!\!\!\!\!\!\!\!\!\!\!\!\!\!\!\!\!\!\!\!\!\!\!\!\!\!\!\!\!\!\!\!\!\!\!\!\!\!\!\!\!\!\!\!\!\!\!\!= h\!\left(\sum_{\nb{\underline i} = c \delta} \rho(\underline i , \underline j) \mathcal{C}(\underline i , I) ;~ \nb{\underline j} = c' \delta , t_1(I) \leq m_1 , \dots , t_p(I) \leq m_p\right) \\
&\!\!\!\!\!\!\!\!\!\!\!\!\!\!\!\!\!\!\!\!\!\!\!\!\!\!\!\!\!\!\!\!\!\!\!\!\!\!\!\!\!\!\!\!\!\!\!\!\!\!\!\!\!\!\!\!\!\!\!\!\!\!\!\!\!\!\!\!\!\!\!\!\!\!\!\!\!\!\!\!\!\!\!\!\!\!\!\!\!\!\!\!\leq h\!\left(\rho(\underline i , \underline j) ;~ \nb{\underline i} = c \delta , \nb{\underline j} = c' \delta\right) + h\!\left(\mathcal{C}(\underline i , I) ;~ \nb{\underline i} = c \delta , t_1(I) \leq m_1 , \dots , t_p(I) \leq m_p\right) \\
&~~~~~~~~~~~~~~~~~~~~~~~~~~~\!\!+ \log{\binom{N + c \delta}{c
\delta}} \label{1.27}
\end{split}
\end{equation}
(gr\^ace à un calcul standard sur les hauteurs). Or
\begin{equation}
\begin{split}
h\!\left(\rho(\underline i , \underline j) ;~ \nb{\underline i} = c \delta , \nb{\underline j} = c' \delta\right) &= \widetilde{h}\!\left(P(\underline A)\right) \\
&\!\!\!\!\!\!\!\!\!\!\!\!\!\!\!\!\!\!\!\!\!\!\!\!\!\!\!\!\!\!\!\!\!\!\!\!\!\!\!\!\!\!\!\!\!\!\leq
\widetilde{h}(P) + \delta \widetilde{h}(\underline A) +
\log\!\left({\!\binom{N + \delta}{\delta} \!\!{\binom{N +
c}{c}}^{\!\!\delta}} \!\!{\binom{N +
c'}{c'}}^{\!\!\delta}\right)\text{(d'apr\`es (\ref{1.29}))} ,
\label{1.28}
\end{split}
\end{equation}
et:
\begin{equation}
\begin{split}
h\!\left(\mathcal{C}(\underline i , I) ; \nb{\underline i} = c \delta , t_1(I) \leq m_1 , \dots , t_p(I) \leq m_p\right) \\
&\!\!\!\!\!\!\!\!\!\!\!\!\!\!\!\!\!\!\!\!\!\!\!\!\!\!\!\!\!\!\!\!\!\!\!\!\!\!\!\!\!\!\!\!\!\!\!\!\!\!\!\!\!\!\!\!\!\!\!\!\!\!\!\!\!\!\!\!\!\!\!\!\!\!\!\!\!\!\!\!\!\!\!\!\leq
\sum_{\ell = 1}^{p} \left(\log{(E_{\ell})} + \log{2}\right)
m_{\ell} + c \delta \!\left(\!g \log{2} + \sum_{\ell = 1}^{p}
\delta_{\ell} \!\left(h({\underline e}_{\ell}) + g_{\ell}
\log{s}\right)\!\right) \label{1.30}
\end{split}
\end{equation}
(d'apr\`es le corollaire \ref{a.15}). \\
La proposition \ref{a.16} en résulte en reportant (\ref{1.28}) et
(\ref{1.30}) dans (\ref{1.27}) et en majorant le produit de
bin\^omes $\binom{N + c \delta}{c \delta} \!\binom{N +
\delta}{\delta} \!{\binom{N + c}{c}}^{\!\delta} \!{\binom{N +
c'}{c'}}^{\!\delta}$ par:
\begin{equation*}
\begin{split}
\binom{N + c \delta}{c \delta} \!\!\binom{N + \delta}{\delta} \!\!{\binom{N + c}{c}}^{\!\!\delta} \!\!{\binom{N + c'}{c'}}^{\!\!\delta} &\leq (N + 1)^{c \delta} (N + 1)^{\delta} (N + 1)^{c \delta} (N + 1)^{c' \delta} \\
&\leq (N + 1)^{\delta (2 c + c' + 1)} .
\end{split}
\end{equation*}
La démonstration est achevée. $~~~~\blacksquare$\vspace{1mm}
\subsection{Multipicités, degrés et hauteurs:}
Rappelons que, lorsque $A$ est un anneau noetherien, $\mathfrak I$
un idéal de $A$ contenu dans le radical de $A$ et $M$ un
$A$-module de type fini, de dimension $d$ $(d \in \mathbb N)$ et
tel que $M \!/ \mathfrak I M$ soit de longueur finie, la
multiplicité (de Samuel) de $M$ relativement à $\mathfrak I$,
notée $e_{\mathfrak I}(M)$, est définie par la formule de Samuel:
$$e_{\mathfrak I}(M) := d! \lim_{n \rightarrow \infty}\!\left(\frac{\longu_{A}\!\left(M \!/ {\mathfrak I}^n M\right)}
{n^d}\right) .$$ En particulier, lorsque $A$ est local, régulier,
d'idéal maximal $\mathfrak M$ et $\mathfrak I$ est $\mathfrak
M$-primaire (donc $\mathfrak I$ est une puissance de $\mathfrak
M$; $\mathfrak I = {\mathfrak M}^r$ pour un certain entier $r \geq
1$), $e_{\mathfrak I}(A)$ est simplement la multiplicité de $A$
relativement à $\mathfrak I$ (c'est-à-dire $e_{\mathfrak I}(A) =
r$). Pour les propriétés de la multiplicité de Samuel, on se
ref\`ere à \cite{Bou}, chapitre $8$, §$7.1$.

Soient maintenant $G \hookrightarrow {\mathbb P}_N$ le groupe
algébrique commutatif, irréductible, du §$4$ et $\mathfrak G$
l'idéal premier homog\`ene de définition de l'adhérence de Zariski
$\Gbar$ de $G$, dans $K[\underline X]$ ($\underline X := (X_0 ,
\dots , X_N)$). Appelons $A$ l'anneau:
$$A := K[\underline X] / \mathfrak{G}$$
et
$$s\: K[\underline X] \to A$$
la surjection canonique de $K[\underline X]$ sur $A$. \\
Pour tout idéal $\mathfrak Q$ de $A$, $s$ induit un isomorphisme
d'anneaux entre $K[\underline X] / s^{- 1}(\mathfrak Q)$ et $A /
\mathfrak Q$. De plus, ces anneaux sont isomorphes aussi en tant
qu'anneaux gradués et, pour un $\delta \in \mathbb N$ donné, les
deux $K$-espaces vectoriels de dimensions finies ${(K[\underline
X] / s^{- 1}(\mathfrak Q))}_{\delta}$ et ${(A / \mathfrak
Q)}_{\delta}$ (formés des éléments homog\`enes de degrés $\delta$
de $K[\underline X] / s^{- 1}(\mathfrak Q)$ et de $A / \mathfrak
Q$ respectivement) sont isomorphes et ont donc la m\^eme
dimension. On a ainsi par définition m\^eme du degré d'un idéal:
\begin{equation*}
\begin{split}
\deg_{A}(\mathfrak Q) &:= \lim_{\delta \to \infty}\!\left(\dim\left(A / \mathfrak Q\right)\right)\!! \frac{\dim_K{\!\left(A / \mathfrak Q\right)}_{\delta}}{\delta^{\dim\left(A / \mathfrak Q\right)}} \\
&~= \lim_{\delta \to \infty}\!\left(\dim\left(K[\underline X] / s^{- 1}(\mathfrak Q)\right)\right)\!! \frac{\dim_K{\!\left(K[\underline X] / s^{- 1}(\mathfrak Q)\right)}_{\delta}}{\delta^{\dim\left(K[\underline X] / s^{- 1}(\mathfrak Q)\right)}} \\
&~= \deg_{K[\underline X]}\!\left(s^{- 1}(\mathfrak Q)\right) ,
\end{split}
\end{equation*}
i.e:
$$\deg_{A}(\mathfrak Q) = \deg_{K[\underline X]}\!\left(s^{- 1}(\mathfrak Q)\right) .$$
En d'autres termes, cela veut dire que le degré est conservé lors
du passage d'un idéal $\mathfrak Q$ de $A$ à l'idéal $s^{-
1}(\mathfrak Q)$ de $K[\underline X]$. Pour avoir une propriété
analogue pour les hauteurs, remarquons que pour tout idéal premier
de $A$ l'idéal $s^{- 1}(\mathfrak Q)$ de $K[\underline X]$ est
aussi premier et posons par définition:
$$h(\mathfrak Q) := h(s^{- 1}(\mathfrak Q)) ,$$
o\`u $h(s^{- 1}(\mathfrak Q))$ désigne la hauteur unitaire d'une
forme de Chow de $s^{- 1}(\mathfrak Q)$ dans $K[\underline X]$. La
hauteur est ainsi conservée lors du passage d'un idéal premier
$\mathfrak Q$ de $A$ à l'idéal premier $s^{- 1}(\mathfrak Q)$ de
$K[\underline X]$. Par conséquent, elle est aussi conservée lors
du passage d'un idéal quelconque $\mathfrak Q$ de $A$ à l'idéal
correspondant $s^{- 1}(\mathfrak Q)$ de $K[\underline X]$.

Le lemme suivant rappelle le lemme $5$ de \cite{Ph3} et donne un
résultat analogue pour la hauteur.
\begin{lemma}\label{a.10}
Soit $\mathfrak I$ un idéal de $A$ engendré par des formes de
degrés $\leq \delta$ ($\delta \in \mathbb N$) et hauteurs $\leq
\eta$ ($\eta > 0$) et soit $d := \rang \mathfrak I$, alors on a:
\begin{equation*}
\!\!\!\!\sum_{\begin{array}{c}
\scriptstyle{\mathfrak P \in \Ass \mathfrak I} , ~\!\septpt{\mbox{minimal}} \\
\scriptstyle{\rang \mathfrak P = d}
\end{array}} \!\!\!\!\!\!\!\!\!\!\!\!\!\!\!\!e_{\mathfrak I A_{\mathfrak P}}\!\!\left(A_{\mathfrak P}\right) . \deg(\mathfrak P) \leq \deg(\mathfrak G) . \delta^d
\end{equation*}
et
\begin{equation*}
\!\!\!\!\sum_{\begin{array}{c}
\scriptstyle{\mathfrak P \in \Ass \mathfrak I} , ~\!\septpt{\mbox{minimal}} \\
\scriptstyle{\rang \mathfrak P = d}
\end{array}} \!\!\!\!\!\!\!\!\!\!\!\!\!\!\!\!\!\!e_{\mathfrak I A_{\mathfrak P}}\!\!\left(A_{\mathfrak P}\right) . h(\mathfrak P) \leq h(\mathfrak G) . \delta^d + g \!\left(\eta + 3 \log{\!\left(\deg(\mathfrak G) . \delta^g + 1\right)}\right) \deg(\mathfrak G) . \delta^{d - 1} \!.
\end{equation*}
\end{lemma}
{\bf Démonstration.---} Soit $\mathfrak I$ un idéal de $A$
engendré par des formes de degrés $\leq \delta$ $(\delta \in
{\mathbb N}^*)$ et de hauteurs $\leq \eta$ $(\eta > 0)$. En
suivant \cite{Ph3} (début de la démonstration du lemme $5$ de
cette référence), on construit des formes $P_1 , \dots , P_d$ de
$\mathfrak I$, de degrés $\leq \delta$, de hauteurs $\leq \eta + 2
d \log{\delta}$ et telles que $\dim(A_{\mathfrak P} / P_1
A_{\mathfrak P} + \dots + P_d A_{\mathfrak P}) = 0$. D'o\`u,
d'apr\`es \cite{Bou} (chapitre $8$, §$7.5$, théor\`eme 1.b):
\begin{equation*}
\begin{split}
e_{\mathfrak I A_{\mathfrak P}}(A_{\mathfrak P}) &\leq e_{\mathfrak I A_{\mathfrak P}}(A_{\mathfrak P} / P_1 A_{\mathfrak
P} + \dots + P_d A_{\mathfrak P}) \\
&\leq \longu_{A_{\mathfrak P}}\!\!\left(A_{\mathfrak P}/\left(P_1
A_{\mathfrak P} + \dots + P_d A_{\mathfrak P}\right)\right) .
\end{split}
\end{equation*}
En appelant $\mathfrak J$ l'idéal de $A$ engendré par les $d$
formes $P_1 , \dots , P_d$, on a donc:
$$e_{\mathfrak I A_{\mathfrak P}}(A_{\mathfrak P}) \leq \longu_{A_{\mathfrak P}}\!\!\left(A_{\mathfrak P}/\mathfrak J A_{\mathfrak P}\right)$$
pour tout $\mathfrak P \in \Ass \mathfrak I$ minimal de rang $d$.
Ce qui donne, en sommant sur tous ces idéaux $\mathfrak P$, les
inégalités:
\begin{align}
\sum_{\begin{array}{c}
\scriptstyle{\mathfrak P \in \Ass \mathfrak I} , ~\!\septpt{\mbox{minimal}} \\
\scriptstyle{\rang \mathfrak P = d}
\end{array}} \!\!\!\!\!\!\!\!\!\!\!\!\!\!\!\!e_{\mathfrak I A_{\mathfrak P}}\!\!\left(A_{\mathfrak P}\right) . \deg(\mathfrak P) &~\leq \sum_{\begin{array}{c}
\scriptstyle{\mathfrak P \in \Ass \mathfrak I} , ~\!\septpt{\mbox{minimal}} \\
\scriptstyle{\rang \mathfrak P = d}
\end{array}} \!\!\!\!\!\!\!\!\!\!\!\!\!\!\!\!\longu_{A_{\mathfrak P}}\!\!\left(A_{\mathfrak P}/\mathfrak J A_{\mathfrak P}\right) . \deg(\mathfrak P) \label{1.14} \\
\intertext{et} \sum_{\begin{array}{c}
\scriptstyle{\mathfrak P \in \Ass \mathfrak I} , ~\!\septpt{\mbox{minimal}} \\
\scriptstyle{\rang \mathfrak P = d}
\end{array}} \!\!\!\!\!\!\!\!\!\!\!\!\!\!\!\!e_{\mathfrak I A_{\mathfrak P}}\!\!\left(A_{\mathfrak P}\right) . h(\mathfrak P) &~\leq \sum_{\begin{array}{c}
\scriptstyle{\mathfrak P \in \Ass \mathfrak I} , ~\!\septpt{\mbox{minimal}} \\
\scriptstyle{\rang \mathfrak P = d}
\end{array}} \!\!\!\!\!\!\!\!\!\!\!\!\!\!\!\!\longu_{A_{\mathfrak P}}\!\!\left(A_{\mathfrak P}/\mathfrak
J A_{\mathfrak P}\right) . h(\mathfrak P) . \label{1.15}
\end{align}
Or, pour $\mathfrak P \in \Ass \mathfrak I$ minimal de rang $d$,
on montre ci-dessous que l'idéal $\mathfrak J A_{\mathfrak P}$ de
$A_{\mathfrak P}$ est primaire:
$$\mathfrak J A_{\mathfrak P} = {\mathfrak Q}^{(\mathfrak P)} A_{\mathfrak P}$$
o\`u ${\mathfrak Q}^{(\mathfrak P)}$ désigne l'unique composante
primaire de $\mathfrak J$ associée à $\mathfrak P$. En effet, en
considérant une décomposition primaire normale $\cap_{i \in I}
{\mathfrak Q}_i = \mathfrak J$ pour $\mathfrak J$, on obtient pour
l'idéal $\mathfrak J A_{\mathfrak P}$ de $A_{\mathfrak P}$ la
décomposition primaire normale:
$$\mathfrak J A_{\mathfrak P} = \bigcap_{\scriptstyle{i \in I , {\mathfrak Q}_i \subset \mathfrak P}} {\mathfrak Q}_i A_{\mathfrak P}$$
(car si ${\mathfrak Q}_i \not\subset \mathfrak P$ alors ${\mathfrak Q}_i A_{\mathfrak P} = A_{\mathfrak P}$ ne compte pas dans l'intersection). \\
Comme, par construction m\^eme des formes $P_1 , \dots , P_d$,
l'idéal $\mathfrak J A_{\mathfrak P}$ de $A_{\mathfrak P}$ est une
intersection compl\`ete, donc pure de rang $d$ (d'apr\`es le
théor\`eme de Cohen-Macauley), toute composante primaire
${\mathfrak Q}_i A_{\mathfrak P}$ $(i \in I , {\mathfrak Q}_i
\subset \mathfrak P)$ de $\mathfrak J A_{\mathfrak P}$ doit \^etre
d'idéal premier associé de rang $d$, c'est-à-dire que pour tout $i
\in I$ tel que ${\mathfrak Q}_i \subset \mathfrak P$:
$$\rang \sqrt{{\mathfrak Q}_i} A_{\mathfrak P} = d .$$
Il s'ensuit que pour tout $i \in I$ tel que ${\mathfrak Q}_i
\subset \mathfrak P$, les deux idéaux premiers $\sqrt{{\mathfrak
Q}_i} A_{\mathfrak P}$ et $\mathfrak P A_{\mathfrak P}$ de
$A_{\mathfrak P}$ sont de m\^eme rang (égal à $d$) et de plus le
premier est inclus dans le second, donc forcément
$\sqrt{{\mathfrak Q}_i} A_{\mathfrak P} = \mathfrak P A_{\mathfrak
P}$, par suite $\sqrt{{\mathfrak Q}_i} = \mathfrak P$. On vient de
montrer que toute composante primaire de $\mathfrak J$ contenue
dans $\mathfrak P$ est d'idéal premier associé $\mathfrak P$.
Comme la décomposition primaire choisie est normale, il ne peut
exister qu'une unique composante primaire de $\mathfrak J$
contenue dans $\mathfrak P$ et, en la notant ${\mathfrak
Q}^{(\mathfrak P)}$ on a: $\mathfrak J A_{\mathfrak P} =
{\mathfrak Q}^{(\mathfrak P)} A_{\mathfrak P}$. Ce qui démontre
notre assertion.

On a ainsi pour tout $\mathfrak P \in {\rm{Ass}} \mathfrak I$,
minimal de rang $d$:
\begin{align}
\longu_{A_{\mathfrak P}}\!\!\left(A_{\mathfrak P}/\mathfrak J A_{\mathfrak P}\right) \deg(\mathfrak P) &=
\longu_{A_{\mathfrak P}}\!\!\left(A_{\mathfrak P}/{\mathfrak Q}^{(\mathfrak P)} A_{\mathfrak P}\right)
\deg(\mathfrak P) = \deg\!\left({\mathfrak Q}^{(\mathfrak P)}\right) \notag \\
\intertext{et} \longu_{A_{\mathfrak P}}\!\!\left(A_{\mathfrak
P}/\mathfrak J A_{\mathfrak P}\right) h(\mathfrak P) &=
\longu_{A_{\mathfrak P}}\!\!\left(A_{\mathfrak P}/{\mathfrak
Q}^{(\mathfrak P)} A_{\mathfrak P}\right) h(\mathfrak P) ~\!=
h\!\left({\mathfrak Q}^{(\mathfrak P)}\right) . \notag
\end{align}
Puis en reportant ces derni\`eres égalités dans (\ref{1.14}) et
(\ref{1.15}):
\begin{align}
\sum_{\begin{array}{c}
\scriptstyle{\mathfrak P \in \Ass \mathfrak I} , ~\!\septpt{\mbox{minimal}} \\
\scriptstyle{\rang \mathfrak P = d}
\end{array}} \!\!\!\!\!\!\!\!\!\!\!\!\!\!\!\!e_{\mathfrak I A_{\mathfrak P}}\!\!\left(A_{\mathfrak P}\right) . \deg(\mathfrak P) &~\leq \sum_{\begin{array}{c}
\scriptstyle{\mathfrak P \in \Ass \mathfrak I} , ~\!\septpt{\mbox{minimal}} \\
\scriptstyle{\rang \mathfrak P = d}
\end{array}} \!\!\!\!\!\!\!\!\!\!\!\!\!\!\!\!\deg\!\left({\mathfrak Q}^{(\mathfrak P)}\right) \label{1.16} \\
\intertext{et} \sum_{\begin{array}{c}
\scriptstyle{\mathfrak P \in \Ass \mathfrak I} , ~\!\septpt{\mbox{minimal}} \\
\scriptstyle{\rang \mathfrak P = d}
\end{array}} \!\!\!\!\!\!\!\!\!\!\!\!\!\!\!\!e_{\mathfrak I A_{\mathfrak P}}\!\!\left(A_{\mathfrak P}\right) . h(\mathfrak P) &~\leq \sum_{\begin{array}{c}
\scriptstyle{\mathfrak P \in \Ass \mathfrak I} , ~\!\septpt{\mbox{minimal}} \\
\scriptstyle{\rang \mathfrak P = d}
\end{array}} \!\!\!\!\!\!\!\!\!\!\!\!\!\!\!\!h\!\left({\mathfrak Q}^{(\mathfrak P)}\right) \label{1.17} .
\end{align}
Soit maintenant $\mathcal O$ un ouvert de Zariski de $K[\underline
X]$ contenant tous les $s^{- 1}(\mathfrak P)$ pour $\mathfrak P
\in {\rm{Ass}} \mathfrak I$, minimal de rang $d$. Comme $\mathfrak
G$ est $\mathcal O$-parfait\footnote{La définition d'un idéal
$\mathcal O$-parfait est rappelée dans \cite{Ph4}3 page $358$.}
(car pour tout $\mathfrak P \in {\rm{Ass}} \mathfrak I$, minimal
de rang $d$, $A_{\mathfrak P}$ est semi-régulier), on déduit du
théor\`eme $6$ de \cite{Ph4}3 (qu'on applique avec $R_m$ remplacé
par $K[\underline X]$, l'idéal homog\`ene $I$ de $R_m$ remplacé
par $\mathfrak G$ et l'idéal $J = (I , p_1 , \dots , p_{\ell})$ de
cette référence remplacé par $s^{- 1}(\mathfrak J) = (\mathfrak G
, P_1 , \dots , P_d)$):
\begin{equation*}
\begin{split}
\left(\sum_{\begin{array}{c}
\scriptstyle{\mathfrak P \in \Ass \mathfrak I} , ~\!\septpt{\mbox{minimal}} \\
\scriptstyle{\rang \mathfrak P = d}
\end{array}} \!\!\!\!\!\!\!\!\!\!\!\!\!\!\!\!\deg\!\left({\mathfrak Q}^{(\mathfrak P)}\right)\right) \delta^{g - d} &\leq \sum_{s^{- 1}(\mathfrak P) \in \mathcal O} \deg\!\left(s^{- 1}\!\left({\mathfrak Q}^{(\mathfrak P)}\right)\right) \delta^{g - d} \\
&\leq d_{\mathcal O}\left(Z\left(s^{- 1}(\mathfrak J)\right)\!; \delta\right) \\
&\leq d_{\mathcal O}\left(Z(\mathfrak G) ; \delta\right) = \deg(\mathfrak G) . \delta^g \\
\intertext{et} \left(\sum_{\begin{array}{c}
\scriptstyle{\mathfrak P \in \Ass \mathfrak I} , ~\!\septpt{\mbox{minimal}} \\
\scriptstyle{\rang \mathfrak P = d}
\end{array}} \!\!\!\!\!\!\!\!\!\!\!\!\!\!\!\!h\!\left({\mathfrak Q}^{(\mathfrak P)}\right)\right) \delta^{g - d + 1} &\leq \sum_{s^{- 1}(\mathfrak P) \in \mathcal O} h\!\left(s^{- 1}\!\left({\mathfrak Q}^{(\mathfrak P)}\right)\right) \delta^{g - d + 1} \\
&\!\!\!\!\!\!\!\!\!\!\!\!\!\!\!\!\!\!\!\!\!\!\!\!\!\!\!\!\!\!\!\!\!\!\!\!\!\!\!\!\!\!\!\!\!\!\!\!\!\!\!\!\!\!\!\!\!\!\!\leq h_{\mathcal O}\left(Z\left(s^{- 1}(\mathfrak J)\right) \!; \delta\right) \\
&\!\!\!\!\!\!\!\!\!\!\!\!\!\!\!\!\!\!\!\!\!\!\!\!\!\!\!\!\!\!\!\!\!\!\!\!\!\!\!\!\!\!\!\!\!\!\!\!\!\!\!\!\!\!\!\!\!\!\!\leq h_{\mathcal O}\left(Z(\mathfrak G) ; \delta\right) + g\left(\eta + 2 d \log{\delta} + \log{\left(d_{\mathcal O}\left(Z(\mathfrak G) ; \delta\right) + 1\right)}\right) d_{\mathcal O}\left(Z(\mathfrak G) ; \delta\right) \\
&\!\!\!\!\!\!\!\!\!\!\!\!\!\!\!\!\!\!\!\!\!\!\!\!\!\!\!\!\!\!\!\!\!\!\!\!\!\!\!\!\!\!\!\!\!\!\!\!\!\!\!\!\!\!\!\!\!\!\!\leq h(\mathfrak G) \delta^{g + 1} + g\left(\eta + 2 d \log{\delta} + \log{\left(\deg(\mathfrak G) \delta^g + 1\right)}\right) \deg(\mathfrak G) \delta^g \\
&\!\!\!\!\!\!\!\!\!\!\!\!\!\!\!\!\!\!\!\!\!\!\!\!\!\!\!\!\!\!\!\!\!\!\!\!\!\!\!\!\!\!\!\!\!\!\!\!\!\!\!\!\!\!\!\!\!\!\!\leq
h(\mathfrak G) \delta^{g + 1} + g\left(\eta + 3
\log{\left(\deg(\mathfrak G) \delta^g + 1\right)}\right)
\deg(\mathfrak G) \delta^g .
\end{split}
\end{equation*}
D'o\`u:
$$
\sum_{\begin{array}{c}
\scriptstyle{\mathfrak P \in \Ass \mathfrak I} , ~\!\septpt{\mbox{minimal}} \\
\scriptstyle{\rang \mathfrak P = d}
\end{array}} \!\!\!\!\!\!\!\!\!\!\!\!\!\!\!\!\deg\!\left({\mathfrak Q}^{(\mathfrak P)}\right) \leq \deg(\mathfrak G) \delta^d
$$
et
$$
\sum_{\begin{array}{c}
\scriptstyle{\mathfrak P \in \Ass \mathfrak I} , ~\!\septpt{\mbox{minimal}} \\
\scriptstyle{\rang \mathfrak P = d}
\end{array}} \!\!\!\!\!\!\!\!\!\!\!\!\!\!\!\!h\!\left({\mathfrak Q}^{(\mathfrak P)}\right) \leq h(\mathfrak G) \delta^d + g\left(\eta + 3 \log{\left(\deg(\mathfrak G) \delta^g + 1\right)}\right) \deg(\mathfrak G) \delta^{d - 1} .
$$
Il ne reste qu'à reporter ces deux derni\`eres inégalités dans
(\ref{1.16}) et (\ref{1.17}) pour aboutir au lemme \ref{a.10}. La
démonstration est achevée. $~~~~\blacksquare$\vspace{1mm}
\subsection{\'Enoncé du lemme de zéros:}
\`A partir de ce paragraphe, on utilisera la notion de la \og
multiplicité d'une sous-variété $V$ de $G$ relativement \`a un
dessous d'escalier $W$ de ${\mathbb N}^g$\fg$~\!\!$ définie dans
\cite{Ph3}, page $1072$. Cette multiplicité est notée $m_W(V)$. On
démontre le théorème suivant:
\begin{theorem}\label{a.13}
Soit $G$ un groupe algébrique commutatif de dimension $g$ $(g \geq
1)$, défini sur un corps de nombres $K$, plongé dans un espace
projectif ${\mathbb P}_N$$(N \geq 1)$ et d'élément neutre $\mathbf
e$ représenté dans ${\mathbb P}_N$ par un syst\`eme de coordonnées
projectives $\underline e = (e_0 : \dots : e_N)$. On désigne par
$K[\underline X]$$(\underline X = (X_0 , \dots , X_N))$ l'anneau
de coordonnées de ${\mathbb P}_N$, par $\mathfrak G$ l'idéal de
définition de l'adhérence de Zariski $\Gbar$ de $G$ dans
$K[\underline X]$ et par $A := K[\underline X] / \mathfrak{G}$
l'anneau de coordonnées de $\Gbar$. Soient aussi $\mathbf x$ un
point fixé de $G$, $\underline A$ une famille de formes de
$K[\underline X , \underline Y]$$(\underline X := (X_0 , \dots ,
X_N) , \underline Y := (Y_0 , \dots , Y_N))$ représentant
l'addition dans $G \plongement {\mathbb P}_N$ au voisinage de
$\{\mathbf e\} \times \{\mathbf x\}$ et $(c , c') \in {{\mathbb
N}^*}^2$ le bidegré des formes constituant $\underline A$. Soient
enfin $W_0 , \dots , W_g$ des dessous d'escaliers finis de
${\mathbb N}^g$ et $P$ une forme de $A$ non identiquement nulle de
degré $\delta$ $(\delta \in {\mathbb N}^*)$. Sous les hypoth\`eses
du §$4.1$ pour $G$, si $P$ s'annule en $\mathbf x$ avec une
multiplicité définie par le dessous d'escalier $W_0 + \dots +
W_g$, alors il existe $1 \leq r \leq g$ et une sous-variété $V$ de
$G$ de dimension $d \leq g - r$, contenant $\mathbf x$,
incompl\`etement définie\footnote{selon la terminologie de
\cite{Ph2}.} dans $G$ par des équations de degrés $c' \delta$ et
satisfaisant:
\begin{equation*}
\begin{split}
m_{W_r}(V) \deg{V} &\leq \deg{G} . (c' \delta)^{g - d} \\
\intertext{et}
m_{W_r}(V) h(V) &\leq h(G) (c' \delta)^{g - d} + g \!\!\left[\sum_{\ell = 1}^{p} \left(\log{(E_{\ell})} + \log{2}\right) t_{\ell}\!\left(\!W_0 + \dots + W_{r - 1}\!\right)\right. \\
&~~~~~~~~~~+ c \delta \!\left(\sum_{\ell = 1}^{p} \delta_{\ell} \!\left(h({\underline e}_{\ell}) + g_{\ell} \log{s}\!\right) + g \log{2}\!\right) + \widetilde{h}(P) + \delta \widetilde{h}(\underline A) \\
&\!\!\!\!\!\!\!\!\!\!\!\!\!\!\!\!\!\!\!\!\!\!\left.+ \delta (2 c +
2 c' + 1) (\log{N} + 1) + 3 \log\!\left((\deg{G}) . (c' \delta)^g
+ 1\right)\!\phantom{\sum_{\ell = 1}^{p} \left(\log{E_{\ell}} +
\log{(c \delta)}\right) t_{\ell}\!\left(\!W_0 + \dots + W_{r -
1}\!\right)}\!\!\!\!\!\!\!\!\!\!\!\!\!\!\!\!\!\!\!\!\!\!\!\!\!\!\!\!\!\!\!\!\!\!\!\!\!\!\!\!\!\!\!\!\!\!\!\!\!\!\!\!\!\!\!\!\!\!\!\!\!\!\!\!\!\!\!\!\!\!\!\!\!\!\!\!\!\!\!\!\!\!\!\!\!\!\!\!\!\!\!\!\!\!\!\!\!\!\!\!\!\!\!\!\!\!\right]\!(\deg{G})(c'
\delta)^{g - d - 1} .
\end{split}
\end{equation*}
De plus, $P$ s'annule sur $V$ avec la multiplicité définie par le
dessous d'escalier $W_0 + \dots + W_r$.
\end{theorem}
{\bf Démonstration.---} On consid\`ere la suite croissante
d'idéaux homog\`enes de $A$:
$$\Delta^{W_0}(P) = {\mathfrak I}_1 \subset {\mathfrak I}_2 \subset \dots \subset {\mathfrak I}_{g + 1} \subset {\mathfrak M}_{\mathbf x}$$
définis par:
$${\mathfrak I}_i := \Delta^{W_0 + \dots + W_{i - 1}}(P) ~~~~~~\text{(pour $1 \leq i \leq g + 1$)}$$
et ${\mathfrak M}_{\mathbf x}$ l'idéal de définition dans $A$ de la sous-variété de $G$ réduite au point $\mathbf x$. \\
Notons que les inclusions d'idéaux ${\mathfrak I}_1 \subset {\mathfrak I}_2 \subset \dots \subset {\mathfrak I}_{g + 1}$ résultent des inclusions de dessous d'escaliers $W_0 \subset W_0 + W_1 \subset \dots \subset W_0 + \dots + W_g$ qui ont lieu puisqu'un dessous d'escalier non vide de ${\mathbb N}^g$ contient l'origine $\underline 0$ de ${\mathbb N}^g$. Par ailleurs l'inclusion ${\mathfrak I}_{g + 1} \subset {\mathfrak M}_{\mathbf x}$ exprime simplement l'hypoth\`ese: ``$P$ s'annule en $\mathbf x$ avec une multiplicité définie par le dessous d'escalier $W_0 + \dots + W_g$''. \\
Pour tout $i = 1 , \dots , g + 1$, appelons $E_i$ le sous-ensemble
algébrique de $G$ défini par l'idéal ${\mathfrak I}_i$ de $A$. On
a alors:
$$G \varsupsetneq E_1 \supset E_2 \supset \dots \supset E_{g + 1} \supset \{\mathbf x\} .$$
On peut choisir, pour tout $i = 1 , \dots , g + 1$, une composante
irréductible isolée $V_i$ de $E_i$ de façon à avoir:
$$G \varsupsetneq V_1 \supset V_2 \supset \dots \supset V_{g + 1} \supset \{\mathbf x\} .$$
En effet, on proc\`ede par récurrence: on commence par choisir $V_{g + 1}$ une composante irréductible isolée quelconque de $E_{g + 1}$ contenant le point $\mathbf x$ (ce qui est possible puisque $\mathbf x \in E_{g + 1}$), puis en supposant que $V_{i + 1}$ est construite (pour un certain $1 \leq i \leq g$), il suffit de prendre pour $V_i$ une composante irréductible isolée de $E_i$ contenant $V_{i + 1}$ (ce qui est possible puisque $V_{i + 1} \subset E_{i + 1} \subset E_i$) et on obtient ainsi des sous-variétés $V_1 , \dots , V_{g + 1}$ de $G$, emboitées comme ci-dessus. \\
Posons maintenant pour tout $i = 1 , \dots , g + 1$:
$$d_i := \dim{V_i} ,$$
on a alors:
$$g > d_1 \geq d_2 \geq \dots \geq d_{g + 1} \geq 0 .$$
Cette derni\`ere chaine d'inégalités entra\^{\i}ne l'existence
d'au moins un entier $1 \leq r \leq g$ pour lequel on a $d_r =
d_{r + 1}$. De plus, en prenant $r$ le plus petit possible ayant
cette propriété, c'est-à-dire:
$$g > d_1 > d_2 > \dots > d_r = d_{r + 1} ,$$
on a
$$d_r = d_{r + 1} \leq g - r .$$
Les deux sous-variétés algébriques $V_r$ et $V_{r + 1}$ de $G$ ont
donc la m\^eme dimension $d := d_r = d_{r + 1} \leq g - r$ et
comme $V_r \supset V_{r + 1}$, alors forcément $V_r = V_{r + 1}$.
En posant:
$$V := V_r = V_{r + 1} ,$$
$V$ est bien une sous-variété algébrique de $G$, de dimension $d \leq g - r$, contenant $\mathbf x$ et incompl\`etement définie dans $G$ par les idéaux ${\mathfrak I}_r$ et ${\mathfrak I}_{r + 1}$ qui sont engendrés par des formes de degrés $c' \delta$ (puisque $\deg{\Delta^{I}P} = c' \delta$ pour tout $I \in {\mathbb N}^g$). Montrons que $P$ s'annule sur $V$ avec la multiplicité définie par le dessous d'escalier $W_0 + \dots + W_r$: \\
D'apr\`es la définition m\^eme de l'idéal ${\mathfrak I}_{r + 1}$, pour tout $I \in W_0 + \dots + W_r$, la forme $\Delta^I P := \Delta_{\mathbf x}^{I}P$ s'annule sur $E_{r + 1}$ donc s'annule -a fortiori- sur $V_{r + 1} = V$; d'o\`u (en se reportant aux définitions \ref{a.14}) $P$ s'annule en tout point de $V \cap \Omega_{\mathbf x}$ avec la multiplicité définie par le dessous d'escalier $W_0 + \dots + W_r$, par suite $P$ s'annule sur presque toute $V$ avec la multiplicité définie par le dessous d'escalier $W_0 + \dots + W_r$ (car $V \setminus V \cap \Omega_{\mathbf x}$ est contenu dans un fermé de $V$). Si on a $\dim{V} \geq 1$, le lemme \ref{a.11} nous permet de conclure que $P$ s'annule sur $V$ toute enti\`ere avec la multiplicité définie par le dessous d'escalier $W_0 + \dots + W_r$; sinon $\dim{V} = 0$, c'est-à-dire $V = \{\mathbf x\}$ et dans ce cas -par hypoth\`ese m\^eme- $P$ s'annule sur $V$ avec la multiplicité définie par le dessous d'escalier $W_0 + \dots + W_g \supset W_0 + \dots + W_r$, donc s'annule -en tout cas- sur $V$ avec la multiplicité définie par le dessous d'escalier $W_0 + \dots + W_r$. Montrons finalement que $V$ satisfait bien les deux inégalités du théor\`eme \ref{a.13}: \\
Soit pour cela $\mathfrak P \in \Ass {\mathfrak I}_r \cap \Ass
{\mathfrak I}_{r + 1}$ le premier définissant $V$ dans $A$.
D'apr\`es le lemme \ref{a.12}, on a:
$$\Delta^{W_r}{\mathfrak I}_r \subset {\mathfrak I}_{r + 1} \subset \mathfrak P ,$$
par conséquent, en appliquant le lemme $7$ de \cite{Ph3}:
\begin{equation}
e_{{\mathfrak I}_r A_{\mathfrak P}}(A_{\mathfrak P}) \geq
m_{W_r}(V) . \label{1.18}
\end{equation}
L'idéal ${\mathfrak I}_r$ de $A$ est engendré par les formes
$\Delta^I P$ $(I \in W_0 + \dots + W_{r - 1})$ qui sont toutes de
degré $c' \delta$ et de hauteurs de Gauss-Weil majorées (gr\^ace à
la proposition \ref{a.16}) par:
\begin{equation*}
\begin{split}
\sum_{\ell = 1}^{p} \!\left(\log{(E_{\ell})} + \log{2}\right) t_{\ell}\!\!\left(\!W_0 + \dots + W_{r - 1}\!\right) \!+ \!c \delta \!\!\left(\!\sum_{\ell = 1}^{p} \!\delta_{\ell} \!\left(h({\underline e}_{\ell}) + g_{\ell} \log{s}\!\right) + g \log{2}\!\!\right) \\
+ \widetilde{h}(P) + \delta \widetilde{h}(\underline A) + \delta
(2 c + c' + 1) \log{(N + 1)} ,
\end{split}
\end{equation*}
puis de hauteurs unitaires majorées (gr\^ace à la relation de
comparaison (\ref{1.21}) entre hauteurs de Gauss-Weil et hauteurs
unitaires) par:
\begin{equation*}
\begin{split}
h\left(\Delta^I P\right) &\leq \widetilde{h}\left(\Delta^I P\right) + \frac{1}{2} \log{\binom{c' \delta + N}{N}} + c' \delta \sum_{j = 1}^{N} \frac{1}{2 j} \\
&\leq \widetilde{h}\left(\Delta^I P\right) + \frac{1}{2} c' \delta \log{(N + 1)} + \frac{1}{2} c' \delta (\log{N} + 1) \\
&\leq \sum_{\ell = 1}^{p} \left(\log{(E_{\ell})} + \log{2}\right) t_{\ell}\!\left(\!W_0 + \dots + W_{r - 1}\!\right) \\
&~~~~~~~~~~~~~~~~~~\!\!\!\!\!+ c \delta \!\!\left(\!\sum_{\ell = 1}^{p} \!\delta_{\ell} \!\left(h({\underline e}_{\ell}) + g_{\ell} \log{s}\!\right) + g \log{2}\!\!\right) + \widetilde{h}(P) + \delta \widetilde{h}(\underline A) \\
&~~~~~~~~~~~~~~~~~~~~~~~~~~~~~~~~~~~~~~~~~~~\!\!+ \delta (2 c + 2
c' + 1) (\log{N} + 1)
\end{split}
\end{equation*}
(pour tout $I \in W_0 + \dots + W_{r - 1}$). \\
Les inégalités du théor\`eme \ref{a.13} découlent de (\ref{1.18})
et de l'application du lemme \ref{a.10} à l'idéal ${\mathfrak
I}_r$. En effet, on a:
\begin{equation*}
\begin{split}
m_{W_r}(V) \deg{V} &\leq e_{{\mathfrak I}_r A_{\mathfrak P}}(A_{\mathfrak P}) \deg{\mathfrak P} ~~~~~~~~~~~~~~~~~~~~\text{(d'apr\`es (\ref{1.18}))} \\
&\leq \sum_{\begin{array}{c}
\scriptstyle{\mathfrak{P'} \in \Ass {\mathfrak I}_r , ~\!\mbox{minimal}} \\
\scriptstyle{\rang \mathfrak{P'} = g - d}
\end{array}} \!\!\!\!\!\!\!\!\!\!\!\!\!\!\!\!e_{{\mathfrak I}_r A_{\mathfrak{P'}}}\!\!\left(A_{\mathfrak{P'}}\right) . \deg \mathfrak{P'} \\
&\leq \deg{G}.(c' \delta)^{g - d} ~~~~~~~~~~~~~~~~~~~~\text{(d'apr\`es le lemme \ref{a.10})} \\
\intertext{et de m\^eme}
m_{W_r}(V) h(V) &\leq e_{{\mathfrak I}_r A_{\mathfrak P}}(A_{\mathfrak P}) h(\mathfrak P) \\
&\leq \sum_{\begin{array}{c}
\scriptstyle{\mathfrak{P'} \in \Ass {\mathfrak I}_r , ~\!\mbox{minimal}} \\
\scriptstyle{\rang \mathfrak{P'} = g - d}
\end{array}} \!\!\!\!\!\!\!\!\!\!\!\!\!\!\!\!e_{{\mathfrak I}_r A_{\mathfrak{P'}}}\!\!\left(A_{\mathfrak{P'}}\right) . h(\mathfrak{P'}) \\
&\leq h(G) (c' \delta)^{g - d} + g \!\!\left[\sum_{\ell = 1}^{p} \left(\log{(E_{\ell})} + \log{2}\right) t_{\ell}\!\left(\!W_0 + \dots + W_{r - 1}\!\right)\right. \\
&~~~~~~~~~~+ c \delta \!\left(\sum_{\ell = 1}^{p} \delta_{\ell} \!\left(h({\underline e}_{\ell}) + g_{\ell} \log{s}\!\right) + g \log{2}\!\right) + \widetilde{h}(P) + \delta \widetilde{h}(\underline A) \\
&\!\!\!\!\!\!\!\!\!\!\!\!\!\!\!\!\!\!\!\!\!\!\left.+ \delta (2 c +
2 c' + 1) (\log{N} + 1) + 3 \log\!\left((\deg{G}) . (c' \delta)^g
+ 1\right)\!\phantom{\sum_{\ell = 1}^{p} \left(\log{E_{\ell}} +
\log{(c \delta)}\right) t_{\ell}\!\left(\!W_0 + \dots + W_{r -
1}\!\right)}\!\!\!\!\!\!\!\!\!\!\!\!\!\!\!\!\!\!\!\!\!\!\!\!\!\!\!\!\!\!\!\!\!\!\!\!\!\!\!\!\!\!\!\!\!\!\!\!\!\!\!\!\!\!\!\!\!\!\!\!\!\!\!\!\!\!\!\!\!\!\!\!\!\!\!\!\!\!\!\!\!\!\!\!\!\!\!\!\!\!\!\!\!\!\!\!\!\!\!\!\!\!\!\!\!\!\right]\!(\deg{G})(c'
\delta)^{g - d - 1} .
\end{split}
\end{equation*}
Ce qui ach\`eve cette démonstration.
$~~~~\blacksquare$\vspace{1mm}
\section{{\bf LEMMES DE Z\'EROS EXPLICITES}}
\subsection{Un lemme de zéros explicite dans un groupe algébrique commutatif $G$ plongé dans un espace projectif ${\mathbb P}_N$:}
Soit $G$ un groupe algébrique commutatif défini sur un corps de
nombres $K$, irréductible, plongé dans ${\mathbb P}_N$$(N \geq
1)$, de dimension $g$$(g \geq 1)$ et d'élément neutre $\mathbf e$
représenté dans ${\mathbb P}_N$ par un syst\`eme de coordonnées
projectives $\underline e = (e_0 : e_1 : \dots : e_N)$. Soient
aussi $\mathbf x$ un point fixé de $G$, $\underline A$ une famille
de formes de $K[\underline X , \underline Y]$$(\underline X :=
(X_0 , \dots , X_N) , \underline Y := (Y_0 , \dots , Y_N))$
représentant l'addition dans $G \plongement {\mathbb P}_N$ au
voisinage de $\{\mathbf e\} \times \{\mathbf x\}$ et $(c , c') \in
{{\mathbb N}^*}^2$ le bidegré des formes constituant $\underline
A$. Nous désignons par $K[\underline X]$$(\underline X := (X_0 ,
\dots , X_N))$ l'anneau de coordonnées de ${\mathbb P}_N$, par
$\mathfrak G$ l'idéal de définition de l'adhérence de Zariski
$\Gbar$ de $G$ dans $K[\underline X]$ et par $A := K[\underline X]
/ \mathfrak G$ l'anneau de coordonnées de $\Gbar$.

Sous les hypoth\`eses de normalisation du §$3$ pour $G$ et
$\mathbf e \in G$, on a une paramétrisation de $G$ au voisinage de
l'origine $\mathbf e$ utilisant comme param\`etres $\underline T
:= (T_1 , \dots , T_g) := (X_1 - e_1 , \dots , X_g - e_g)$ et
représentée par un monomorphisme $\varphi$ de $A$ dans
$K[[\underline T]]$. De plus, en écrivant pour tout $i \in \{0 ,
\dots , N\}$:
$$\varphi(X_i) = \sum_{I \in {\mathbb N}^g} a_{I}^{(i)} {\underline T}^I ,$$
les nombres $a_{I}^{(i)}$$(I \in {\mathbb N}^g , i = 0 , \dots ,
N)$ de $K$ vérifient (d'apr\`es le lemme \ref{a.6}), pour toute
place $v$ de $K$, la relation:
\begin{equation}
{\nb{a_{I}^{(i)}}}_v \leq {E_v(N , G , \underline e)}^{\nb{I}} .
H_v(\underline e) ~~~~ (I \in {\mathbb N}^g , i = 0 , \dots , N) ,
\label{1.6}
\end{equation}
o\`u les $E_v$$(v \in M_K)$ sont les expressions (supérieures ou
égales à $1$) dépendant de $N$, $G$ et $\underline e$ définies
par:
\begin{align}
E_v &:= {H_v(\underline e)}^{2 (d(G) - 1)} \prod_{j = g + 1}^{N} \!\!\!\left({H_v({\widetilde P}_j)}^2 {\max\!\left\{1 , {\nb{\frac{1}{\frac{\partial P_j}{\partial X_j}(e_1 , \dots , e_g , e_j)}}}_v\right\}}^{\!\!2}\right) \notag \\
\intertext{si $v$ est finie et}
E_v &:= 8 g {d(G)}^3 (d(G) + 1)^{2 (g + 1)} . {H_v(\underline e)}^{2 (d(G) - 1)} \notag \\
&~~~~~~~~~~~~~~~~~~~~~~~~~~\!\!\times \prod_{j = g + 1}^{N}
\!\!\!\left({H_v({\widetilde P}_j)}^2 {\max\!\left\{1 ,
{\nb{\frac{1}{\frac{\partial P_j}{\partial X_j}(e_1 , \dots , e_g
, e_j)}}}_v\right\}}^{\!\!2}\right) \notag
\end{align}
si $v$ est infinie. \\
Les relations (\ref{1.6}) montrent bien que les hypoth\`eses du
§$4.1$ pour $G$ sont satisfaites avec $p = 1$, $n_1 = N$,
$\delta_1 = 1$, ${\underline e}_1 = \underline e$, $t_1$ la
fonction additive de ${\mathbb N}^g$ dans ${\mathbb R}^+$ définie
par:
$$t_1(I) = \nb{I} ~~~~~~ \text{(pour tout $I \in {\mathbb N}^g$)}$$
et $E_{1 , v} = E_v$ pour toute place $v$ de $K$. Ceci nous permet
donc d'appliquer le théor\`eme \ref{a.13} afin d'obtenir -sous les
hypoth\`eses de normalisation faites sur $G$- un lemme de zéros
explicite dans $G$. Avant cette application, nous posons:
$$E := \prod_{v \in M_K} E_{v}^{\frac{[K_v : {\mathbb Q}_v]}{[K : \mathbb Q]}}$$
($E$ correspond à $E_1$ du §$4.1$) et nous majorons dans le lemme
qui suit $\log{E}$ (ou $\log{E} + \log{4}$ m\^eme) par une
expression qui dépend seulement de $N$, $g$, $d(G)$, $h(G)$ et
$h(\mathbf e)$ (dans le but de se débarrasser des polyn\^omes
$P_j$$(g + 1 \leq j \leq N)$ figurant dans les $E_v$$(v \in
M_K)$).
\begin{lemma}\label{a.9}
On a:
\begin{equation*}
\begin{split}
\log{(E)} + \log{4} &\leq 4 (N - g) h(G) + \left[2 (N - g + 1) h(\mathbf e) + 4 (N - g)\right] d(G) \\
&~~~~~~~~\!\!\!\!\!\!\!+ (N - g + 1) (2 g + 5) (\log{(d(G))} + 1)
- 2 (N - g + 1) h(\mathbf e) .
\end{split}
\end{equation*}
\end{lemma}
{\bf Démonstration.---} Par définition m\^eme de $E$ et des
$E_v$$(v \in M_K)$ on a:
\begin{equation*}
\begin{split}
\log{(E)} &= \sum_{v \in M_K} \frac{[K_v : {\mathbb Q}_v]}{[K : \mathbb Q]} \log{(E_v)} \\
&= 2 \sum_{j = g + 1}^{N} \widetilde{h}({\widetilde P}_j) + 2 \sum_{j = g + 1}^{N} h\!\left(1 : \frac{1}{\frac{\partial P_j}{\partial X_j}(e_1 , \dots , e_g , e_j)}\right) + 2 (d(G) - 1) h(\mathbf e) \\
&~~~~~~~~~~~~~~~~~~~~~~\!\!+ \log{(8 g)} + 3 \log{(d(G))} + 2 (g +
1) \log{(d(G) + 1)} .
\end{split}
\end{equation*}
Or, on peut montrer en raisonnant place par place qu'on a pour
tout $j \in \{g + 1 , \dots , N\}$:
\begin{equation*}
\begin{split}
h\!\left(1 : \frac{1}{\frac{\partial P_j}{\partial X_j}(e_1 , \dots , e_g , e_j)}\right) &= h\!\left(1 : \frac{\partial P_j}{\partial X_j}(e_1 , \dots , e_g , e_j)\right) \\
&\!\!\!\!\!\!\!\!\!\!\!\!\!\!\!\!\!\!\!\!\!\!\!\!\!\!\!\!\!\!\!\!\!\!\!\!\!\!\!\!\!\!\!\!\!\!\!\!\!\!\!\!\!\!\!\leq \widetilde{h}(P_j) + \left({d°}_{\rm{tot}}P_j - 1\right) h(1 : e_1 : \dots : e_g : e_j) + \log{\!\left(\!{d°}_{\rm{tot}}P_j \binom{{d°}_{\rm{tot}}P_j + g}{g + 1}\!\right)} \\
&\!\!\!\!\!\!\!\!\!\!\!\!\!\!\!\!\!\!\!\!\!\!\!\!\!\!\!\!\!\!\!\!\!\!\!\!\!\!\!\!\!\!\!\!\!\!\!\!\!\!\!\!\!\!\!\leq
\widetilde{h}({\widetilde P}_j) + (d(G) - 1) h(\mathbf e) + (g +
2) \log{(d(G))}
\end{split}
\end{equation*}
car: $\widetilde{h}(P_j) = \widetilde{h}({\widetilde P}_j)$,
${d°}_{\rm{tot}}P_j \leq {d°}_{\rm{tot}}{\widetilde P}_j = d(\pi_j
(G)) \leq d(G)$ et le coefficient binomial
$\binom{{d°}_{\rm{tot}}P_j + g}{g + 1}$ étant majoré par
$({d°}_{\rm{tot}}P_j)^{g + 1} \leq {d(G)}^{g + 1}$. D'o\`u:
\begin{equation*}
\begin{split}
\log{(E)} + \log{4} \leq 4 \sum_{j = g + 1}^{N} \!\!\widetilde{h}({\widetilde P}_j) + 2 (N - g + 1) (d(G) - 1)
h(\mathbf e) ~~~~~~~~~~~~~~~~~~~~~~~~~~~~~~~\\
\!\!\!\!\!\!\!\!\!\!\!\!\!\!\!+ (N - g) (2 g + 4) \log{(d(G))} + \log{(8 g)} + 3 \log{(d(G))} + (2 g + 2) \log{(d(G) + 1)}
+ \log{4}\\
\leq 4 \sum_{j = g + 1}^{N} \!\!\widetilde{h}({\widetilde P}_j) +
2 (N - g + 1) (d(G) - 1) h(\mathbf e) + (N - g + 1) (2 g +
5)(\log{(d(G))} + 1)
\end{split}
\end{equation*}
car: $(N - g) (2 g + 4) \log{d(G)} + \log{(8 g)} + 3 \log{(d(G))} + (2 g + 2) \log{(d(G) + 1)} + \log{4}$ est majorée
par $(N - g + 1) (2 g + 5) (\log{(d(G))} + 1)$. Pour finir, nous montrons ci-dessous que les formes
${\widetilde{P}}_j$ ($g + 1 \leq j \leq N$) sont toutes de hauteur de Gauss-Weil majorée par $h(G) + d(G)$.
Ces majorations reportées dans l'estimation d'avant pour $\log{(E)} + \log{4}$ aboutissent immédiatement au lemme
\ref{a.9} et ach\`event cette démonstration. \\
Pour $j \in \{g + 1 , \dots , N\}$, montrons qu'on a effectivement
$\widetilde{h}({\widetilde{P}}_j) \leq h(G) + d(G)$. D'apr\`es
\cite{Ph4}3 (page $347$) on a:
$$h\left(\pi_j(G)\right) = h({\widetilde P}_j) + {d°}_{\rm{tot}} {\widetilde P}_j . h({\mathbb P}_g) ,$$
donc
\begin{equation*}
\begin{split}
h({\widetilde P}_j) &= h\left(\pi_j(G)\right) - {d°}_{\rm{tot}} {\widetilde P}_j . h({\mathbb P}_g) \\
&\leq h(G) - {d°}_{\rm{tot}} {\widetilde P}_j . h({\mathbb P}_g) .
\end{split}
\end{equation*}
D'o\`u, d'apr\`es la relation de comparaison (\ref{1.21}) entre
les deux hauteurs $\widetilde h$ et $h$:
\begin{equation*}
\begin{split}
\widetilde{h}({\widetilde P}_j) &\leq h({\widetilde P}_j) + (g + 1) {d°}_{\rm{tot}} {\widetilde P}_j . \log{2} \\
&\leq h(G) - {d°}_{\rm{tot}} {\widetilde P}_j . h({\mathbb P}_g) +
(g + 1) {d°}_{\rm{tot}} {\widetilde P}_j . \log{2} ,
\end{split}
\end{equation*}
i.e
\begin{equation}
\widetilde{h}({\widetilde P}_j) \leq h(G) + {d°}_{\rm{tot}}
{\widetilde P}_j \left((g + 1) \log{2} - h({\mathbb P}_g)\right) .
\label{1.11}
\end{equation}
Posons, provisoirement, pour tout $n \in \mathbb N$:
\begin{equation*}
\begin{split}
\psi(n) &= (n + 1) \log{2} - h({\mathbb P}_n) \\
&= (n + 1) \log{2} - \sum_{m = 1}^{n} \sum_{k = 1}^{m} \frac{1}{2 k} ~~~~ \text{(voir \cite{Ph4}3 page $346$)} \\
&= (n + 1) \log{2} - \sum_{k = 1}^{n} \sum_{m = k}^{n} \frac{1}{2 k} \\
&= (n + 1) \log{2} - \sum_{k = 1}^{n} \frac{n + 1 - k}{2 k} \\
&= (n + 1) \log{2} - \frac{n + 1}{2} \sum_{k = 1}^{n} \frac{1}{k} + \frac{n}{2} \\
&= \frac{n + 1}{2} \left(2 \log{2} + 1 - \sum_{k = 1}^{n}
\frac{1}{k}\right) - \frac{1}{2} .
\end{split}
\end{equation*}
De cette derni\`ere égalité, on vérifie facilement que $\psi(0)$,
$\psi(1)$, $\psi(2)$, $\psi(3)$ et $\psi(4)$ sont inférieurs à $1$
et que $\psi(n) < 0$ pour tout $n \geq 5$. Donc $\psi(n) \leq 1$
pour tout $n \in \mathbb N$, et en particulier $\psi(g) \leq 1$.
Ainsi, de la relation (\ref{1.11}) s'ensuit:
\begin{equation*}
\begin{split}
\widetilde{h}({\widetilde P}_j) &\leq h(G) + {d°}_{\rm{tot}} {\widetilde P}_j \\
&\leq h(G) + d(G) .
\end{split}
\end{equation*}
Ce qui compl\`ete la démonstration. $~~~~\blacksquare$\vspace{1mm}

En appliquant le théor\`eme \ref{a.13} à $G$ et en tenant compte
de la majoration du lemme \ref{a.9}, on obient le théor\`eme
suivant:
\begin{theorem}\label{a.17}
Soit $G$ un groupe algébrique commutatif de dimension $g$ $(g \geq
1)$, défini sur un corps de nombres $K$, plongé dans un espace
projectif ${\mathbb P}_N$$(N \geq 1)$ et d'élément neutre $\mathbf
e$ représenté dans ${\mathbb P}_N$ par un syst\`eme de coordonnées
projectives $\underline e = (e_0 : \dots : e_N)$. On désigne par
$K[\underline X]$$(\underline X = (X_0 , \dots , X_N))$ l'anneau
de coordonnées de ${\mathbb P}_N$, par $\mathfrak G$ l'idéal de
définition de l'adhérence de Zariski $\Gbar$ de $G$ dans
$K[\underline X]$ et par $A := K[\underline X] / \mathfrak{G}$
l'anneau de coordonnées de $\Gbar$. Soient aussi $\mathbf x$ un
point fixé de $G$, $\underline A$ une famille de formes de
$K[\underline X , \underline Y]$$(\underline X := (X_0 , \dots ,
X_N) , \underline Y := (Y_0 , \dots , Y_N))$ représentant
l'addition dans $G \plongement {\mathbb P}_N$ au voisinage de
$\{\mathbf e\} \times \{\mathbf x\}$ et $(c , c') \in {{\mathbb
N}^*}^2$ le bidegré des formes constituant $\underline A$. Soient
enfin $W_0 , \dots , W_g$ des dessous d'escaliers finis de
${\mathbb N}^g$ et $P$ une forme de $A$ non identiquement nulle de
degré $\delta$$(\delta \in {\mathbb N}^*)$. Sous les hypoth\`eses
de normalisation:
\begin{description}
\item[1)] $e_0 = 1$ et l'espace tangent $T_{\mathbf e} G$ de $G$
en $\mathbf e$ a pour équations:
$$X_{g + 1} = \dots = X_N = 0 ;$$
\item[2)] Un point générique $(y_0 : y_1 : \dots : y_N)$ de $G
\hookrightarrow {\mathbb P}_N$ vérifie:
$$\mbox{deg}{\mbox{tr}}_K K(y_0 , \dots , y_g) = \mbox{deg}{\mbox{tr}}_K K(y_0 , \dots , y_N) = g + 1 ,$$
(c'est-à-dire que $y_0 , \dots , y_g$ sont $K$-algébriquement
indépendants);
\end{description}
si $P$ s'annule en $\mathbf x$ avec une multiplicité définie par
le dessous d'escalier $W_0 + \dots + W_g$ (voir les définitions
\ref{a.14}), alors il existe $1 \leq r \leq g$ et une sous-variété
$V$ de $G$ de dimension $d \leq g - r$, contenant $\mathbf x$,
incompl\`etement définie\footnote{selon la terminologie de
\cite{Ph2}.} dans $G$ par des équations de degrés $c' \delta$ et
satisfaisant:
\begin{equation*}
\begin{split}
m_{W_r}(V) \deg{V} &\leq \deg{G} . (c' \delta)^{g - d} \\
\intertext{et}
m_{W_r}(V) h(V) &\leq h(G) (c' \delta)^{g - d} + g \!\!\left[f(N , G , \mathbf e) H\!\left(\!W_0 + \dots + W_{r - 1}\!\right) \right. \\
&~\!+ c \delta \left(h(\mathbf e) + g \log{2}\right) + \widetilde{h}(P) + \delta \widetilde{h}(\underline A) + \delta (2 d°\Asoul + 1) (\log{N} + 1) \\
&~~~~~~~~~~~~~~~~~~~~~\!\left.+ 3 \log\left((\deg{G}) . (c'
\delta)^g + 1\right)\!\right]\!(\deg{G})(c' \delta)^{g - d - 1} ,
\end{split}
\end{equation*}
o\`u
\begin{equation*}
\begin{split}
f(N , G , \mathbf e) &:= 4 (N - g) h(G) + \left[2 (N - g + 1) h(\mathbf e) + 4 (N - g)\right] d(G) \\
&~~~~~~~~~~~~~\!+ (N - g + 1) (2 g + 5) (\log{d(G)} + 1) - 2 (N -
g + 1) h(\mathbf e)
\end{split}
\end{equation*}
et
$$H\!\left(\!W_0 + \dots + W_{r - 1}\!\right) := \max_{I \in W_0 + \dots + W_{r - 1}} \nb{I} .$$
De plus, $P$ s'annule sur $V$ avec la multiplicité définie par le
dessous d'escalier $W_0 + \dots + W_r$.
\end{theorem}
\subsection{Un lemme de zéros explicite dans un groupe produit $G = G_1 \times \dots \times G_p$ plongé dans un
espace multiprojectif $\mathbb P = {\mathbb P}_{n_1} \times \dots
\times {\mathbb P}_{n_p}$:} Soient $G_1 , \dots , G_p$$(p \in
{\mathbb N}^*)$ des groupes algébriques commutatifs définis sur un
m\^eme corps de nombres $K$, irréductibles et plongés
respectivement dans des espaces projectifs ${\mathbb P}_{n_1} ,
\dots , {\mathbb P}_{n_p}$$(n_1 , \dots , n_p \in {\mathbb N}^*)$.
Pour tout $\l = 1 , \dots , p$, nous appelons $g_\l$$(g_\l \geq
1)$ la dimension et ${\mathbf e}_\l$ l'élément neutre de $G_\l$
représenté dans ${\mathbb P}_{n_\l}$ par un syst\`eme de
coordonnées projectives $\esoul_{\l} = (e_{\l 0} : \dots : e_{\l
n_\l})$. Soient aussi pour tout $\l \in \{1 , \dots , p\}$,
${\mathbf x}_\l$ un point fixé de $G_{\l}$ et ${\underline A}_\l$
une famille de formes de $K[{\underline X}_{\l} , {\underline
Y}_{\l}]$$({\underline X}_{\l} = (X_{\l 0} , \dots , X_{\l n_\l})
, \Ysoul_{\l} = (Y_{\l 0} , \dots , Y_{\l n_\l}))$ représentant
l'addition dans $G_\l \plongement {\mathbb P}_{n_\l}$ au voisinage
de $\{{\mathbf e}_\l\} \times \{{\mathbf x}_\l\}$, que nous
supposons constituée de formes de bidegré $(c , c') \in {{\mathbb
N}^*}^2$ indépendant de $\l$. Appelons $G$ le groupe produit $G :=
G_1 \times \dots \times G_p$ plongé dans l'espace multiprojectif
$\mathbb P := {\mathbb P}_{n_1} \times \dots \times {\mathbb
P}_{n_p}$, de dimension $g := g_1 + \dots + g_p$, d'élément neutre
$\mathbf e := ({\mathbf e}_1 , \dots , {\mathbf e}_p)$ et
contenant le point $\mathbf x := ({\mathbf x}_1 , \dots , {\mathbf
x}_p)$. Nous posons ${\underline{\underline X}} := ({\underline
X}_1 , \dots , {\underline X}_p)$ et ${\underline{\underline A}}
:= ({\underline A}_1 , \dots , {\underline A}_p)$ la famille de
formes représentant l'addition dans $G \plongement \mathbb P$ au
voisinage de $\{\mathbf e\} \times \{\mathbf x\}$ et nous
désignons par ${\mathfrak G}_\l$$(\l = 1 , \dots , p)$ l'idéal de
définition dans $K[{\underline X}_\l]$ de l'adhérence de Zariski
$\Gbar_\l$ de $G_\l$ et par $A_\l := K[\Xsoul_\l] / {\mathfrak
G}_\l$$(\l = 1 , \dots , p)$ l'anneau de coordonnées de
$\Gbar_\l$. Nous désignons finalement par $\mathfrak G$ l'idéal de
définition dans $K[\Xsoull]$ de l'adhérence de Zariski $\Gbar$ de
$G$ et par $A := K[\Xsoull] / \mathfrak{G}$ l'anneau de
coordonnées de $\Gbar$.

Par ailleurs, nous posons pour tous réels strictement positifs
donnés $\epsilon , \delta_1 , \dots , \delta_p$ et pour
$\underline{\delta} := (\delta_1 , \dots , \delta_p)$:
$$W(\underline{\delta} , \epsilon) := \left\{(\alphasoul_1 , \dots , \alphasoul_p) \in {\mathbb N}^{g_1 + \dots + g_p} / \frac{\nb{\alphasoul_1}}{\delta_1} + \dots + \frac{\nb{\alphasoul_p}}{\delta_p} < \epsilon\right\}$$
qui est un dessous d'escalier fini de ${\mathbb N}^g$.

En se mettant dans la situation décrite par les hypoth\`eses de
normalisation du §$3$ pour chaque groupe $G_\l$$(\l = 1 , \dots ,
p)$ et ${\mathbf e}_\l \in G_\l$, le groupe $G_{\l} \plongement
{\mathbb P}_{n_\l}$ est incompl\`etement défini par des équations
de la forme: ${\widetilde{P}}_{\l j}(X_{\l 0} , \dots , X_{\l g} ,
X_{\l j}) = 0$ $(g_\l + 1 \leq j \leq n_\l)$, o\`u les
${\widetilde{P}}_{\l j}$$(g_{\l} + 1 \leq j \leq n_\l)$ sont des
formes, non identiquement nulles de $K[X_{\l 0} , \dots , X_{\l
g_\l} , X_j]$ ayant chacune un de ces coffecients égal à $1$.
Posons aussi $P_{\l j}$$(\l = 1 , \dots , p ; j = g_{\l} + 1 ,
\dots , n_\l)$ le polyn\^ome obtenu en substituant dans
${\widetilde{P}}_{\l j}$, l'indétérminée  $X_0$ par $1$. On a
d'apr\`es le §$3$, une paramétrisation de $G_\l$ au voisinage de
son élément neutre ${\mathbf e}_\l$ utilisant comme param\`etres
$\Tsoul_\l = (T_{\l 1} , \dots , T_{\l g_\l}) := (X_{\l 1} - e_{\l
1} , \dots , X_{\l g_\l} - e_{\l g_\l})$ et représentée par un
monomorphisme $\varphi_\l$ de $A_\l$ dans $K[[\Tsoul_\l]]$. De
plus, en écrivant pour tout $i \in \{0 , \dots , n_\l\}$:
$$\varphi_\l(X_{\l i}) = \sum_{I_\l \in {\mathbb N}^{g_\l}} a_{I_\l}^{(\l i)} \Tsoul_{\l}^{I_\l} ,$$
les nombres $a_{I_\l}^{(\l i)}$$(I_\l \in {\mathbb N}^{g_\l} , i =
0 , \dots , n_\l)$ de $K$ vérifient, pour toute place $v$ de $K$
(d'apr\`es le lemme \ref{a.6}), la relation:
\begin{equation}
{\nb{a_{I_\l}^{(\l i)}}}_v \leq {F_{\l , v}(n_\l , G_\l ,
\esoul_\l)}^{\nb{I_\l}} H_v(\esoul_\l) ~~~~~~ (I_\l \in {\mathbb
N}^{g_\l} , i = 0 , \dots , n_\l) \label{1.12}
\end{equation}
o\`u les $F_{\l , v}$$(v \in M_K)$ sont les expressions
(supérieures ou égales à $1$) dépendant de $n_\l , G_\l$ et
$\esoul_\l$ définies par:
\begin{align}
F_{\l , v} &:= {H_v(\esoul_\l)}^{2 (d(G_\l) - 1)} \prod_{j = g_\l + 1}^{n_\l} \!\!\!\left({H_v({\widetilde P}_{\l j})}^2 {\max\!\left\{1 , {\nb{\frac{1}{\frac{\partial P_{\l j}}{\partial X_{\l j}}(e_{\l 1} , \dots , e_{\l g_\l} , e_{\l j})}}}_v\right\}}^{\!\!2}\right) \notag \\
\intertext{si $v$ est finie et}
F_{\l , v} &:= 8 g_\l {d(G_\l)}^3 (d(G_\l) + 1)^{2 (g_\l + 1)} . {H_v(\esoul_\l)}^{2 (d(G_\l) - 1)} \notag \\
&~~~~~~~~~~~~~~~~~~~~\!\!\!\!\times \prod_{j = g_\l + 1}^{n_\l}
\!\!\!\left({H_v({\widetilde P}_{\l j})}^2 {\max\!\left\{1 ,
{\nb{\frac{1}{\frac{\partial P_{\l j}}{\partial X_{\l j}}(e_{\l 1}
, \dots , e_{\l g_\l} , e_{\l j})}}}_v\right\}}^{\!\!2}\right)
\notag
\end{align}
si $v$ est infinie. \\
On pose aussi:
$$F_{\l} := \prod_{v \in M_K} \!\!F_{\l , v}^{\frac{[K_v : {\mathbb Q}_v]}{[K : \mathbb Q]}} .$$
Les paramétrisations de $G_1 , \dots , G_p$ au voisinage de
${\mathbf e}_1 , \dots , {\mathbf e}_p$ respectivement induisent
évidemment une paramétrisation du groupe produit $G$ au voisinage
de son élément neutre $\mathbf e$, utilisant comme param\`etres
$\Tsoull := (\Tsoul_1 , \dots , \Tsoul_p) = (T_{1 1} , \dots ,
T_{1 g_1} ,$\\$ \dots , T_{p 1} , \dots , T_{p g_p})$ et c'est
cette paramétrisation qu'on utilise.

Notons enfin que les définitions \ref{a.14} concernant
l'annulation d'une forme de $A$ sur un sous-ensemble algébrique
$E$ de $G$ avec une multiplicité définie par un dessous d'escalier
$W$ de ${\mathbb N}^g$ s'étendent d'une mani\`ere évidente à notre
cas $G \plongement {\mathbb P}_{n_1} \times \dots \times {\mathbb
P}_{n_p}$. On a le théor\`eme suivant:
\begin{theorem}[Notre lemme de Roth]\label{a.19}
Sous toutes les hypoth\`eses précédentes, soient $\epsilon$ un
réel strictement positif et $P \in A$ une forme non identiquement
nulle de multidegré $\underline{\delta} = (\delta_1 , \dots ,
\delta_p)$. On suppose $\epsilon \leq 1$, $\delta_\l \geq g_\l +
1$ (pour $\l = 1 , \dots , p$) et:
\begin{equation}
\frac{\delta_\l}{\delta_{\l + 1}} > \left(\frac{p
c'}{\epsilon}\right)^g d(G_1) \dots d(G_p) , ~~~~~~ \l = 1 , \dots
, p - 1 . \label{1.33}
\end{equation}
Alors, si $P$ s'annule en $\mathbf x$ avec une multiplicité
définie par le dessous d'escalier $W(\underline{\delta} , g
\epsilon)$ de ${\mathbb N}^g$, il existe une sous-variété propre
$V$ de $G$ définie sur $K$, contenant $\mathbf x$, dont toutes les
composantes irréductibles sur $\Kbar$ sont des sous-variétés
produit dans $\mathbb P$, on a $V = \bigcup_{\sigma \in
{\rm{Gal}}(\Kbar / K)}\sigma(\widetilde{V})$ o\`u $\widetilde{V} =
V_1 \times \dots \times V_p$ et en posant $D :=
\frac{d(V)}{d(\widetilde{V})}$ (qui est -en d'autres termes- le
nombre de composantes irréductibles de $V$ sur $\Kbar$), on a:
\begin{align}
D d(V_1) \dots d(V_p) &\leq \left(\frac{g c'}{\epsilon}\right)^{g
- \dim{V}} d(G_1) \dots d(G_p) \notag \intertext{et} D
\frac{d(V_1)}{d(G_1)} \dots \frac{d(V_p)}{d(G_p)} \sum_{\l =
1}^{p} \delta_\l \frac{h(V_\l)}{d(V_\l)} &\leq \left(\frac{(g + 1)
c'}{\epsilon}\right)^{g - \dim{V}} \!\left[\widetilde{h}(P) +
\sum_{\l = 1}^{p} R_\l \delta_\l + S\right] \notag
\end{align}
o\`u pour tout $\l = 1 , \dots , p$, on a noté $R_\l$ l'expression
dépendant de $n_\l , G_\l , {\mathbf e}_\l$ et $\Asoul_\l$ définie
par:
\begin{equation*}
\begin{split}
R_\l &:= \frac{1}{g_\l + 1} \frac{h(G_\l)}{d(G_\l)} + (g - 1) \epsilon \left\{4(n_\l - g_\l) h(G_\l) \right. \\
&~~~~~~~~~~~~\!\!+ \left[2(n_\l - g_\l + 1) h({\mathbf e}_\l) + 4(n_\l - g_\l)\right] d(G_\l) \\
&~~~~~~~~~~~~\!\!+\left.(n_\l - g_\l + 1) (2 g_\l + 5) (\log{(d(G_\l))} + 1) - 2(n_\l - g_\l + 1) h({\mathbf e}_\l)\right\} \\
&~~~~~~~~~~~~\!\!+ c\left(h({\mathbf e}_\l) + g_\l \log{2}\right)
+ \widetilde{h}(\Asoul_\l) + (3 c + 3 c' + 4) \log{(n_\l + 1)}
\end{split}
\end{equation*}
et $S$ l'expression:
$$S := 3 \sum_{\l = 1}^{p} \log{\left(d(G_\l)\right)} + g(3 \log{p} + 3 \log{c'} + c \log{2}) + 2 c + 2 c' + 4 .$$
De plus, $P$ s'annule sur $V$ avec la multiplicité définie par le
dessous d'escalier $W(\underline{\delta} , \epsilon)$.
\end{theorem}
{\bf Démonstration.---} Afin de pouvoir utiliser le théor\`eme
\ref{a.13}, on plonge $\mathbb{P}$ dans un espace projectif
${\mathbb{P}}_N$ par le plongement de Segre-Veronese
$$
\begin{array}{rcl}
\rho\!: ~~~~~~~\mathbb{P}~~~~~~~~~ & \plongement & ~~~~~~{\mathbb P}_N \\
(\Xsoul_1 , \dots , \Xsoul_p) & \mapsto &
{\left(\Xsoul_{1}^{\alphasoul_1} \dots
\Xsoul_{p}^{\alphasoul_p}\right)}_{\nb{\alphasoul_1} = \delta_1 ,
\dots , \nb{\alphasoul_p} = \delta_p}
\end{array} ,
$$
$N$ est alors égal au cardinal de l'ensemble
$$\mathcal{N} := \left\{\alphasoull = (\alphasoul_1 , \dots , \alphasoul_p) \in {\mathbb N}^{n_1 + 1} \times \dots \times {\mathbb N}^{n_p + 1} / \nb{\alphasoul_1} = \delta_1 , \dots , \nb{\alphasoul_p} = \delta_p\right\} ,$$
donc:
$$N = \binom{n_1 + \delta_1}{\delta_1} \dots \binom{n_p + \delta_p}{\delta_p} .$$
Appelons $\Zsoul = {\left(Z^{(\alphasoull)}\right)}_{\alphasoull
\in \mathcal{N}}$ les coordonnées de ${\mathbb P}_N$. A toute
forme homog\`ene $Q$ sur ${\mathbb P}_N$ (i.e $Q \in K[\Zsoul]$)
de degré $\delta \in \mathbb N$ est associée la forme
$Q\left(\Xsoull^{\alphasoull} , \alphasoull \in
\mathcal{N}\right)$ (obtenue en substituant
$\Xsoull^{\alphasoull}$ dans $Q$ à chaque coordonnée
$Z^{(\alphasoull)}$$(\alphasoull \in \mathcal{N})$ de ${\mathbb
P}_N$), qui est une forme multihomog\`ene sur $\mathbb P$ de degré
$(\delta \delta_1 , \dots , \delta \delta_p)$. Ainsi, une forme
linéaire sur ${\mathbb P}_N$ représente une forme multihomog\`ene
de degré $(\delta_1 , \dots , \delta_p)$ sur $\mathbb P$. Et en
particulier notre forme $P \in A$ du théor\`eme \ref{a.19} est
représentée par une forme linéaire $L$ sur ${\mathbb P}_N$.

On obtient à partir des monomorphismes $\varphi_1 , \dots ,
\varphi_p$ de paramétrisation des groupes $G_1 , \dots , G_p$
respectivement, au voisinage de leurs origines, un monomorphisme
$\varphi$ de paramétrisation de $G \plongement {\mathbb P}_N$ au
voisinage de $\mathbf e$ défini par:
$$
\begin{array}{rcl}
\varphi\: K[\Zsoul] / \mathfrak{I}\left(\rho(G)\right) & \rightarrow & ~~~~~~~~~~~~~~~K[[\Tsoull]] \\
Q~~~~~~ & \mapsto & \varphi(Q) := \left(\varphi_1 \otimes \dots
\otimes \varphi_p\right)\!\!\left(Q(\Xsoull^{\alphasoull} ,
\alphasoull \in \mathcal{N})\right)
\end{array} .
$$
En posant (comme dans le corollaire \ref{a.15}) pour tout $\l \in
\{1 , \dots , p\}$ et tout $\alphasoul_{\l} \in {\mathbb
N}^{n_{\l} + 1}$:
$$\varphi_{\l}\left(\Xsoul_{\l}^{\alphasoul_{\l}}\right) = \sum_{I_{\l} \in {\mathbb N}^{g_{\l}}} {\mathcal C}_{\l}\!\left(\alphasoul_{\l} , I_{\l}\right) {\Tsoul}_{\l}^{I_{\l}} ~~~~ \left({\mathcal C}_{\l}\!\left(\alphasoul_{\l} , I_{\l}\right) \in K\right) ,$$
on a pour tout $\alphasoull = (\alphasoul_1 , \dots ,
\alphasoul_p) \in \mathcal{N}$:
\begin{equation*}
\begin{split}
\varphi\left(Z^{(\alphasoull)}\right) &= \varphi_1\left(\Xsoul_{1}^{\alphasoul_1}\right) \dots \varphi_p\left(\Xsoul_{p}^{\alphasoul_p}\right) \\
&= \prod_{\l = 1}^{p} \sum_{I_{\l} \in {\mathbb N}^{g_{\l}}} \!\!{\mathcal C}_{\l}\!\left(\alphasoul_{\l} , I_{\l}\right) {\Tsoul}_{\l}^{I_{\l}} \\
&= \sum_{I = (I_1 , \dots , I_p) \in {\mathbb N}^{g_1} \times
\dots \times {\mathbb N}^{g_p} = {\mathbb N}^g} \!\!\!\!{\mathcal
C}_1\!\left(\alphasoul_1 , I_1\right) \dots {\mathcal
C}_p\!\left(\alphasoul_p , I_p\right) \Tsoull^I ,
\end{split}
\end{equation*}
c'est-à-dire:
\begin{equation}
\varphi\left(Z^{(\alphasoull)}\right) = \sum_{I \in {\mathbb N}^g}
a_{I}^{(\alphasoull)} \Tsoull^I \label{1.31}
\end{equation}
avec les $a_{I}^{(\alphasoull)}$$(I \in {\mathbb N}^g)$ sont les
nombres de $K$ définis par:
$$a_{I}^{(\alphasoull)} := {\mathcal C}_1\!\left(\alphasoul_1 , I_1\right) \dots {\mathcal C}_p\!\left(\alphasoul_p , I_p\right)$$
quand $I = (I_1 , \dots , I_p) \in {\mathbb N}^{g_1} \times \dots \times {\mathbb N}^{g_p} = {\mathbb N}^g$. \\

Pour pouvoir appliquer le théor\`eme \ref{a.13} à $G \plongement {\mathbb P}_N$, nous allons montrer dans ce qui suit que ces nombres $a_{I}^{(\alphasoull)}$$(\alphasoull \in \mathcal{N} , I \in {\mathbb N}^g)$ satisfont bien la relation (\ref{1.25}) du §$4.1$ pour un certain $s > 0$, certaines expressions $E_{\l , v} \geq 1$$(\l = 1 , \dots , p , v \in M_K)$ et certaines fonctions additives $t_{\l}$$(\l = 1 , \dots , p)$ de ${\mathbb N}^g$ dans ${\mathbb R}^+$, à détérminer. \\
Pour tout $\l = 1 , \dots , p$, la paramétrisation considérée pour
le groupe $G_{\l} \plongement {\mathbb P}_{n_{\l}}$ au voisinage
de ${\mathbf e}_{\l}$ satisfait clairement les hypoth\`eses du
§$4.1$ avec $p$ remplacé par $1$, $n_1$ remplacé par $n_{\l}$,
$g_1$ remplacé par $g_{\l}$ (ou m\^eme quelconque), $\delta_1$
remplacé par $1$, $s$ remplacé par $1$, ${\underline e}_1$
remplacé par ${\underline e}_{\l}$, $t_1$ remplacée par la
longueur d'un $g_{\l}$-uplet de ${\mathbb N}^{g_{\l}}$, $\Tsoul$
remplacé par $\Tsoul_{\l}$ et $E_{1 , v}$$(v \in M_K)$ remplacées
par les $F_{\l , v}$ de la relation (\ref{1.12}). Ce qui nous
permet d'appliquer le corollaire \ref{a.15} à $G_{\l} \plongement
{\mathbb P}_{n_{\l}}$ et on obtient: pour tout $\alphasoul_{\l}
\in {\mathbb N}^{n_{\l} + 1}$ tel que $\nb{\alphasoul_{\l}} =
\delta_{\l}$, tout $I_{\l} \in {\mathbb N}^{g_{\l}}$ et toute
place $v$ de $K$:
\begin{align}
{\nb{{\mathcal C}_{\l}\!\left(\alphasoul_{\l} , I_{\l}\right)}}_v &\leq 2^{\nb{I_{\l}} + g_{\l} \delta_{\l}} F_{\l , v}^{\nb{I_{\l}}} {H_v(\esoul_{\l})}^{\delta_{\l}} \notag \\
\intertext{si $v$ est infinie et:} {\nb{{\mathcal
C}_{\l}\!\left(\alphasoul_{\l} , I_{\l}\right)}}_v &\leq F_{\l ,
v}^{\nb{I_{\l}}} {H_v(\esoul_{\l})}^{\delta_{\l}} \notag
\end{align}
si $v$ est finie. \\
Définissons pour tout $\l \in \{1 , \dots , p\}$ et toute place
$v$ de $K$:
$$s := 2 ,~\! s_v := \!\begin{cases} 1 & \text{si $v$ est finie} \\ 2 & \text{si $v$ est infinie} \end{cases} \!,~\! E_{\l , v} := s_v F_{\l , v} ,~\! E_{\l} := \!\!\prod_{v \in M_K} \!\!E_{\l , v}^{\frac{[K_v : {\mathbb Q}_v]}{[K : \mathbb Q]}} \!\!= 2 F_{\l}$$
$$
\begin{array}{rcl}
\!\!\mbox{et}~~~~~~~~~~~ t_{\l}\: ~~~~~~~~~~~~~~~~~~~{\mathbb N}^g  ~~~~~~~~~~~~~~~~~~& \rightarrow & ~~~~~{\mathbb R}^+ \\
I = (I_1 , \dots , I_g) \in {\mathbb N}^{g_1} \times \dots \times
{\mathbb N}^{g_p} = {\mathbb N}^g & \mapsto & t_{\l}(I) :=
\nb{I_{\l}}
\end{array}~~~~
$$
(les $t_{\l}$$(\l = 1 , \dots , p)$ sont clairement des fonctions additives satisfaisant $t_1(I) + \dots + t_p(I) = \nb{I}$ pour tout $I \in {\mathbb N}^g$). \\
Avec toutes ces définitions, les relations précédentes estimant
les valeurs absolues $v$-adiques des nombres ${\mathcal
C}_{\l}\!\left(\alphasoul_{\l} , I_{\l}\right)$$(\l ,
\alphasoul_{\l} , I_{\l})$ s'écrivent:
$${\nb{{\mathcal C}_{\l}\!\left(\alphasoul_{\l} , I_{\l}\right)}}_v \leq s_{v}^{g_{\l} \delta_{\l}} E_{\ell , v}^{\nb{I_{\l}}} {H_v(\esoul_{\l})}^{\delta_{\l}} ~~~~~~ \text{(pour tous $\l , \alphasoul_{\l} , I_{\l}$ et $v$)} .$$
D'o\`u, pour tout $\alphasoull \in \mathcal{N}$, tout $I \in
{\mathbb N}^g$ et toute place $v$ de $K$:
\begin{equation}
{\nb{a_{I}^{(\alphasoull)}}}_v \leq s_{v}^{g_1 \delta_1 + \dots +
g_p \delta_p} E_{1 , v}^{t_1(I)} \dots E_{p , v}^{t_p(I)}
{H_v(\esoul_1)}^{\delta_1} \dots {H_v(\esoul_p)}^{\delta_p} .
\label{1.13}
\end{equation}
On a ainsi vérifié (relations (\ref{1.31}) et (\ref{1.13})) les hypoth\`eses du §$4.1$ pour la paramétrisation $\varphi$ de $G \plongement {\mathbb P}_N$ au voisinage de $\mathbf e$, ce qui nous autorise enfin à utiliser le théor\`eme \ref{a.13} pour $G \plongement {\mathbb P}_N$. \\
Dans les hypoth\`eses du théor\`eme \ref{a.19}, $P$ s'annule en
$\mathbf x$ avec la multiplicité définie par le dessous d'escalier
$W(\underline{\delta} , g \epsilon)$. Pour se ramener à la
situation du théor\`eme \ref{a.13}, on pose $W_0 , \dots , W_g$
les dessous d'escaliers de ${\mathbb N}^g$:
$$W_0 = \{0\} ~~\mbox{et}~~ W_1 = \dots = W_g = W(\underline{\delta} , \epsilon),$$
(on a bien $W_0 + \dots + W_g = W(\underline{\delta} , g
\epsilon)$). Dans cette application du théor\`eme \ref{a.13}, la
forme linéaire $L$ sur ${\mathbb P}_N$ représentant notre forme
$P$ sur $\mathbb P$ a les m\^emes coefficients que $P$ donc aussi
la m\^eme hauteur de Gauss-Weil que $P$. Par ailleurs, comme la
famille de formes $\Asoull = (\Asoul_1 , \dots , \Asoul_p)$
représentant l'addition dans $G \plongement \mathbb P$ au
voisinage de $\{\mathbf e\} \times \{\mathbf x\}$ est constituée
de formes, toutes de m\^eme bidegré $(c , c')$, l'addition dans $G
\plongement {\mathbb P}_N$ au voisinage de $\{\mathbf e\} \times
\{\mathbf x\}$ peut \^etre représentée par une famille
$\underline{B} = \left(B^{(\alphasoull)} , \alphasoull \in
\mathcal N\right)$ de formes de $K[\Zsoul]$ de bidegré $(c , c')$
et tel que $\underline{B}$ a\"{\i}t la m\^eme hauteur de
Gauss-Weil que la famille de formes $\Asoull^{\alphasoull}$,
$\alphasoull \in \mathcal N$ de $K[\Xsoull]$ (en fait chaque forme
$B^{(\alphasoull)}$$(\alphasoull \in \mathcal N)$ représente
simplement la forme $\Asoull^{\alphasoull}$ sur ${\mathbb P}_N$).
On a, gr\^ace à un calcul standard sur les hauteurs:
\begin{equation*}
\begin{split}
\widetilde{h}(\underline{B}) &= \widetilde{h}\left(\Asoull^{\alphasoull} , \alphasoull \in \mathcal N\right) \\
&\!\!\!\!\!\!\!\!\leq \delta_1 \widetilde{h}(\Asoul_1) + \dots +
\delta_p \widetilde{h}(\Asoul_p) + (c + c') \left[\delta_1
\log{(n_1 + 1)} + \dots + \delta_p \log{(n_p + 1)}\right] .
\end{split}
\end{equation*}
La conclusion du théor\`eme \ref{a.13} affirme alors l'existence
d'une sous-variété propre $V$ de $G$ contenant $\mathbf x$ et
satisfaisant:
\begin{equation}
m_{W_r}(V) d_{{\mathbb P}_N}\!\left(\rho(V)\right) \leq
d_{{\mathbb P}_N}\!\left(\rho(G)\right) {c'}^{g - \dim{V}}
\label{1.39}
\end{equation}
et
\begin{equation}
\begin{split}
m_{W_r}(V) h_{{\mathbb P}_N}\!\left(\rho(V)\right) &~\leq~
h_{{\mathbb P}_N}\!\left(\rho(G)\right) {c'}^{g - \dim{V}}\\
&\!\!\!\!\!\!\!\!\!\!\!\!\!\!\!\!\!\!\!\!\!\!\!\!\!\!\!\!\!\!\!\!\!\!\!\!\!\!\!\!\!\!\!\!\!\!\!\!\!\!\!\!+
g \!\!\left[\sum_{\ell = 1}^{p} \left(\log{(E_{\ell})} +
\log{2}\right) t_{\ell}\!\left(\!W_0 + \dots + W_{r -
1}\!\right)\right. + c \!\left(\sum_{\ell = 1}^{p} \delta_{\ell}
\!\left(h({\mathbf e}_{\ell})
+ g_{\ell} \log{2}\!\right) + g \log{2}\!\right)\\
&\!\!\!\!\!\!\!\!\!\!\!\!\!\!\!\!\!\!\!\!\!\!\!\!\!\!\!\!\!\!\!\!\!\!\!\!\!\!\!\!\!\!\!\!\!\!\!\!\!\!\!\left.
+ \widetilde{h}(L) + \widetilde{h}(\underline B) + (2 c + 2 c' +
1) (\log{N} + 1) + 3 \log\!\left(d_{{\mathbb
P}_N}\!\left(\rho(G)\right) . {c'}^g +
1\right)\!\phantom{\sum_{\ell = 1}^{p} \left(\log{E_{\ell}} +
\log{(c \delta)}\right) t_{\ell}\!\left(\!W_0 + \dots + W_{r -
1}\!\right)}\!\!\!\!\!\!\!\!\!\!\!\!\!\!\!\!\!\!\!\!\!\!\!\!\!\!\!\!\!\!\!\!\!\!\!\!\!\!\!\!\!\!\!\!\!\!\!\!\!\!\!\!\!\!\!\!\!\!\!\!\!\!\!\!\!\!\!
\!\!\!\!\!\!\!\!\!\!\!\!\!\!\!\!\!\!\!\!\!\!\!\!\!\!\!\!\!\!\!\!\!\!\!\!\!\!\!\!\!\!\!\right]\!d_{{\mathbb
P}_N}\!\left(\rho(G)\right){c'}^{g - \dim{V} - 1}
\end{split} \label{1.40}
\end{equation}
pour un certain entier $1 \leq r \leq g - \dim{V}$. \\
De plus, $P$ s'annule sur $V$ avec la multiplicité définie par le dessous d'escalier $W_0 + \dots + W_r \supset W(\underline{\delta} , \epsilon)$ (donc $P$ s'annule -à fortiori- sur $V$ avec la multiplicité définie par le dessous d'escalier $W(\underline{\delta} , \epsilon)$). \\
On obtient les inégalités du théor\`eme \ref{a.19} à partir des
relations (\ref{1.39}) et (\ref{1.40}), mais avant cela nous
montrons d'abord que les composantes irréductibles de $V$ sur
$\Kbar$ sont toutes des sous-variétés produit de $\mathbb P$. Soit
$\widetilde{V}$ une composante irréductible de $V$ sur $\Kbar$.
Pour tout $\l = 1 , \dots , p$, soient $\pi_{\l}$ la projection de
$\mathbb P$ sur ${\mathbb P}_{n_{\l}}$, $\mathbf{y} = ({\mathbf
y}_1 , \dots , {\mathbf y}_p)$ un point suffisamment général de
$\widetilde{V}$ et $X_{\l}$ le sous-ensemble algébrique de
$\mathbb P$:
$$X_{\l} := \widetilde{V} \cap \left(\{{\mathbf y}_1\} \times \dots \times \{{\mathbf y}_{\l - 1}\} \times
{\mathbb P}_{n_{\l}} \times \{{\mathbf y}_{\l + 1}\} \times \dots
\times \{{\mathbf y}_p\}\right) .$$ Il est clair qu'on a pour tout
$\l = 1 , \dots , p$: $\dim{X_{\l}} \leq
\dim{\pi_{\l}(\widetilde{V})}$. De plus, en supposant que
$\widetilde{V}$ n'est pas le produit de ses projections sur les
différents ${\mathbb P}_{n_{\l}}$, il existe $1 \leq \l < p$ tel
que l'on a\"{\i}t:
$$\dim{X_{\l}} < \dim{\pi_{\l}(\widetilde{V})} .$$
On choisit $\l$ minimal ayant cette propriété, donc pour tout $j =
1 , \dots , \l - 1$, on a:
$$\dim{X_j} = \dim{\pi_j(\widetilde{V})} ,$$
qui est équivaut à:
$$\{{\mathbf y}_1\} \times \dots \times \{{\mathbf y}_{j - 1}\} \times \pi_j(\widetilde{V}) \times \{{\mathbf y}_{j + 1}\} \times \dots \times \{{\mathbf y}_p\} \subset \widetilde{V} .$$
Pour $j = 1 , \dots , \l$, on pose $\alpha_j :=
\dim{\pi_j(\widetilde{V})}$ et $\beta_j := g_j - \dim{X_j}$, alors
on a $\alpha_j + \beta_j = g_j$ pour tout $j \in \{1 , \dots , \l
- 1\}$ et $\alpha_{\l} + \beta_{\l} > g_{\l}$. On compl\`ete
$\alpha_1 , \dots , \alpha_{\l}$ et $\beta_1 , \dots , \beta_{\l}$
pour avoir des $p$-uplets $\underline{\alpha}$ et
$\underline{\beta}$ dans ${\mathbb N}^p$ satisfaisant
$\nb{\alphasoul} = \dim{V}$, $\nb{\underline{\beta}} = g -
\dim{V}$, $d_{\alpha}(V) \neq 0$, $m_{W_r}(V) \geq
\delta_{1}^{\beta_1} \dots \delta_{p}^{\beta_p}
\epsilon^{\nb{\underline{\beta}}} = \delta_{1}^{\beta_1} \dots
\delta_{p}^{\beta_p} \epsilon^{g - \dim{V}}$ et $\alphasoul +
\betasoul$ maximal pour l'ordre lexicographique\footnote{Celà veut
dire qu'on prend $\underline{\alpha}$ et $\underline{\beta}$ pour
que l'on a\"{\i}t: $\sum_{j = 1}^{i} (\alpha_j + \beta_j) \geq
\sum_{j = 1}^{i} g_j$ pour tout $i = 1 , \dots , p$.}\label{B}
(ceci est possible d'apr\`es la minimalité de $\l$). Comme
maintenant $\nb{\alphasoul} + \nb{\underline{\beta}} = g = g_1 +
\dots + g_p$, $\alpha_1 + \beta_1 = g_1 , \dots , \alpha_{\l - 1}
+ \beta_{\l - 1} = g_{\l - 1}$ et $\alpha_{\l} + \beta_{\l} >
g_{\l}$, il doit exister un indice $\l < k \leq p$ pour lequel on
a $\alpha_k + \beta_k < g_k$. Prenons $k$ minimal pour cette
propriété et posons pour tout $i = 1 , \dots , p$:
$$s_i := \begin{cases}
\alpha_i + \beta_i - g_i & \text{si $i \neq  \ell$ et $i \neq k$} \\
\alpha_i + \beta_i - g_i - 1 & \text{si $i = \ell$} \\
\alpha_i + \beta_i - g_i + 1 & \text{si $i = k$}
\end{cases} .$$
On vérifie aisemment gr\^ace au choix de $\underline{\alpha}$ et
$\underline{\beta}$ maximaux pour l'ordre lexicographique (voir la
note d'en bas de la page \pageref{B}) et gr\^ace à la minimalité
de $k$, qu'on a $\sum_{j = 1}^{i} s_j \geq 0$ pour tout $i = 1 ,
\dots , p$. Par suite on a:
\begin{equation*}
\begin{split}
\delta_{1}^{\alpha_1 + \beta_1} \dots \delta_{p}^{\alpha_p + \beta_p} &= \delta_{1}^{g_1} \dots \delta_{p}^{g_p} \frac{\delta_{\l}}{\delta_k} \prod_{i = 1}^{p} \delta_{i}^{s_i} \\
&= \delta_{1}^{g_1} \dots \delta_{p}^{g_p} \frac{\delta_{\l}}{\delta_k} \prod_{i = 1}^{p} \left(\frac{\delta_i}{\delta_{i + 1}}\right)^{\sum_{j = 1}^{i} s_j} ~~~~ \text{(avec $\delta_{p + 1} = 1$)} \\
&\geq \delta_{1}^{g_1} \dots \delta_{p}^{g_p}
\frac{\delta_{\l}}{\delta_{\l + 1}} ~~~~ \text{(car les $\delta_i$
décroissent)} ,
\end{split}
\end{equation*}
et puis:
\begin{equation*}
\begin{split}
\delta_{1}^{g_1} \dots \delta_{p}^{g_p} \frac{\delta_{\l}}{\delta_{\l + 1}} \epsilon^{g - \dim{V}} &\leq \delta_{1}^{\alpha_1} \dots \delta_{p}^{\alpha_p} . \delta_{1}^{\beta_1} \dots \delta_{p}^{\beta_p} \epsilon^{g - \dim{V}} \\
&\leq d_{{\mathbb P}_N}\!\left(\rho(V)\right) m_{W_r}(V) \\
&\leq d_{{\mathbb P}_N}\!\left(\rho(G)\right) {c'}^{g - \dim{V}} ~~~~\text{(d'apr\`es (\ref{1.39}))} \\
&\leq \binom{g}{g_1 \dots g_p} d(G_1) \dots d(G_p) \delta_{1}^{g_1} \dots {\delta}_{p}^{g_p} {c'}^{g - \dim{V}} \\
&\leq \delta_{1}^{g_1} \dots {\delta}_{p}^{g_p} d(G_1) \dots
d(G_p) (p c')^g
\end{split}
\end{equation*}
o\`u la deuxi\`eme inégalité de cette série d'inégalités vient du
fait qu'on a:\\ $d_{{\mathbb P}_N}\!\left(\rho(V)\right) \geq
\delta_{1}^{\alpha_1} \dots {\delta}_{p}^{\alpha_p}$ car
$d_{\alphasoul}(V) \neq 0$ (en effet, on a:\\ $\displaystyle
d_{{\mathbb P}_N}\!\left(\rho(V)\right) = (\dim{V})!
\!\!\!\!\sum_{\underline{\gamma} \in {\mathbb N}^{p} /
\nb{\underline{\gamma}} = \dim{V}}
\!\!\!\!d_{\underline{\gamma}}(V) \left(\delta_{1}^{\gamma_1} /
\gamma_1!\right) \dots \left(\delta_{p}^{\gamma_p} /
\gamma_p!\right) \geq \delta_{1}^{\alpha_1} \dots
{\delta}_{p}^{\alpha_p}$). D'o\`u:
$$\frac{\delta_{\l}}{\delta_{\l + 1}} \leq \left(\frac{p c'}{\epsilon}\right)^{\!\!\!g} d(G_1) \dots d(G_p) ,$$
ce qui contredit notre hypoth\`ese (\ref{1.33}) et montre ainsi que $\widetilde{V}$ est bien une sous-variété produit de $\mathbb P$. \\
Maintenant il ne reste qu'à déduire des inégalités (\ref{1.39}) et
(\ref{1.40}), celles du théor\`eme \ref{a.19}. Pour cela, nous
allons estimer chacune des quantités $d_{{\mathbb
P}_N}\!\left(\rho(G)\right) , h_{{\mathbb
P}_N}\!\left(\rho(G)\right)$, \\ $d_{{\mathbb
P}_N}\!\left(\rho(V)\right) , h_{{\mathbb
P}_N}\!\left(\rho(V)\right) , m_{W_r}(V) , \log{(E_{\l})} +
\log{2}$ $(\l = 1 , \dots , p)$, $\log{N} + 1$ et \\
$\log{\left(d_{{\mathbb P}_N}\!\left(\rho(G)\right) {c'}^g +
1\right)}$ indépendamment du plongement $\rho$, de $N$ et de
${\mathbb P}_N$. En désignant par $H_g(G ; \underline{\delta})$ et
$H_a(G ; \underline{\delta})$ (resp $H_g(V ; \underline{\delta})$
et $H_a(V ; \underline{\delta})$) les polyn\^omes multidegré et
multihauteur de $G$ (resp de $V$) définis dans \cite{Ph2} (§$3$),
on a:
$$d_{{\mathbb P}_N}\!\left(\rho(G)\right) = H_g(G ; \underline{\delta}) = \binom{g}{g_1 \dots g_p} d(G_1) \dots d(G_p) \delta_{1}^{g_1} \dots \delta_{p}^{g_p}$$ (car $G = G_1 \times \dots \times G_p$) et
\begin{equation*}
\begin{split}
h_{{\mathbb P}_N}\!\left(\rho(G)\right) &= H_a(G ; \underline{\delta}) \\
&= \binom{g}{g_1 \dots g_p} d(G_1) \dots d(G_p) \delta_{1}^{g_1}
\dots \delta_{p}^{g_p} \sum_{\l = 1}^{p} \frac{g + 1}{g_{\l} + 1}
\frac{h(G_{\l})}{d(G_{\l})} \delta_{\l} .
\end{split}
\end{equation*}
De m\^eme, en posant $\widetilde{V} = V_1 \times \dots \times V_p$
et $D := d_{\mathbb P}(V) / d_{\mathbb P}(\widetilde{V})$ (qui est
le nombre de composantes irréductibles de $V$ sur $\Kbar$), on a:
\begin{equation*}
\begin{split}
d_{{\mathbb P}_N}\!\left(\rho(V)\right) &= D d_{{\mathbb P}_N}\!\left(\rho(\widetilde{V})\right) \\
&= D \binom{\dim{V}}{\dim{V_1} \dots \dim{V_p}} d(V_1) \dots
d(V_p) \delta_{1}^{\dim{V_1}} \dots \delta_{p}^{\dim{V_p}}
\end{split}
\end{equation*}
et
\begin{equation*}
\begin{split}
h_{{\mathbb P}_N}\!\left(\rho(V)\right) &= D h_{{\mathbb P}_N}\!\left(\rho(\widetilde{V})\right) \\
&\!\!\!\!\!\!\!\!\!\!\!\!\!\!\!\!\!\!\!\!\!\!\!\!\!= D
\binom{\dim{V}}{\dim{V_1} \dots \dim{V_p}} d(V_1) \dots d(V_p)
\delta_{1}^{\dim{V_1}} \dots \delta_{p}^{\dim{V_p}} \!\sum_{\l =
1}^{p} \frac{\dim{V} + 1}{\dim{V_{\l}} + 1}
\frac{h(V_{\l})}{d(V_{\l})} \delta_{\l}
\end{split}
\end{equation*}
(en remarquant que $\dim{V} = \dim{\widetilde{V}}$). \\
Minorons maintenant la quantité $m_{W_r}(V)$ (remarquer que $W_r =
W_1$ est indépendant de $r$). On a par définition m\^eme (voir
\cite{Ph3}, page 1072):
$$m_{W_r}(V) := (g - \dim{V})! \!\!\!\!\!\!\!\!\!\!\!\!\!\!\!\!\!\!\!\max_{\begin{array}{cc}
\scriptstyle{\mathbf{y} \in V} \\
\scriptstyle{1 \leq i_1 < \dots < i_{g - \dim{V}} \leq g}
\end{array}}
\!\!\!\!\!\!\!\!\!\!\!\!\!\!\!\!\!\!\!{\rm{Vol}}\!\left({\mathbb{R}}_{+}^{g
- \dim{V}} \setminus C_{\scriptstyle{i_1 , \dots , i_{g -
\dim{V}}}}(W_r)\right)
$$
o\`u le maximum est pris pour les points $\mathbf{y}$ de $V$ et les $1 \leq i_1 < \dots < i_{g - \dim{V}} \leq g$ tels que les vecteurs de la base canonique de ${\mathbb N}^g$: $\varepsilon_{i_1} , \dots , \varepsilon_{i_{g - \dim{V}}}$ compl\`ete $T_{\mathbf{y}}V$ (l'espace tangent de $V$ en $\mathbf{y}$) dans ${\mathbb C}^g$ (en désignant par $(\varepsilon_1 , \dots , \varepsilon_g)$ la base canonique de ${\mathbb C}^g$) et o\`u $C_{\scriptstyle{i_1 , \dots , i_{g - \dim{V}}}}(W_r)$ désigne l'enveloppe convexe de la trace de l'escalier ${\mathbb N}^g \setminus W_r$ sur le $\mathbb R$-espace vectoriel de dimension $g - \dim{V}$: $\mathbb{R} \varepsilon_{i_1} + \dots + \mathbb{R} \varepsilon_{i_{g - \dim{V}}}$. \\
Il est clair qu'on a $m_{W_r}(V) = m_{W_r}(\widetilde{V})$ et
comme $\widetilde{V} = V_1 \times \dots \times V_p$, pour tout
$\mathbf{y} = ({\mathbf{y}}_1 , \dots , {\mathbf{y}}_p) \in
\widetilde{V}$ on a: $T_{\mathbf{y}}\widetilde{V} =
T_{{\mathbf{y}}_1}V_1 \times \dots \times T_{{\mathbf{y}}_p}V_p$.
Donc, Pour completer $T_{\mathbf{y}}\widetilde{V}$ dans ${\mathbb
C}^g$, on peut compléter chacune des bases des
$T_{{\mathbf{y}}_{\l}}V_{\l}$$(\l = 1 , \dots , p)$ dans ${\mathbb
C}^{g_{\l}}$. Ainsi les $i_1 , \dots , i_{g - \dim{V}}$
correspondant à cette complétion satisfont clairement:
$${\rm{Vol}}\!\left({\mathbb{R}}_{+}^{g - \dim{\widetilde{V}}} \setminus C_{\scriptstyle{i_1 , \dots , i_{g - \dim{V}}}}(W_r)\right) = \frac{\delta_{1}^{g_1 - \dim{V_1}} \dots \delta_{p}^{g_p - \dim{V_p}} \epsilon^{g - \dim{V}}}{(g - \dim{V})!} ,$$
d'o\`u:
$$m_{W_r}(V) = m_{W_r}(\widetilde{V}) \geq \epsilon^{g - \dim{V}} \delta_{1}^{g_1 - \dim{V_1}} \dots \delta_{p}^{g_p - \dim{V_p}} .$$
En substituant toutes ces estimations dans (\ref{1.39}) et en
majorant ensuite \\ $\binom{g}{g_1 \dots g_p} /
\binom{\dim{V}}{\dim{V_1} \dots \dim{V_p}}$ par $g^{g - \dim{V}}$,
la premi\`ere inégalité du théor\`eme \ref{a.19} en résulte. Afin
d'obtenir la deuxi\`eme inégalité du théor\`eme \ref{a.19}, on
substitue dans (\ref{1.40}) les estimations précédentes et on
utilise de plus les majorations suivantes: pour $\l = 1 , \dots ,
p$:
\begin{equation*}
\begin{split}
\log{(E_{\l})} + \log{2} &= \log{(F_{\l})} + \log{4} ~~~~\text{(car $E_{\l} = 2 F_{\l}$)} \\
&\leq 4(n_{\l} - g_{\l}) h(G_{\l}) + [2(n_{\l} - g_{\l} + 1) h({\mathbf{e}}_{\l}) + 4(n_{\l} - g_{\l})] d(G_{\l}) \\
&\quad+ (n_{\l} - g_{\l} + 1)(2 g_{\l} + 5)(\log{(d(G_{\l}))} + 1)
- 2(n_{\l} - g_{\l} + 1) h({\mathbf{e}}_{\l})
\end{split}
\end{equation*}
(d'apr\`es le lemme \ref{a.9} appliqué au groupe $G_{\l}$),
$$t_{\l}(W_0 + \dots + W_{r - 1}) = t_{\l}\left(W(\underline{\delta} ; (r - 1) \epsilon)\right) \leq (r - 1) \epsilon \delta_{\l} \leq (g - 1) \epsilon \delta_{\l} ,$$
$$\log{N} \leq \delta_1 \log{(n_1 + 1)} + \dots + \delta_p \log{(n_p + 1)} ,$$
\begin{equation*}
\begin{split}
\log{\!\left(d_{{\mathbb P}_N}\!\left(\rho(G)\right) {c'}^g + 1\right)} &\leq \log{\!\left(d_{{\mathbb P}_N}\!\left(\rho(G)\right) {c'}^g\right)} + 1 \\
&\!\!\!\!\!\!\!\!\!\!\!\!\!\!\!\!\!\leq \log{\!\left(\binom{g}{g_1 \dots g_p} d(G_1) \dots d(G_p) \delta_{1}^{g_1} \dots \delta_{p}^{g_p}\right)} + g \log{c'} + 1 \\
&\!\!\!\!\!\!\!\!\!\!\!\!\!\!\!\!\!\leq \log{\!\left(p^g d(G_1) \dots d(G_p) (n_1 + 1)^{\delta_1} \dots (n_p + 1)^{\delta_p}\right)} + g \log{c'} + 1 \\
&\!\!\!\!\!\!\!\!\!\!\!\!\!\!\!\!\!\leq g \log{p} + \log{(d(G_1))} + \dots + \log{(d(G_p))} \\
&~~~~~~~~~~+ \delta_1 \log{(n_1 + 1)} + \dots + \delta_p \log{(n_p
+ 1)} + g \log{c'} + 1 .
\end{split}
\end{equation*}
o\`u l'avant derni\`ere inégalité est obtenue gr\^ace à
l'hypoth\`ese $\delta_{\l} \geq g_{\l} + 1$ pour $\l = 1 , \dots ,
p$, qui entraine $\delta_{\l}^{g_{\l}} \leq (g_{\l} +
1)^{\delta_{\l}} \leq (n_{\l} + 1)^{\delta_{\l}}$. Apr\`es avoir
supposé (sans restreindre la généralité) que $\dim{V_{\l}} \leq
g_{\l} - 1$, $\forall \l = 1 , \dots , p$, on majore enfin: \\ $(g
+ 1) \binom{g}{g_1 \dots g_p} / (\dim{V} + 1)
\binom{\dim{V}}{\dim{V_1} \dots \dim{V_p}}$ par $(g + 1)^{g -
\dim{V}} \frac{1}{g_1 \dots g_p}$. Le reste consiste en des
majorations numériques triviales. La démonstration est achevée.
$~~~~\blacksquare$\vspace{1mm}
\subsubsection{{\bf Application aux puissances du groupe multiplicatif}}
En prenant dans notre th\'eor\`eme ci-dessus $G_1 , \dots , G_p$
des puissances du groupe multiplicatif ${\mathbb{G}}_m$:
${\mathbb{G}}_1 = {\mathbb{G}}_{m}^{n_1} , \dots , {\mathbb{G}}_p
= {\mathbb{G}}_{m}^{n_p}$, on obtient le corollaire suivant:
\begin{corollary}
Soit $K$ un corps de nombres, $0 < \epsilon \leq 1$ un r\'eel, $p
, n_1 , \dots , n_p$ des entiers strictement positifs et $n$ la
somme des entiers $n_1 , \dots , n_p$. Pour tout $1 \leq i \leq
p$, on plonge le groupe alg\'ebrique ${\mathbb{G}}_{m}^{n_i}$ dans
l'espace projectif ${\mathbb P}_{n_i}$ de la mani\`ere suivante:
$$
\begin{array}{ccc}
{\mathbb{G}}_{m}^{n_i} & \hookrightarrow & {\mathbb P}_{n_i} \\
(x_1 , \dots , x_{n_i}) & \mapsto & (1 : x_1 : \dots : x_{n_i})
\end{array} .
$$
Ces plongements induisent un plongement du groupe produit $G :=
{\mathbb{G}}_{m}^{n_1} \times \dots \times {\mathbb{G}}_{m}^{n_p}$
dans l'espace multiprojectif $\mathbb{P} := {\mathbb P}_{n_1}
\times \dots \times {\mathbb P}_{n_p}$ auquel est associ\'e un
anneau multigradu\'e que l'on d\'esigne par $K[\Xsoull]$ (avec
$\Xsoull := (\Xsoul_1 , \dots , \Xsoul_p)$ et $\Xsoul_i := (X_{i
0} , \dots ,
X_{i n_i})$ pour tout $1 \leq i \leq p$).\\
Soient aussi $P$ une forme non identiquement nulle de $K[\Xsoull]$
de multidegr\'e $\underline{\delta} = (\delta_1 , \dots ,
\delta_p) \in {\mathbb{N}}^p$ et $\x$ un point de $G(K)
\plongement {\mathbb P}(K)$. On suppose que l'on a pour tout $1
\leq i \leq p$: $\delta_i \geq n_i + 1$ et
$$\frac{\delta_i}{\delta_{i + 1}} ~>~ \left(\frac{p}{\epsilon}\right)^n .$$
Alors si $P$ s'annule en $\x$ avec la multiplicit\'e d\'efinie par
le dessous d'escalier $W(\underline{\delta} , n \epsilon)$ de
${\mathbb{N}}^n$, il existe une sous-vari\'et\'e propre $V$ de $G$
d\'efinie sur $K$, contenant $\x$, dont toutes les composantes
irr\'eductibles sur $\overline{K}$ sont des sous-vari\'et\'es
produit dans $\mathbb{P}$ et si l'on d\'esigne par $N$ le nombre
de ces derni\`eres composantes et par $V_1 \times \dots \times
V_p$ l'une de ces composantes, on a les in\'egalit\'es:
\begin{eqnarray*}
N d(V_1) \dots d(V_p) & ~\leq~ &
\left(\frac{n}{\epsilon}\right)^{n - \dim{V}} \\
N d(V_1) \dots d(V_p) \sum_{i = 1}^{p} \delta_i
\frac{h(V_i)}{d(V_i)} & ~\leq~ & \left(\frac{n +
1}{\epsilon}\right)^{\!\!n - \dim{V}} \!\!\left\{h(P) + \left(n
\epsilon + 2\right) \sum_{i = 1}^{p} (2 n_i + 5) \delta_i \right.\\
& ~ &\left.\phantom{\sum_{i = 1}^{p}}\!\!\!\!\!\!\!\!+ 3 (n + 1)
(\log{p} + 2)\right\} .
\end{eqnarray*}
De plus, $P$ s'annule sur $V$ avec la multiplicité définie par le
dessous d'escalier $W(\underline{\delta} , \epsilon)$.
\end{corollary}
{\bf Démonstration.---} Appliquons le théorème \ref{a.19} aux
groupes algébriques $G_i := {\mathbb G}_{m}^{n_i}$ plongés dans
les espaces projectifs ${\mathbb P}_{n_i}$ $(1 \leq i \leq p)$
comme expliqué dans l'énoncé du corollaire. Ces groupes
algébriques $G_1 , \dots , G_p$ sont de dimensions respectives
$g_1 := n_1 , \dots , g_p := n_p$; leur produit $G := G_1 \times
\dots \times G_p$ est donc de dimension $g := n_1 + \dots + n_p =
n$. De plus, le degré de chaque groupe $G_i \plongement {\mathbb
P}_{n_i}$ est bien égal à $1$ et sa hauteur (comme $G_i$ est dense
dans ${\mathbb P}_{n_i}$ au sens de la topologie de Zariski) est
égale à:
\begin{eqnarray*}
h(G_i) & = & h({\mathbb P}_{n_i}) \\
& = & \sum_{\ell = 1}^{n_i} \sum_{k = 1}^{\ell} \frac{1}{2 k} ~~~~
(\text{voir \cite{Ph4}3, page 346}) \\
& = & \sum_{k = 1}^{n_i} \sum_{\ell = k}^{n_i} \frac{1}{2 k} \\
& = & \sum_{k = 1}^{n_i} \frac{n_i + 1 - k}{2 k} \\
& = & \frac{n_i + 1}{2} \sum_{k = 1}^{n_i} \frac{1}{k} -
\frac{n_i}{2} \\
& = & \frac{n_i + 1}{2} \sum_{k = 2}^{n_i + 1} \frac{1}{k} \\
& \leq & \frac{1}{2} (n_i + 1) \log(n_i + 1) .
\end{eqnarray*}
Par ailleurs, l'addition dans chaque groupe $G_i \plongement
{\mathbb P}_{n_i}$ est représentée par la famille de formes
$\Asoul_i$ de $K[\Xsoul_i , \Ysoul_i]$ (avec $\Xsoul_i := (X_{i 0}
, \dots , X_{i n_i})$ et $\Ysoul_i := (Y_{i 0} , \dots , Y_{i
n_i})$) définie par:
$$\Asoul_i(\Xsoul_i , \Ysoul_i) ~=~ (X_{i 0} Y_{i 0} , \dots , X_{i n_i} Y_{i n_i}) .$$
La famille $\Asoul_i$ est donc constituée de formes de bidegré
$(\c , \c') = (1 , 1)$ et elle est de hauteur de Gauss-Weil
nulle.\\
Enfin, l'élément neutre ${\mathbf{e}}_i$ de $G_i$ est représenté
dans ${\mathbb P}_{n_i}$ par le système de coordonnées projectives
$(1 : \dots : 1)$, donc $\mathbf{e_i}$ est de hauteur de
Gauss-Weil nulle aussi.\\
Le reste ne sont que des majorations numériques triviales.
$~~~~\blacksquare$

\end{document}